\documentclass[11pt,a4paper]{amsart}
\usepackage{amsmath}
\usepackage{amscd}
\usepackage{amssymb}

\numberwithin{equation}{section}

\theoremstyle{plain} 
\newtheorem{thm}{Theorem}[section]
\newtheorem{cor}[thm]{Corollary}
\newtheorem{lem}[thm]{Lemma}
\newtheorem{prop}[thm]{Proposition}

\newtheorem{defn}[thm]{Definition}

\newtheorem{rem}[thm]{Remark}

\theoremstyle{remark}

\newcommand{\thmref}[1]{Theorem~\ref{#1}}
\newcommand{\propref}[1]{Proposition~\ref{#1}}

\newcommand{\lemref}[1]{Lemma~\ref{#1}}
\newcommand{\corref}[1]{Corollary~\ref{#1}}

\newsavebox{\SmallMathBox}

\def\pdo{\psi{\rm do}}

\def\Ci{C^\infty}

\def\dd{\partial}

\def\Di{D\kern -.65em /}
\def\Dii{D\kern -.45em /}
\def\di{{\dd}\kern -.55em /}
\def\dii{{\dd}\kern -.40em /}

\def\na{\nabla}
\def\noi{\noindent}

\def\to{\rightarrow}
\def\too{\longrightarrow}

\def\mtoo{\longmapsto}

\def\re{{\rm Re}}

\def\NN{{\bf N}}

\def\pp{{\bf p}}

\def\Aa{{\mathcal A}}

\def\Cc{{\mathcal C}}
\def\Dd{{\mathcal D}}
\def\Ee{{\mathcal E}}
\def\Ff{{\mathcal F}}
\def\Gg{{\mathcal G}}

\def\Kk{{\mathcal K}}

\def\Mm{{\mathcal M}}
\def\Nn{{\mathcal N}}

\def\Pp{{\mathcal P}}

\def\Rr{{\mathcal R}}
\def\Sss{{\mathcal S}}
\def\Tt{{\mathcal T}}
\def\Uu{{\mathcal U}}
\def\Vv{{\mathcal V}}
\def\Ww{{\mathcal W}}

\def\={\cong}
\def\>{\supset}
\def\<{\subset}
\def\ii{^{-1}}
\def\si{^{-s}}

\def\pp{^{\perp}}
\def\12{\frac{1}{2}}
\def\2{\Dd}
\def\3{\Nn}
\def\4{\Rr}
\def\6{\cup}
\def\8{\otimes}
\def\0{^{\circ}}
\def\){\hfill{\ \qed}\enddemo}

\def\a{\alpha}
\def\A{\Aa}
\def\b{\beta}

\def\C{\mathbb{C}}
\def\Cn{\mathbb{C}^{n}}

\def\d{\delta}

\def\e{\varepsilon}

\def\g{\gamma}
\def\G{\Gamma}

\def\la{\lambda}

\def\N{\NN}
\def\o{\infty}
\def\p{\pi}

\def\Si{\Sigma}

\def\ch{\mbox{\rm ch}}

\def\dom{\mbox{\rm dom\,}}

\def\End{\mbox{\rm End}}

\def\Gr2n{\mbox{${\rm Gr}(\Cn\oplus\Cn)$}}

\def\Grk2n{\mbox{${\rm Gr}_{k}(\Cn\oplus\Cn)$}}
\def\Grk{\mbox{${\rm Gr}_{k}$}}
\def\Gr{\mbox{${\rm Gr}$}}

\def\Hom{\mbox{\rm Hom}}

\def\index{\mbox{\rm index\,}}
\def\Ind{\mbox{\rm Ind\,}}

\def\LIM{\mbox{\rm LIM}}

\def\Ker{\mbox{\rm Ker}}

\def\ind{\mbox{\rm ind\,}}

\def\ran{\mbox{\rm ran\,}}

\def\Si{S\kern -.65em /}

\def\str{\mbox{\rm Str}}

\def\tr{\mbox{\rm tr\,}}
\def\Tr{\mbox{\rm Tr\,}}

\def\Af{\mathbb{A}}
\def\Bf{\mathbb{B}}
\def\Cf{\mathbb{C}}
\def\Df{\mathbb{D}}
\def\Ef{\mathbb{E}}
\def\Fff{\mathbb{F}}

\def\Hf{\mathbb{H}}
\def\If{\mathbb{I}}

\def\Kf{\mathbb{K}}

\def\Nf{\mathbb{N}}

\def\Pf{\mathbb{P}}

\def\Rf{\mathbb{R}}

\def\Tf{\mathbb{T}}

\def\Zf{\mathbb{Z}}

\def\wPf{{\widehat \Pf}}

\def\weta{{\widehat \eta}}

\def\t{\tilde}

\def\ord{\mbox{\rm ord}}

\def\Pker{\mbox{$P_{\Kf {\rm er}}$}}

\def\Ps{\textsf{P}}
\def\Ds{\textsf{D}}
\def\Rs{\textsf{R}}
\def\Gs{\textsf{G}}
\def\Hs{\textsf{H}}
\def\Ks{\textsf{K}}
\def\Ts{\textsf{T}}
\def\Ss{\textsf{S}}
\def\Rt{\textsf{R}^{\pi_{*}(\Ef|\wPf)}}
\def\DP{\textsf{D}_{\Pf}}
\def\DP1{\textsf{D}_{\Pf_1}}
\def\DP2{\textsf{D}_{\Pf_2}}

\def\nE{\na^{\pi_{*}(\Ef)}}

\def\wPf{\widehat{\Pf}}

\begin{document}

\title[]{Eta Forms and the Chern Character\vskip 5mm {\Small {\rm
S. Scott}}}

\address{King's College, University of London.}

\maketitle

\vskip 5mm \noi {\small ABSTRACT: The semi-topological nature of
the eta-invariant of a self-adjoint elliptic differential operator
derives from a relative identification with a Chern character.
This remarkable semi-locality property of the eta-invariant can be
seen in spectral flow formulae and many other applications
\cite{AtPaSi2,AtPaSi3,BiCh89,DaZh97,Lo84}. In this paper we prove
two geometric index theorems for a family of first-order elliptic
operators over a manifold with boundary by computing eta form
representatives for the Chern character classes of the index
bundle. The eta forms occur as  relative and regularized traces on
infinite-dimensional vector bundles realized as the limiting
values of superconnection character forms. The formulas are
non-local and general, they do not require spin structures,
compatibility with Clifford actions, or dimensional restrictions.}

\vskip 15mm

\section{Introduction}


Let $M \too B$ be a smooth Riemannian fibration with fibre
diffeomorphic to a closed even-dimensional spin manifold $X$ and
let $\Ef=\Ef^+\oplus\Ef^-$ be a graded vector bundle of Clifford
modules over $M$. Let $\Ds = \{D_z \ | \ z\in B\}$ be a family of
compatible Dirac operators acting on the space of $\Ci$ sections
$\Gamma(M,\Ef)$. The kernels $\Ker(D_z) =
\Ker(D^{+}_z)\oplus\Ker(D^{-}_z)$ are $\Zf_2$-graded by the
kernels of the chiral Dirac operators $D^{+}_z$, and if assumed to
vary smoothly with $z$ form a finite rank superbundle $\Ker(\Ds) =
\Ker(\Ds^+)\oplus\Ker(\Ds^-)$ on $B$. The eta-form $\eta\in
\Aa^{\rm odd}(B) = \sum_i \G(B,\wedge^{2i+1}T^* B)$ of the Bismut
superconnection $\Af_t$ is the odd-degree differential form
introduced by Bismut \cite{Bi86}, Bismut-Cheeger \cite{BiCh89},
Berline-Verne \cite{BeVe90}
\begin{equation*}
    \eta = \frac{1}{\sqrt{\pi}}\int^{\o}_0
    \str(\dot{\Af_t}e^{-\Af^{2}_t}) \ dt \
\end{equation*}
as a canonical transgression form for the local families index
formula \cite{Bi86,BeGeVe92}
\begin{equation}\label{e:localindexform}
\ch(\Af^{{\small {\rm Ker}(\Ds)}}) = \int_{M/B}
\widehat{A}(M/B)\ch^{\prime}(\Ef) - d\eta \ .
\end{equation}
Here $\ch(\Af^{{\small {\rm Ker}(\Ds)}}) = \str(e^{-({\Af^{\small
{\rm Ker}(\Ds)}})^2}) \in \Aa(B)$ is a Chern character form for
the index bundle $$\Ind(\Ds) = [{\rm Ker}(\Ds^+)] - [{\rm
Ker}(\Ds^-)] \in K(B) \ ,$$ and hence \eqref{e:localindexform}
implies the cohomological Atiyah-Singer families index theorem in
${\rm H}^{\bullet}(B)$
\begin{equation}\label{e:familiesindex}
\ch(\Ind(\Ds) ) = \int_{M/B} \widehat{A}(M/B)\ch^{\prime}(\Ef) \ .
\end{equation}
Thus, though not local, $\eta$ defines via
\eqref{e:localindexform}  a secondary topological invariant in the
form of a generalized Chern-Simons form for superconnections.

Extensions to compatible Dirac operators associated to a fibration
of manifolds with boundary $\pi:M \stackrel{X}{\too} B$ were
achieved by Bismut-Cheeger \cite{BiCh89} and Melrose-Piazza
\cite{MePi97}. For $\Ds$ to have a well-defined index non-local
boundary conditions must be imposed. The induced boundary
fibration $\pi^N : N = \dd M \too B$ of closed Riemannian spin
manifolds with fibre $Y_z = \dd X_z$ defines a family of
self-adjoint compatible Dirac operators $\dd = \{\partial^z \ | \
z\in B\}$ acting on the infinite-dimensional bundle
$\pi_*^N(\Ef^0)$ over $B$ whose fibre at $z\in B$ is the space of
sections over $Y_z$. $\dd$ occurs in the degree zero component of
the induced Bismut superconnection $\Bf_t = (\Af_t)_{|N}$ on the
boundary. For simplicity we assume the operators $\dd^z$ are
invertible. Let $\Pi_> =\{\Pi_{>}^z \ | \ z\in B\}$ be the smooth
family of boundary projections onto the direct sum of the positive
eigenspaces of $\dd^z$. Then by restricting the domain of $\Ds$ to
those sections whose boundary values lie in the kernel of
$\Pi_{>}$ we obtain a smooth family of Fredholm operators
$\Ds_{\Pi_{>}}$ for which one can aim to construct analogues of
\eqref{e:localindexform}, \eqref{e:familiesindex}. The index form
$\int_{M/B} \widehat{A}(M/B)\ch^{\prime}(\Ef)$ is no longer
closed, rather a correcting non-local boundary Eta form term is
required
\begin{equation}\label{e:boundaryetaform}
    \weta_{\Pi_{\geq}} = \frac{1}{\sqrt{\pi}}\int^{\o}_0
    \Tr_{{\rm even}}(\dot{\Bf_t}e^{-\Bf^{2}_t}) \ dt \ ,
\end{equation}
where $\Tr_{{\rm even}}$ is the trace over forms of even degree.
The cohomological APS families index theorem states that in
$H^{\bullet}(B)$
\begin{equation}\label{e:familiesindex2}
\ch(\Ind(\Ds_{\Pi_{>}}) ) = \int_{M/B}
\widehat{A}(M/B)\ch^{\prime}(\Ef) - \frac{\weta_{\Pi_{>}}}{2}\ .
\end{equation}
This was proved by Bismut-Cheeger \cite{BiCh90} by an adiabatic
limit argument identifying the index bundle with that for a family
of Dirac operators over a closed manifold with singular metric
formed by attaching a cone to $X$. A conceptually simpler proof
and generalization to spectral sections was subsequently given by
Melrose-Piazza \cite{MePi97} using the $b$-calculus --- a spectral
section is a smooth family of pseudodifferential operator ($\pdo$)
boundary projections $\Pp = \{P_z \ | \ z\in B\}$ differing from
$\Pi_{>}$ by a family of finite-rank operators. Melrose-Piazza
based their proof on the observation that $\Pp$ occurs as the
spectral projection $\Pi^{\dd_1}_{>}$ of a family $\dd_1$ of
boundary $\pdo$s differing from $\dd$ by a family of finite-rank
$\pdo$s of order 0, and thus generalized \eqref{e:familiesindex2}
to the case where $\Ds_{\Pp}$ replaces $\Ds_{\Pi_>}$ and
$\weta_{\Pi_>}$ is replaced by an analogously defined eta form
$\weta_{\Pp}$.

The corresponding transgression formula for
\eqref{e:familiesindex2} for any spectral section $\Pp$ implies
that as elements of $\Aa(B)$
\begin{equation*}
 d \weta_{\Pp} = 2\int_{\dd M/B} \widehat{A}(\dd M
 /B)\ch^{\prime}(\Ef) \ .
\end{equation*}
The relative eta form $\weta_{\Pp_1} - \weta_{\Pp_2}$ is thus
closed, and defines a cohomology class in $H^{\bullet}(B)$ equal
to the relative Chern character
\begin{equation}\label{e:MPrelchch}
\ch(\Ind(\Ds_{\Pp_1}))  - \ch(\Ind(\Ds_{\Pp_2})) =
\ch(\Ind(\Pp_2,\Pp_1) ) \ ,
\end{equation}
where $(\Pp_2,\Pp_1)$ is a canonical family of boundary Fredholm
operators defined by the spectral sections. The equality  in
\eqref{e:MPrelchch} follows from K-theory arguments
\cite{DaZh97,MePi97,Sc03}.

\subsection{Theorem (I): Relative Boundary Eta forms}

To explain our constructions we begin with the case where $B$ is
just a single point. Then \eqref{e:familiesindex2} reduces to the
classical APS formula for the index of the single operator
$D_{\Pi^{\dd}_{>}}$
\begin{equation*}
    \ind(D_{\Pi^{\dd}_{>}}) = \int_X \widehat{A}(X)\ch(E) +
    \frac{\weta(\dd)}{2} \ .
\end{equation*}
Now $\Bf = t^{1/2}\dd$ and \eqref{e:boundaryetaform} becomes the
usual single operator eta-invariant
\begin{eqnarray}\label{e:etainvariant}
    \weta(\dd) & = & \frac{1}{\sqrt{\pi}}\int^{\o}_0
    t^{-1/2}\Tr(\dd e^{-t\dd^2}) \ dt \nonumber \\
        & = & \Tr(\dd|\dd|^{-s-1})|_{s=0}^{{\rm mer}} \ ,
\end{eqnarray}
the superscript indicating the meromorphically continued trace
evaluated at $s=0$. We may rewrite \eqref{e:etainvariant} as the
regularized trace
\begin{equation}\label{e:etareg}
 \weta(\dd) = \Tr((\Pi^{\dd}_{>} - \Pi^{\dd}_{<})|\dd|^{-s})|_{s=0}^{{\rm
 mer}}
\end{equation}
of the involution $\Pi^{\dd}_{>} - \Pi^{\dd}_{<}$ defined by the
order 0 $\pdo$ projections
$$\Pi^{\dd}_{>} = \frac{1}{2}(I + \dd|\dd|\ii) \ , \hskip 10mm
\Pi^{\dd}_{<}  = \frac{1}{2}(I - \dd|\dd|\ii) =
(\Pi^{\dd}_{>})\pp$$ onto the positive and negative spectral
subspaces of $\dd$.

Consider the infinite Grassmannian $\Gr_{\o}(\dd)$ parameterizing
$\pdo$ projections $P$ such that $P - \Pi_{>}$ is an element of
the algebra $\Psi^{-\o}(E^0)$ of smoothing operators on
$\G(Y,E^0)$. If $P_1, P_2\in \Gr_{\o}(\dd)$, then $P_1 - P_2 \in
\Psi^{-\o}(E^0)$ is represented by a smooth kernel $k(y_1, y_2)$
with trace $\Tr(P_1 - P_2) = \int_Y \tr (k(y,y)) \ dy$, and so the
relative variant of \eqref{e:etareg} exists without regularization
\begin{equation}\label{e:releta}
    \weta(P_1,P_2) = \Tr\left((P_1 - P^{\perp}_1) - (P_2 -
    P^{\perp}_2)\right) \ ,
\end{equation}
where $P_i^{\perp} = I - P_i$.

\begin{lem}
For spectral sections $\Pi_{>}^{\dd_1},\Pi_{>}^{\dd_2}$ the
additivity formula holds
\begin{equation}\label{e:releta3}
\weta(\Pi_{>}^{\dd_1},\Pi_{>}^{\dd_2}) = \weta(\dd_1) -
\weta(\dd_2) \ .
\end{equation}
\end{lem}
\vskip 2mm

\noi One can prove \eqref{e:releta3} directly from elementary
properties of the Wodzicki residue trace. Alternatively, rewriting
\eqref{e:releta} as
\begin{equation}\label{e:releta4}
    \frac{\weta(P_1,P_2)}{2} = \Tr(P_1 - P_2) = \ind(P_1,P_2) \ ,
\end{equation}
where $(P_1,P_2) := P_2\circ P_1 : W_1 \to W_2$ and $W_i =
\ran(P_i)$, it is equivalent to the relative index formula
\begin{equation*}
\frac{\weta(P_2,P_1)}{2} =  \ind(D_{P_1}) - \ind(D_{P_2})  \ ,
\end{equation*}
which coincides with \eqref{e:MPrelchch} for the case $B= pt$.
Notice that, unlike $\weta(\dd_i)$, the relative eta-invariant is
a homotopy invariant.

To generalize this to smooth families of APS-type boundary
problems the $\pdo$ projection $P_i$ must be replaced by a
spectral section --- or, more generally, by a Grassmann section
$\Pp_i = \{P_{i,z} \ | \ z\in B\}$, defined to be a smooth section
of the fibration with fibre $\Gr_{\o}(\dd^z)$ at $z\in B$. $\Pp_i$
defines a smooth infinite-dimensional subbundle $\Ww_i$ of
$\pi_{*}^N(\Ef^0)$ with a canonical connection with curvature
$(\Pp_i\cdot\na^{\pi_{*}^N(\Ef^0)}\cdot\Pp_i)^2\in\Aa^2(B,\End(\Ww_i))$,
which by extending by zero we can consider as a 2-form
endomorphism of $\pi_{*}^N(\Ef^0)$, where $\na^{\pi_{*}^N(\Ef^0)}$
is the Bismut connection on $\pi_{*}^N(\Ef^0)$. Generalizing
\eqref{e:releta4} we define the relative Eta-form by
\begin{equation*}
\eta(\Pp_1,\Pp_2) =
\Tr\left(e^{-(\Pp_1\cdot\na^{\pi_{*}^N(\Ef^0)}\cdot\Pp_1)^2} -
e^{-(\Pp_2\cdot\na^{\pi_{*}^N(\Ef^0)}\cdot\Pp_2)^2}\right) \in
\Aa^{2\bullet}(B) \ .
\end{equation*}
Here the individual exponential operators
$e^{-(\Pp_i\cdot\na^{\pi_{*}^N(\Ef^0)}\cdot\Pp_i)^2}$ are families
of unbounded $\pdo$s of positive order and so not of trace class.
The relative exponential operator, however, is a smooth family of
smoothing operators with coefficients in $\Aa^{2\bullet}(B)$ and
so has a well-defined fibrewise trace.

\vskip 2mm

\noi {\bf Theorem (I)} {\it Let $\Ds$ be a smooth family of
first-order elliptic operators. In $H^{2\bullet}(B)$
\begin{equation}\label{e:tI1}
\ch(\Ind(\Ds_{\Pp}) )  = \eta(\Pp(\Ds),\Pp) \ .
\end{equation}
Equivalently,
\begin{equation}\label{e:tI2}
\ch(\Ind(\Ds_{\Pp_1})) - \ch(\Ind(\Ds_{\Pp_2}))  =
\eta(\Pp_2,\Pp_1) \ .
\end{equation}
One has in $\Aa(B)$
\begin{equation*}
\eta(\Pp_2,\Pp_1)_{[2k]} =
\frac{(-1)}{k!}^k\Tr\left((\Pp_1\cdot\na^{\pi_{*}^N(\Ef^0)}\cdot\Pp_1)^{2k}
- (\Pp_2\cdot\na^{\pi_{*}^N(\Ef^0)}\cdot\Pp_2)^{2k}\right) \ ,
\end{equation*}
where $\omega_{[j]}\in\Aa^j(B)$ is the $j$-form component of
$\omega\in\Aa(B)$.} \vskip 2mm

Here, the Calderon section  $\Pp(\Ds)$ is a Grassmann section
canonically associated to the infinite-dimensional subbundle
$\Ker(\Ds)$ of $\pi_*(\Ef)$.

For $k=0$, $(\Pp_i\cdot\na^{\pi_{*}^N(\Ef^0)}\cdot\Pp_i)^0 =
\Pp_i$, and hence the degree zero component of \eqref{e:tI2}
coincides with the classical pointwise index identity
\eqref{e:releta4}.

For related relative eta form formulae in the context of gerbes we
refer to the work of Lott \cite{Lo02}.

Theorem (I) and the Melrose-Piazza Theorem \cite{MePi97}, stated
in \eqref{e:familiesindex2}, imply:
\begin{cor}
If $\Ds$ is a family of compatible Dirac operators and $\Pp$ is a
spectral section, then in $H^{2\bullet}(B)$
\begin{equation}\label{e:IndMP1}
\eta(\Pp ,\Pp(\Ds) ) =  \int_{M/B}
\widehat{A}(M/B)\ch^{\prime}(\Ef) - \frac{\weta_{\Pp}}{2} \ .
\end{equation}
Generalizing \eqref{e:releta3}, one has
\begin{equation}\label{e:IndMP2}
\eta(\Pp_2 ,\Pp_1) =  \frac{\weta_{\Pp_1}}{2} -
\frac{\weta_{\Pp_2}}{2} \ ,
\end{equation}
where $\weta_{\Pp_i} = \weta(\dd_i)$ are the eta-forms defined
similarly to \eqref{e:boundaryetaform}.
\end{cor}

\vskip 2mm The index formulas \eqref{e:familiesindex2},
\eqref{e:IndMP1} refer specifically to a family of compatible
Dirac operators, though Dai-Zhang \cite{DaZh98} generalize to
families of generalized Dirac operators with spectral section
boundary conditions. It is worth pointing out that Theorem (I),
and Theorem (II), below, hold for any smooth family of first-order
elliptic differential operators and any Grassmann section, and, in
fact, may be extended to significantly more general classes of
boundary condition and to families of higher order elliptic
operators. [Such generality is needed, for example, in local
anomaly computations for supersymmetric branes.]

\subsection{Theorem (II): Relative Interior Eta forms}
Theorem (I) may, in the general sense explained in Section (1.3),
be seen as a Quillen-Chern-Weil model for Theorem(II), which is
the principal result of this paper.

Each of \cite{BiCh90}, \cite{MePi97} and Theorem(I), transform the
family of boundary problems to a K-theory equivalent family of
elliptic operators over a closed, or complete, manifold to which
the methods of \eqref{e:familiesindex} can be adapted.
Theorem(II), in contrast, computes the heat trace of a
superconnection constructed directly on the infinite-dimensional
bundle $\pi_*(\Ef|\Pp)$ of domains of the family of boundary
problems. Like \cite{BiCh90}, \cite{MePi97}, this relies on a
generalized $\pdo$ boundary calculus, here the singular Green
operator (sgo) calculus developed by Grubb-Seeley \cite{GrSe95},
\cite{Gr00}.

There are a number of technical difficulties inherent in a direct
assault on the families index, starting with the bundle structure
of $\pi_*(\Ef|\Pp)$. We show that this is inherited in a natural
way from the structure of $\pi_*(\Ef)$ and the boundary bundle
$\Ww$, realized explicitly by a smooth family  of sgo projections
$\Ps$ on $\pi_*(\Ef)$. For Grassmann sections $\Pp_i$, $i=1,2$,
the corresponding vertical sgo projections $\Ps_i $ induce
superconnections on $\pi_*(\Ef|\Pp_i)$ with curvature forms
$\Rs_i\in\Aa^2(B,\End(\pi_*(\Ef|\Pp_i)))$ with restricted
differential operator coefficients. The relative curvature form
$\Rs_1 - \Rs_2$ is then trace-class in $\Aa(B,\End(\pi_*(\Ef)))$,
and there is a generalized  relative eta-form defined by the
supertrace
\begin{equation*}
\eta^{[M]}(\Ps_1,\Ps_2) = \str\left(e^{-\Rs_1} - e^{-\Rs_2}\right)
\in \Aa^{2\bullet}(B) \ .
\end{equation*}
This is not a closed form, in general, rather the Chern character
is modified by zeta-trace correction terms. These arise as
follows. First, consider over a closed manifold a classical
pseudodifferential operator $F$, and let $\Delta$ be an elliptic
$\pdo$ of order $\ord(\Delta)>0$ of even parity, acting on the
space of sections of a $\Zf_2$-graded vector bundle. Then from
Grubb-Seeley \cite{GrSe95}, for $m >> 0 $ there is a resolvent
trace expansion as $\la\to\o$
\begin{equation}\label{e:resexp}
\str(F(\Delta - \la)^{-m}) \sim \sum_{j=-n}^{\o}c_j
(-\la)^{\frac{{\tiny \ord}(F)-j}{{\tiny \ord}(\Delta)}-m} +
\sum_{l=0}^{\o}(c^{\prime}_l \log(-\la) + c_l^{\prime\prime})
(-\la)^{-l-m} \ .
\end{equation}
It follows that the spectral zeta-function $\str(F\Delta\si)$,
defined for $\re(s) > m+\ord(F)/\ord(\Delta)$, extends
meromorphically to all of $\C$ with the pole structure
\begin{equation}\label{e:Deltapoles}
\G(s)\,\str(F\Delta\si) \sim \sum_{j=-n}^{\o}\frac{\t{c}_j}{ s+
\frac{j-{\tiny \ord}(F)}{{\tiny \ord}(\Delta)}} -
\frac{\str(F\Pi_0(\Delta))}{s} +
\sum_{l=0}^{\o}\left(\frac{\t{c}^{\prime}_l}{(s+ l)^2} +
\frac{\t{c}^{\prime\prime}_l}{s+ l}\right) \ ,
\end{equation}
where $\Pi_0(\Delta)$ is the orthogonal projection onto
$\Ker(\Delta)$. The coefficients differ from those in
\eqref{e:resexp} by universal constants. In particular, the
coefficient $\t{c}^{\prime}_0$ determines the Guillemin-Wodzicki
residue supertrace ${\rm res}(F)$ of $F$; one has
\begin{equation*}
{\rm res}(F) = \ord(\Delta) \cdot {\rm
Res}_{s=0}(\str(F\Delta\si)|^{{\rm mer}}) = \ord(\Delta)\cdot
\t{c}^{\prime}_0 \ .
\end{equation*}
This is independent of the choice of regularizing operator
$\Delta$ and defines the projectively unique supertrace functional
on the algebra of classical $\pdo$s \cite{Wo84}, \cite{Gu85}.

The  {\em pseudo}-trace $\tau_{\Delta}(F)$, on the other hand, is
defined to be the sum of constant term coefficients
$$\t{c}_{{\tiny \ord}(F)} + \t{c}_0^{\prime\prime} = c_{{\tiny
\ord}(F)} + c_0^{\prime\prime} \  $$ in the Laurent expansion of
$\str(F\Delta\si)|^{{\rm mer}}$ around zero. Thus
\begin{equation}\label{e:regtrace}
\tau_{\Delta}(F) = \LIM_{s\to 0} \, \str\left(F(\Delta\si +
\Pi_0(\Delta))\right) \ .
\end{equation}
This is not independent of the choice of regularizing operator
$\Delta$, and is not a supertrace functional on the space of
classical $\pdo$s. An example of such a pseudo-trace is the eta
invariant in \eqref{e:etareg}
$$\weta(\dd) = \tau_{|\dd|}(\Pi^{\dd}_{>} -
\Pi^{\dd}_{<}) \ , $$ with the trivial grading. The functional
$\tau_{\Delta}$ does, however, restrict to a supertrace
independent of the regularizing operator on certain subclasses of
operators, and on which the Wodzicki residue vanishes. In
particular it coincides with the $L^2$-supertrace on operators of
order less than $-\dim(X)$, and restricts to a supertrace on
non-integer order $\pdo$s \cite{KoVi94}. Otherwise the dependence
on the choice of regularizing operator
\begin{equation*}
\tau_{\Delta_1,\Delta_2}(F) := \LIM_{s\to 0}
\,\str\left(F(\Delta_1\si -\Delta_2\si + \Pi_0(\Delta_1) -
\Pi_0(\Delta_2))\right)
\end{equation*}
is determined by the Wodzicki residue trace of the operator
$F(\log \Delta_1 - \log \Delta_2)$ \cite{KoVi94}.

On closed manifolds this is essentially well understood. On
manifolds with boundary less is known, but the expansions
\eqref{e:resexp}, \eqref{e:Deltapoles} are known to hold when $F$
is a differential operator and the regularizing operator is an
APS-type Dirac-Laplacian \cite{Gr02}.\vskip 1mm

We have: \vskip 3mm

 \noi {\bf Theorem (II)} {\it In $H^*(B)$
\begin{equation}\label{e:tII2}
\ch(\Ind(\Ds_{\Pp_1})) - \ch(\Ind(\Ds_{\Pp_2})) =
 \eta^{[M]}(\Ps_1,\Ps_2) + \sum_{k\geq 0}\frac{k+1}{k!} \
 \tau_{\Delta_1,\Delta_2}(\Rs^k) \ .
 \end{equation}
where $\Rs$ is the curvature operator associated to any Grassmann
section, and $\Delta_i = \Ds^2_{\Pp_i}$ is the Dirac-Laplacian.
Equivalently,
\begin{equation}\label{e:tII1}
\ch(\Ind(\Ds_{\Pp})) = \eta^{[M]}(\Ps,\Ps(\Ds))  + \sum_{k\geq
0}\frac{k+1}{k!} \ \tau_{\Delta_{\Pp},\Delta_{\Pp(\Ds)}}(\Rs^k) \
.
\end{equation}
}

\vskip 2mm

The regularized pseudo-trace $\tau_{\Delta_1,\Delta_2}(\Rs^k)$
means pointwise on the operator coefficient of $\Rs^k$. The
construction of the connection forms means that the order of
$\Ps_1, \Ps_2$ is reversed from $\Pp_1, \Pp_2$ in Theorem(I). The
statement that the curvature operator $\Rs$ for any Grassmann
section can be used in the relative regularized trace term is
simply that $\tau_{\Delta_1,\Delta_2}(S)$ vanishes if $S$ is
smooth family of generalized smoothing operators, in the sense
defined in Section(4.1).

Applying the Melrose-Piazza Theorem \eqref{e:familiesindex2} we
obtain:

\begin{cor} If $\Ds$ is a family of compatible Dirac operators
and $\Pp$ is a spectral section, then in $H^{2\bullet}(B)$
\begin{equation}\label{e:II=MP}
\eta^{[M]}(\Ps,\Ps(\Ds)) + \sum_{k\geq 0}\frac{k+1}{k!} \
\tau_{\Delta_{\Pp},\Delta_{\Pp(\Ds)}}(\Rs^k)  = \int_{M/B}
\widehat{A}(M/B)\ch^{\prime}(\Ef)  - \frac{\weta_{\Pp}}{2} \ .
\end{equation}
\end{cor}

\noi Taking differences of \eqref{e:II=MP}
one is left with the interior analogue of
\eqref{e:IndMP2}:

\begin{cor} For spectral sections $\Pp_1, \Pp_2$ one has in
$H^*(B)$
\begin{equation*}
\eta^{[M]}(\Ps_1,\Ps_2) + \sum_{k\geq 0} \frac{k+1}{k!} \
\tau_{\Delta_1,\Delta_2}(\Rs^k) = \frac{\weta_{\Pp_2}}{2} -
\frac{\weta_{\Pp_1}}{2} \ .
\end{equation*}
\end{cor}

\vskip 4mm

\subsection{Relation to Quillen-Chern-Weil Theory} A consequence of
using the computational device of relative supertraces is that in
the simplest non-trivial case the constructions in Thm(I) reduce
to a superconnection formulation of Chern-Weil theory on the
classifying space ${\rm BGl}(\o)$. This is interesting because
{\it a priori} there is no total Chern character form as the trace
of an operator $e^{-\na^2}$ on the universal bundle $\Ee_1$ over
${\rm BGl}(\o)$ in the semi-ring of vector bundles ${\rm
Vect}(\o)$ --- however, there is a Chern character form at the
level of K-theory, defined on virtual bundles.

By way of example and to explain some of the methodology of our
approach we pause for a moment to outline this further.

\subsubsection{Construction of Heat Operators}
The heat operators we shall use are defined via holomorphic
functional calculus, thus placing primacy on the resolvent
operator. This works in Quillen-Chern-Weil theory as follows.
Consider a complex super vector bundle $\Ee= \Ee^+\oplus \Ee^-$
over a finite-dimensional manifold $B$, and let $L$ be an odd
Hermitian endomorphism of $\Ee$. A superconnection \cite{Qu85} on
$\Ee$ adapted to $L$ is an odd-parity differential operator $\Af$
on the graded complex $\Aa(B,\Ee) = \Aa(B)\otimes\G(B,\Ee)$, with
$\otimes=\widehat{\otimes}$ the graded tensor product, such that
\begin{equation*}
    \Af(\omega.s) = dw.s + (-1)^{|\omega|}\omega.\Af s
\end{equation*}
for $\omega\in\Aa(B), s\in\G(B,\Ee)$, and with $\Af_{[0]}=L$,
where $\Af_{[i]}$ is the component of $\Af$ which raises form
degree by $i$. The curvature of $\Af$ is the even-parity element
$\Fff = \Af^2$ of $\Aa(B,\End(\Ee))$.

It is convenient to utilize the scaling operator $\d_t\in
\Aa(B,\End(\Ee))$ of \cite{BeVe90} defined for $t>0$ by $\d_t
(\omega_{[i]}) = t^{-i/2}\omega_{[i]}$ for $\omega_{[i]} \in
\Aa^i(B),$ extended to endomorphisms of $\Aa(B,\Ee)$ by
\begin{equation}\label{e:delta-t2}
\d_t (T) = \d_t \cdot T \cdot \d_t\ii \ .
\end{equation}
We then have the scaled superconnection
\begin{equation*}
\Af_t = t^{1/2} \d_t (\Af) = t^{1/2}L + \Af_{[1]} +
t^{-1/2}\Af_{[2]} + \ldots
\end{equation*}
and curvature operator $\Fff_t = t^{1/2} \d_t (\Fff) = tL^2 +
\Ff_t$, where $\Ff_t$ raises exterior degree.

The heat operator
 $$e^{-\Fff_t}\in \Aa(B,\End(\Ee))$$ is defined by
the contour integral
\begin{equation}\label{e:heat}
e^{-\Fff_t} = \frac{i}{2\pi}\int_{\Cc} e^{-\la} \ (\Fff_t - \la
I)\ii \ d\la \ ,
\end{equation}
where $I\in \Aa(B,\End(\Ee))$ is the identity automorphism and
$\Cc$ is a contour surrounding the positive real axis $\Rf_+$
coming in on a ray with argument in $(0,\pi/2)$, encircling the
origin, and leaving on a ray with argument in $(-\pi/2,0)$. The
resolvent operator $(\Fff_t - \la I)\ii $ is defined for any
$\la\in\Cf\backslash\Rf_+$. To see this, observe that
$(\Fff_t)_{[0]} = tL^2$ is a positive Hermitian operator and hence
$(tL^2 - \la I)\ii $ is defined for any
$\la\in\Cf\backslash\Rf_+$, while the remaining terms are of
non-zero form degree. Precisely,
\begin{eqnarray}
(\Fff_t - \la I)\ii & = & (\Ff_t + tL^2 - \la I)\ii \nonumber\\
& = & (tL^2 - \la I)\ii (\Ff_t(tL^2 - \la I)\ii + I)\ii \nonumber\\
& = & \sum_{k=0}^{\dim B}(-1)^k (tL^2 - \la I)\ii (\Ff_t(tL^2 -
\la I)\ii)^k  \label{e:resolvent} \ ,
\end{eqnarray}
where \eqref{e:resolvent} is the finite algebraic Neumann series
expansion since $\Ff_t(tL^2 - \la I)\ii$ is nilpotent. The heat
trace asymptotics of the Chern character form $\ch(\Af_t) =
\str(e^{-\Fff_t})$ can thus be analyzed from the resolvent trace
asymptotics.

If $\Ee$ is finite-rank this is elementary; one then has in
$\Aa(B)$
\begin{equation}\label{e:chchform}
    \ch(\Af_t)
     = \frac{i}{2\pi}\int_{\Cc} e^{-\la} \ \str\left((\Fff_t - \la
I)\ii\right) \ d\la \ ,
\end{equation}
and taking $\Af_t = t^{1/2}L + \na^{\Ee}$, where $\na^{\Ee} =
\na^{\Ee^+}\oplus\na^{\Ee^-}$ is a connection with curvature
2-form $R = R_+\oplus R_-$, we can set $t$ to zero in
\eqref{e:resolvent} so that
 $(\Fff_0 - \la I)\ii  =
\sum_{k=0}^{\dim B}(-1)^k (-\la)^{-k-1} R^k$ and continuity yields
\begin{equation*}
\lim_{t\to 0}\ch(\Af_t) =  \sum_{k\geq 0} (-1)^k
\left(\frac{i}{2\pi}\int_{\Cc} e^{-\la} \, (-\la)^{-k-1} \
d\la\right) \str(R^k) =  \sum_{k\geq 0} \frac{1}{k!}
 \left[\tr(R_+^k) - \tr(R_-^k)\right] \ ,
\end{equation*}
and thus the transgression formula
becomes the classical Chern-Weil formula.


The interest in this methodology is where $L$ is replaced by a
smooth family of first-order elliptic operators. $\Fff_t$ is then
a vertical differential operator valued form and the integral
\eqref{e:chchform} must be regularized by repeated integration by
parts, replacing the resolvent by $(\Fff_t - \la I)^{-m}$ for
large enough $m$ to ensure a trace class operator. To analyze
$\ch(\Af_t)$ for small $t$ one computes the  asymptotics as
$\la\to \o$ of generalized resolvent traces of the form
$\str(F(\Delta - \la I)^{-k})$ where $F$ is a differential
operator and $\Delta$ a generalized Laplacian. As $t\too 0$ in
$\ch(\Af_t)$ this yields regularized trace terms and Wodzicki
residue traces \cite{Wo84}
--- computed for closed manifolds in \cite{Sc02}. In Thm(I) and
Thm(II) residue traces cancel in the relative heat traces except
for the regularized zeta trace correction terms. See
\cite{PaRo02} for a discussion of Wodzicki residues and
determinant bundle curvature.

\subsubsection{Comparing Chern Forms on Subbundles}

Between these two extremes there is an intermediate case worthy of
consideration with $\Ee^{\pm}$ Hilbert bundles with structure
group reduced to ${\rm Gl_1}$, the space of invertible operators
of the form $I + g$ with $g$ of trace-class. First, let $H =
H^+\oplus H^-$ be a graded Hilbert space where the orthogonal
projections $\Pi_+, I - \Pi_+$ onto $H^{\pm}$ are of infinite
rank. The infinite Grassmannian $\Gr_1(H)$ parameterizing
projections $P$ on $H$ for which $P-\Pi_+$ is trace-class has
homotopy type $\Zf\times {\rm BGl}(\o)$. By a theorem of Quillen
\cite{Qu87}, the de-Rham theorem holds for $\Gr_1(H)$ and so the
Chern character cohomology classes $\ch_k\in H^{2k}(\Gr_1(H),\Rf)$
of the canonical K-classes can be represented by differential
forms. The universal bundle $\Ee_1$, with fibre $\ran(P)$ at
$P\in\Gr_1(H)$, considered as a sub-Hilbert bundle of the trivial
bundle $H\times\Gr_1(H)$ inherits the canonical connection $P\cdot
d\cdot P$ with curvature 2-form $R=
PdPdP\in\Aa^2(\Gr_1(H),\End(\Ee_1))$. Since $dP$ is trace class
the form $\omega_{k} = \Tr(PdP^{2k})\in\Aa^{2k}(\Gr_1(H))$ is
well-defined, and up to a constant it represents the primitive
generator $\ch_k$ for $k\geq 1$. Although $e^{-R}$ is a
well-defined bounded operator, it is not trace-class, and so there
is no total Chern character form. On the other hand, $e^{-R} -
\Pi_+$ is of trace class with
\begin{equation}\label{e:Grchch}
 \Tr(e^{-R} - \Pi_+) = \index(\Pi_+\circ P) + \sum_{k\geq
1}\frac{(-1)^k}{k!} \omega_k \in \Aa(\Gr_1(H)) \  ,
\end{equation}
which is the relative Chern character $\Tr(e^{-R} - e^{-0_+})$ of
the virtual bundle $\Ee_1 - \Ee^{{\rm triv}}$, with $\Ee^{{\rm
triv}}$ the trivial subbundle $H^+\times \Gr_1(H)$ and $0_+$ the
zero operator on $H^+$.

\begin{thm}\label{t:schatten}
Let $\Ee=\Ee^+\oplus\Ee^-$ be a Hilbert bundle with connection
$\nabla^{\Ee}$ with structure group ${\rm Gl_1}$. Let
$\Gr_1(\Ee)\too B$ be the fibration with fibre the
infinite-Grassmannian $Gr_1(\Ee_z)$. Let $\Pf_1,\Pf_2$ be smooth
sections of $\Gr_1(\Ee)$, defining smooth subbundles $\Ww_1,
\Ww_2$ of $\Ee$ endowed with connections with curvature $R_i =
(\Pf_i\cdot\nabla^{\Ee}\cdot\Pf_i)^2$. Then there are canonical
superconnections $\Af_{i,t}$ associated to $\Ww_i$ and a
transgression formula for $t >0$ on $\Aa(B)$
$$\str(e^{-\Af_{1,t}^2} - e^{-\Af_{2,t}^2}) =  \index(P_1,P_2) +
\sum_{k=1}^{\dim B} \frac{(-1)^k}{k!}
 \Tr(R_1^k - R_2^k) + d\b_t \ .$$
The individual heat operators $e^{-\Af_{1,t}^2}$ are not trace
class. Here $\index(P_1,P_2)$ is the pointwise index of a
canonical family $(\Pf_1,\Pf_2) : \Ww_1 \too \Ww_2$ of generalized
Toeplitz operators with index bundle $\Ind(\Pf_1,\Pf_2)\in K(B)$.
Setting $\eta(\Pf_1,\Pf_2) = \Tr(e^{-R_1} - e^{-R_2}) \ ,$ one has
in $H^*(B)$
\begin{equation*}
\ch(\Ind(\Pf_1,\Pf_2))  =  \eta(\Pf_1,\Pf_2)  \ .
\end{equation*}
\end{thm}

\vskip 3mm

The proof of this theorem, of which \eqref{e:Grchch} is obviously
a special case, can be extracted from the proof of Theorem (I).
\thmref{t:schatten} settles the question posed by Freed
\cite{Fr88} on the relation of the Schatten class Chern character
forms and the constructions in Bismut \cite{Bi86}, Bismut-Freed
\cite{BiFr85}. In particular, it provides an explicit
superconnection construction of generators of the cohomology of
the loop group ${\rm LG}$ of a compact group $G$.

\subsection{Curvature of the Determinant Line Bundle}
In the latter paper \cite{BiFr85} the curvature of the determinant
bundle of a family of Dirac operators is identified with the
2-form component of the Chern character form. Likewise here, an
immediate consequence of the transgression formulas for Thm(I),
Thm(II) and \thmref{t:schatten} are the corresponding determinant
line bundle curvature formulas (see \cite{Sc03}). In particular,
since the central extension ${\rm \widehat{LG}} \too {\rm LG}$ of
the loop group is the determinant line bundle of the family
$\{g\Pi_{>} g\ii\circ\Pi_{>} \ | \ g\in {\rm LG}\}$, the
superconnection curvature form computed here is precisely the
canonical symplectic 2-form computed in Pressley-Segal
\cite{PrSe86} defining the central extension.



\section{Preliminaries}

Let $X$ be a compact manifold with boundary $Y = \dd X$ and let
$\pi : M\too B$ be a smooth Riemannian fibration of manifolds with
fibre $X_z$ diffeomorphic to $X$, restricting on the boundary $N =
\dd M$ to a smooth fibration $\pi^N : N\too B$ of closed manifolds
with fibre $Y_z$ diffeomorphic to $Y$. Let $\Ef=\Ef^+\oplus\Ef^-$
be a vertical super bundle over $M$ with metric and let
$$\Ef^{'} = \Ef_{|\N}$$
be the restricted vertical bundle over $N$ graded as
\begin{equation*}
\Ef^{'} = \Ef^0\oplus\Ef^1
\end{equation*}
where
\begin{equation*}
 \Ef^0 = \Ef^{+}_{|\N} \ \ \ {\rm and} \ \ \ \Ef^1 = \Ef^{-}_{|\N}
 \ .
\end{equation*}
We assume a collar neighborhood $$\Uu= [0,1) \times N$$ of the
boundary such that each bundle  splits isometrically. Thus if
$g^{\Ef^+}, g^{\Ef^0}$ denote the metrics on $\Ef^{+}, \Ef^0$,
then
\begin{equation}\label{e:collarE}
\Ef^{+}_{|\Uu} = p^*(\Ef^0) \cong [0,1) \times \Ef^0 \ , \hskip
10mm g^{\Ef^+}_{|\Uu} = p^*(g^{\Ef^0}) \ ,
\end{equation}
where $$p : \Uu \too N$$ is the canonical projection map, and so
forth.
In particular, the tangent bundle $T(M/B)$ along the fibres of $M$
is assumed to be oriented and endowed with a metric $g^M$ such
that
\begin{equation}\label{e:collarmetric}
g^{M}_{|\Uu} = du^2 + g^N \ ,
\end{equation}
where  $g^N$ is the induced metric on $T(N/B)$.

Associated to $\Ef$ is the graded infinite-dimensional $\Ci$
Frechet bundle over $B$
\begin{equation*}
\pi_{*}(\Ef) = \pi_{*}(\Ef^+)\oplus\pi_{*}(\Ef^-)
\end{equation*}
 with fibre $\Gamma(X_z,E_z)$ at
$z\in B$, where $E_z = \Ef_{|X_z}$, defined by
\begin{equation*}
\Gamma(B,\pi_{*}(\Ef)) = \Gamma(M,\Ef) \ .
\end{equation*}
More generally, we have the complex of $\Ci$ forms on $B$ with
values in $\pi_{*}(\Ef)$
\begin{equation*}
\Aa(B,\pi_{*}(\Ef)) =  \G(M,\pi^{*}(\wedge T^{*}B) \otimes \Ef)
 \ ,
\end{equation*}
where $\otimes = \widehat{\otimes}$ is the graded tensor product.
Likewise over $N= \dd M$ there is the bundle
\begin{equation*}
\pi^N_{*}(\Ef^{'}) = \pi^N_{*}(\Ef^0)\oplus\pi^N_{*}(\Ef^1)
\end{equation*}
with fibre $\Gamma(Y_z,N_z)$ and the complex
\begin{equation*}
\Aa(B,\pi_{*}(\Ef^{'})) = \Aa(N,(\pi^N)^{*}(\wedge T^{*}B) \otimes
\Ef^{'}) \ .
\end{equation*}

An element of $\Aa(B,\End(\pi_{*}(\Ef)))$ is an endomorphism of
$\Aa(B,\pi_{*}(\Ef)) $ which supercommutes  with the (pull-back)
action of $\Aa(B)$ by exterior multiplication. The inclusion
$$\Aa(M,\End(\Ef))\subset \Aa(B,\End(\pi_{*}(\Ef)))$$ is strict.
Specifically, a vertical differential operator
$\Df\in\Aa(B,\End(\pi_{*}(\Ef)))$ on $\Ef$ is a differential
operator on $\Aa(B,\pi_{*}(\Ef)) $ which supercommutes with
$\Aa(B)$. Then $\Df$ is a smooth section of the bundle $\Dd(\Ef)$
over $B$ with fibre the space of differential operators on
$\Gamma(X_z,E_z)$ and, by definition, a smooth family of elliptic
differential operators on $\Ef$ if elliptic on each fibre. As such
$\Df$ decomposes as
\begin{equation*}
\Df = \sum_{i=0}^{\dim B} \Df_{[i]} \ ,
\end{equation*}
where $\Df_{[i]} \in \Aa^i(B,\Dd(\Ef))$ raises form degree by $i$.

Over the closed boundary $N$ we have the graded bundle
$\Psi(\Ef^{'})$ of boundary $\pdo$s with fibre the space
$\Psi(Y_z,E^{'}_z)$ of $\pdo$s on $\Gamma(Y_z,E^{'}_z)$. Let
$\Psi^r(\Ef_N)$ be the subbundle of $\pdo$s of order $r$. A smooth
family of $\pdo$s on $\Ef^{'}$ is a smooth section $\Tf = \{T_z\in
\Psi^r(Y_z,E_z) \ | \ z\in B\}$ of $\Psi(\Ef^{'})$, identified
with a $\pdo$ on $\Gamma(N,\Ef^{'})$ which commutes with the
pull-back action of $\Ci(B)$.

The restriction of sections of $\Ef$ over $M$ to boundary sections
of $\Ef^{'}$ over $N=\dd M$ defines an even parity vertical bundle
map
\begin{equation*}
\g \in  \Aa(B,\Hom(\pi_{*}(\Ef),\pi^N_{*}(\Ef^{'})) \ ,
\end{equation*}
$$\g(\omega\otimes\psi) =  \omega\otimes\psi_{|N}\ , $$
where $\omega\in\Aa(B), \psi\in \G(B,\pi_*(\Ef))$.

We consider a smooth family of formally self-adjoint first-order
elliptic differential operators $\Ds = \{D_z \ | \ z\in B\} \in
\Aa^0(B,\Dd(\Ef))$ of odd parity with respect to the grading on
$\Aa(B,\End(\pi_*(\Ef)))$, so
\begin{equation*}
\Ds =
\begin{bmatrix}
  0 & \Ds^- \\
  \Ds^+ & 0 \\
\end{bmatrix} \ ,
\end{equation*}
where $\Ds^{\pm} = \{D^{\pm}_z \ | \ z\in
B\}\in\Aa(B,\Hom(\Ef^{\pm},\Ef^{\mp}))$ are formal adjoints --
that is, symmetric for the induced metrics on $\pi_*(\Ef^{\pm})$
for sections with support disjoint from from boundary. $\Ds$
extends to all of $\Aa(B,\Dd(\Ef))$ by $$\Ds(\omega. s) =
(-1)^{|\omega|}\omega.\Ds s \ , $$ where $\omega \in
\G(M,\pi^*(\wedge T^* B)), s \in \G(M,\Ef)$.

We assume $\Ds^{\pm}$ have product structure in the collar $\Uu$
\begin{equation}\label{e:product-type}
\Ds^{\pm}_{|\Uu} = \Upsilon_{\pm}\left(\frac{\dd}{\dd u} +
\dd_{\pm}\right) \ ,
\end{equation}
where $$\dd_+ = \{\dd_{+,z} \ | \ z\in B\}\in
\Gamma(B,\Dd(\Ef^{0})) \ , \hskip 15 mm \dd_- = \{\dd_{-,z} \ | \
z\in B\}\in \Gamma(B,\Dd(\Ef^{1})) $$ are smooth boundary families
of self-adjoint first-order elliptic differential operators on
$\Ef^{'}$ and $\Upsilon^+ = -(\Upsilon^-)^* : \Ef_{|\Uu} \too
\Ef_{|\Uu}$ is a unitary bundle isomorphism. Then $\Ds$ has
product structure in $\Uu$
\begin{equation}\label{e:Dscollar}
\Ds_{|\Uu} = \begin{bmatrix}
  0 & \Upsilon_{-} \\
  \Upsilon_{+} & 0 \\
\end{bmatrix}\left(\frac{\dd}{\dd u} + \begin{bmatrix}
\dd_{+} & 0 \\
  0 & \dd_- \\
\end{bmatrix}\right) \ .
\end{equation}
For example, the isometric splitting of the metric
\eqref{e:collarmetric} ensures this in the case of a family of
Dirac operators associated to a fibration of spin manifolds, as in
\cite{BiCh90,MePi97}.

\vskip 1mm

Vertical APS-boundary problems for $\Ds$ arise as follows. For the
moment we restrict our attention to $\Ds$ acting on
$\G(B,\pi_*(\Ef)) = \Aa^0(B,\pi_*(\Ef))$. Associated to $\Ds$ is
the subspace of interior solutions
\begin{equation*}
\Ker(\Ds) = \{\psi\in \G(B,\pi_{*}(\Ef)) \ | \ \Ds\psi = 0 \ \
{\rm in} \ \ M\backslash N\} \ .
\end{equation*}
Since $\Ds$ is a bundle endomorphism of $\pi_{*}(\Ef)$,  a
vertical differential operator on $\Ef$, then $\Ker(\Ds)$ is
fibered by the $\Ci$ kernels $\Ker(D_z)  = \{\psi_z\in \G(X_z,E_x)
\ | \ D_z\psi_z = 0 \ {\rm in} \ X_z\backslash Y_z\}$. The space
of vertical Cauchy data is defined by restriction of interior
solutions to the boundary
\begin{equation*}
\Ww(\Ds) : = \g \Ker(\Ds) = \{\xi\in \pi_{*}(\Ef^{'}) \ | \ \xi =
\g\phi, \ \phi \in \Ker(\Ds)\} \ ,
\end{equation*}
and fibered $\Ww(\Ds) = \cup_z W(D_z)$ by the pointwise Cauchy
data spaces $W(D_z) = \g\Ker(D_z)$ of classical boundary $\pdo$
theory. These spaces are naturally graded as
\begin{equation*}
\Ker(\Ds) = \Ker(\Ds^+)\oplus \Ker(\Ds^-) \ , \hskip 10mm \Ww(\Ds)
= \Ww(\Ds^+)\oplus \Ww(\Ds^-)
\end{equation*}
with respect to the gradings on $\pi_{*}(\Ef), \pi_{*}(\Ef^{'})$.
For manifolds with boundary,  $\Ker(D_z)$ is infinite-dimensional,
contrasting with elliptic regularity on closed manifolds. On the
other hand, and again unlike closed manifolds, the kernels
$\Ker(D_z)$ vary smoothly with $z$:
\begin{prop}\label{p:infinitebundle}
$\Ker(\Ds)$ and $\Ww(\Ds)$ are, respectively, smooth subbundles of
$\pi_{*}(\Ef)$ and $\pi^N_{*}(\Ef^{'})$.\footnote{By a $C^k$
vector bundle $\Ee\to B$ of infinite-rank we mean in the weak
sense: for $z\in B$ there is an open set $U_z$ and a
fibre-preserving $C^k$-diffeomorphism $\Ee_{|U_z} \cong U_z \times
\Ee_z$, where $\Ee_z$ is the fibre of $\Ee$ at $z$.} The
restriction map defines a vector bundle isomorphism over $B$
\begin{equation}\label{e:Ker=H(D)}
\Ker(\Ds) \cong \Ww(\Ds)  \ .
\end{equation}
Equivalently, $\Ker(\Ds^{\pm})$  are smooth subbundles of
$\pi_{*}(\Ef^{\pm})$, $\Ww(\Ds^+)$ a subbundle of
$\pi^N_{*}(\Ef^{0}), $ and $\Ww(\Ds^-)$  a subbundle of
$\pi^N_{*}(\Ef^{1}), $ with isomorphisms $\Ker(\Ds^{\pm}) \cong
\Ww(\Ds^{\pm}) $.
\end{prop}
\begin{proof}
The assertion is that there is an exact sequence
\begin{equation}\label{e:Kerexact}
0 \too \Ker(\Ds) \too \pi_{*}(\Ef) \stackrel{\Ds}{\too}
\pi_{*}(\Ef) \too 0
\end{equation}
and that it is canonically split. (This holds also on all Sobolev
completions \cite{Sc03}.) Let $\widehat{M} = M \cup_N(-M)\to B$ be
the fibration of closed double manifolds constructed from $M$,
with fibre $\widehat{X}_z = X_z \cup_{Y_z} (-X_z)$. With the
product structure \eqref{e:Dscollar}, $\Ds$ extends by a standard
argument to an invertible vertical first-order differential
operator $\widehat{\Ds}$ on $\G(\widehat{M},\widehat{\Ef})$, where
$\widehat{\Ef}_{|M} = \Ef$ and $\widehat{\Ds}_{|M} = \Ds$. Using
the continuous extension operator $\bf{e}:\pi_{*}(\Ef)\to
\pi^{L^2}_{*}(\widehat{\Ef})$ to the bundle of $L^2$ sections of
$\widehat{\Ef}$, where $\bf{e}(\psi) = \psi $ on $M$ and zero
elsewhere, and the restriction operator
$\bf{r}:\pi_{*}(\widehat{\Ef})\too \pi_{*}(\Ef)$, we define
\begin{equation*}
\Gs := {\bf r}\widehat{\Ds}\ii {\bf e}\in\G(B,\End(\pi_{*}(\Ef)))
\ .
\end{equation*}
Since $\widehat{\Ds}\widehat{\Ds}\ii = \If$ on
$\G(\widehat{M},\widehat{\Ef})$, with $\If$ the vertical identity
operator, by locality we have $\Ds\cdot\Gs = \If$ on $\G(M,\Ef)$
and hence $\Gs$ is a vertical right inverse of $\Ds$, and the
required splitting of \eqref{e:Kerexact}.

On the other hand, we can use the inverse $\widehat{\Ds}\ii$ to
define the vertical Poisson operator
\begin{equation*}
\Kk = {\bf r}\widehat{\Ds}\ii \widehat{\g}^* \Upsilon \in
\G(B,\Hom(\pi^{N}_{*}(\Ef^{'}),\pi_{*}(\Ef))) \ .
\end{equation*}
where $\widehat{\g}$ is the restriction map from $\widehat{M}$ to
$N$. Then $\Kk = \{K_z \ | \ z\in B\}$ is a smooth family of
operators taking boundary fields into interior fields coinciding
fibrewise with the classical Poisson operator $K_z$ of $D_z$
\cite{Gr96,Gr99}. Hence $\Kk$ is a bundle map with range in
$\Ker(\Ds)$, and which restricts to a fibrewise isomorphism
\begin{equation}\label{e:Kkisom}
\Kk : \Ww(\Ds)\stackrel{\cong}{\too} \Ker(\Ds)
\end{equation}
that by construction depends smoothly on $z\in B$. It is inverse
to $\g$ on $\Ker(\Ds)$. The vertical operator
$$P_{{\rm Ker}} = \Kk\circ \g : \pi_{*}(\Ef) \too \pi_{*}(\Ef) \ ,$$ is
therefore a canonical smooth family of projections on $\pi_{*}
(\Ef)$ with range $\Ker(\Ds)$, while
\begin{equation}\label{e:Pker-ort}
\Gs\cdot \Ds = \If - \Kk\circ \g = \If - P_{{\rm Ker}} = P_{{\rm
Ker}}\pp \in \G(B,\End(\pi_*(\Ef)))
\end{equation}
is the complementary vertical projection with range
$\Ker(\Ds)\pp$. The first equality in \eqref{e:Pker-ort} follows
from the distributional Green's formula (the single operator case
proved in \cite{Gr99} applies directly to $\Ds$). On the other
hand, the restriction of $\Kk$ to the bundle of boundary sections
\begin{equation}\label{e:P(D)}
P(\Ds) := \g\circ \Kk  \in \G(B, \Psi^0(\Ef^{'}))
\end{equation}
is the vertical Calderon projector with range $\Ww(\Ds)$. $P(\Ds)
= \{P(D_z) \ | \ z\in B\}$ is a smooth family of order 0 $\pdo$
projections restricting on each fibre to the usual Calderon
projector with range $\g\Ker(D_z) = W(D_z)$.

To see that $\Ww(\Ds)$ is a vector subbundle of
$\pi^{N}_{*}(\Ef^{'})$, because $\Pp(\Ds)$ (and $\Pker$) is a
smooth family and invertibility is an open condition in the space
of Fredholm operators, it is sufficient to show that $ P(D_z)
\circ P(D_{b^{'}}) : W(D_z)\to W(D_{b^{'}})$ is Fredholm for $z$
near $z^{'}$ in $B$, since $ P(D_z) \circ P(D_{b^{'}})$ is the
identity when $z=z^{'}$. But if $\Pi_{>,z}, \Pi_{>,z^{'}}$ are the
APS spectral projections at $z$ and $z^{'}$, it is a well known
spectral flow property that
$\Pi_{>,z}\circ\Pi_{>,z^{'}}:\ran(\Pi_{>,z^{'}})\to
\ran(\Pi_{>,z})$ is invertible for $z$ close to $z^{'}$; for,
since $\dd_z, \dd_{z^{'}}$ are self-adjoint elliptic, perturbation
theory asserts the spectrum is locally smooth, so $\dim
\Ker(\dd_z)$ is locally constant. On the other hand, it is well
known that $P(D_z) - \Pi_{>,z}$ and $P(D_{z^{'}}) - \Pi_{>,z^{'}}$
are smoothing operators \cite{Sc95}. Hence $ P(D_z) \circ
P(D_{z^{'}})$ is also locally invertible and we thus reach the
conclusion.

Finally, the bundle structure on $\Ker(\Ds)$ is immediate from
that on $\Ww(\Ds)$ and the fibre preserving isomorphism
\eqref{e:Kkisom}.
\end{proof}

Evidently, the (orthogonalized) Calderon section $P(\Ds)$ splits
with the grading as
\begin{equation}\label{e:PDsplits}
    P(\Ds) = \begin{bmatrix}
      P(\Ds^+) & 0 \\
      0 & P(\Ds^-) \\
    \end{bmatrix} \ ,
\end{equation}
where with obvious notation
\begin{equation*}
P(\Ds^+) := \g\circ \Kk^+  \in \G(B, \Psi^0(\Ef^0)), \hskip 10mm
\Kk^+ = {\bf r}\widehat{\Ds^+}\ii \widehat{\g}^* \Upsilon \in
\G(B,\Hom(\pi^{N}_{*}(\Ef^{0}),\pi_{*}(\Ef^+))) \ ,
\end{equation*}
and similarly $P(\Ds^-)$, are the vertical Calderon projections
for $\Ds^{\pm}$. In this way we obtain the bundle structure of
$\Ww(\Ds^{\pm})$. One has
\begin{equation}\label{e:adj}
P(\Ds^-)  = \Upsilon_+ \cdot P(\Ds^+)\pp\cdot\Upsilon_+^* = -
\Upsilon_+ \cdot P(\Ds^+)\pp\cdot \Upsilon_- \ .
\end{equation}

The (orthogonalized) Calderon section $P(\Ds)$ defines a canonical
section in $\G(B,\Psi^0(\Ef^{'})$, which by \eqref{e:adj} is
determined by the section $P(\Ds^+)$ in $\G(B,\Psi^0(\Ef^{0})$.
Relative to  $P(\Ds^+)$ we define the subspace of $\Gamma(B,\Psi^0
(\Ef^0))$ of {\em Grassmann sections}
\begin{equation*}
\Gr(\pi^{N}_{*}(\Ef^0)) = \Gamma(B,\Gr_{\o}(\dd_+)) \subset
\G(B,\End(\pi_*^N(\Ef^0)))\ ,
\end{equation*}
where $\Gr_{\o}(\dd_+)\to B$ is the fibration with fibre at $z\in
B$ the smooth Grassmannian
\begin{equation*}
  \Gr_{\o}(\dd_{+,z}) = \{P_{+,z}\in\Psi^0(Y_z , E^0_z) \ | \
  P_{+,z}^* = P_{+,z}, \ P^2_{+,z} = P_{+,z}, \ P_{+,z} - P(D^+_z)
  \in \Psi^{-\o}(Y_z , E^0_z)\} \ .
\end{equation*}
Thus a Grassmann section is a vertical projection $\Pp^+ =
\{P_{+,z} \ | \ z\in B\}$ on $\pi_{*}^{N}(\Ef^0)$, defining a
smooth family of of $\pdo$ projections, such that
\begin{equation*}
  \Pp^+ = P(\Ds^+) + \verb"S" \ ,
\end{equation*}
where $\verb"S" \in \Gamma(B,\Psi^{-\o}(\Ef^0))$ is a smooth
family of smoothing operators. This means $\verb"S"$ has kernel
$$k\in \Gamma(N\times_{\pi^N}N , \Ef^0\boxtimes \Ef^0)\ ,$$ where
$N\times_{\pi^N}N $ is the fibre product with respect to the
projection maps $pr_i : N\times N \to N$ onto the $i^{th}$
component, and $\Ef^0\boxtimes\Ef^0 =  pr_{1}^{*}(\Ef^0)\otimes
pr_{2}^{*}(\Ef^0)$.

Notice that due to boundary spectral flow the family $\Pi_{>} =
\{\Pi^{z}_{>} \ | \ z\in B\}$ is only locally continuous. However,
when $\dim\Ker(D_{Y_z})$ is constant, as \cite{BiCh90}, then
$\Pi_{>} \in \Gr(\pi^{N}_{*}(\Ef^0))$.

The primary purpose of $\Gr(\pi^{N}_{*}(\Ef^0))$ is  a parameter
space of self-adjoint global boundary conditions for the the
family $\Ds$. To this end, we associate to each $\Pp^+\in
\Gr(\pi^{N}_{*}(\Ef^0))$ the Grassmann section for $\Ds$
\begin{equation*}
    \Pp = \begin{bmatrix}
      \Pp^+ & 0 \\
      0 & \Pp^- \\
    \end{bmatrix} \ ,
\end{equation*}
where
\begin{equation}\label{e:Padj}
\Pp^-  =  - \Upsilon_+ \cdot (\Pp^+)\pp\cdot \Upsilon_- \ ,
\end{equation}
defining the smooth family of odd-parity APS-type global boundary
problems
\begin{equation*}
  \Ds_{\Pp} = \begin{bmatrix}
    0 & \Ds^{-}_{\Pp^-}  \\
    \Ds^{+}_{\Pp^+} & 0 \\
  \end{bmatrix}: \pi_{*}(\Ef | \Pp) \too \pi_{*}(\Ef) \ .
\end{equation*}
This means $$\Ds_{\Pp} = \Ds \ , \hskip 10mm \Ds^{+}_{\Pp^+} =
\Ds^{+} \ , \hskip 10mm \Ds^{-}_{\Pp^-} = \Ds^{-}$$ as operators,
but with graded domain
\begin{equation*}
   \pi_{*}(\Ef | \Pp) = \pi_{*}(\Ef^+ | \Pp^+) \oplus
   \pi_{*}(\Ef^- |\Pp^-)\ ,
\end{equation*}
where
\begin{equation*}
\pi_{*}(\Ef | \Pp) = \Ker(\Pp \circ \g :\pi_{*}(\Ef) \too
\pi^N_{*}(\Ef^{'}))
\end{equation*}
and
\begin{equation*}
\pi_{*}(\Ef^+ | \Pp^+) = \Ker(\Pp^+ \ \circ \ \g :\pi_{*}(\Ef^+)
\too \pi^N_{*}(\Ef^0))\ ,\hskip 5mm  \pi_{*}(\Ef^- |\Pp^-) =
\Ker(\Pp^- \ \circ \ \g :\pi_{*}(\Ef^-) \too \pi^N_{*}(\Ef^1))\ .
\end{equation*}
\vskip 1mm Over $X_z$, $\Ds_{\Pp}$ restricts to the self-adjoint
APS-type operator $$D_{P_z} := (D_z)_{P_z} :\dom(D_{P_z}) \too
\G(X_z,E_z) \ ,$$ in the usual single operator sense.

The definition \eqref{e:Padj} ensures $\Ds^{-}_{\Pp^-}$ is the
adjoint family to $\Ds^{+}_{\Pp^+}$, in the obvious sense. Note
that
\begin{equation*}\label{e:Ppsadj}
\Pp  =  - \Upsilon\cdot \Pp \pp\cdot \Upsilon \ .
\end{equation*}
$\Pp$ extends naturally to the even parity boundary operator
\begin{equation*}
\Pp = 1\otimes\Pp\in \Aa(B,\End(\Ef^{'})) \ , \hskip 10mm
\Pp(\omega\otimes s) = \omega\otimes \Pp s \ ,
\end{equation*}
and which hence extends the smooth family of APS-type boundary
problems to
\begin{equation*}
\Ds_{\Pp } : \Aa(B,\pi_{*}(\Ef|\Pp)) \too \Aa(B,\pi_{*}(\Ef)) \ ,
\end{equation*}
where $\Aa(B,\pi_{*}(\Ef|\Pp)) = \Aa(M,\pi^*(\wedge T^* B))\otimes
\G(M,\Ef|\Pp))$ is the subspace of sections annihilated by
$\Pp\circ \g$, with $\g$ the extended restriction map
\begin{equation}\label{e:restriction}
\g : \Aa(B,\pi_{*}(\Ef)) \too \Aa(B,\pi^{N}_{*}(\Ef^{'})) \hskip
10mm \g(\omega\otimes s) = \omega\otimes \g s \ .
\end{equation}

\vskip 4mm

A Grassmann section $\Pp = \{P_z \ | \ z\in B\}$ distinguishes the
graded boundary subspace
\begin{equation*}
\Ww = \Ww^+ \oplus \Ww^- \subset \pi_{*}^N(\Ef^{'})
\end{equation*}
with
\begin{equation*}
\Ww = \ran(\Pp) , \hskip 10mm \Ww^+ = \ran(\Pp_0) , \hskip 10mm
\Ww^- = \ran(\Pp^-) \ ,
\end{equation*}
and fibre $W_z = \ran(P_z)$ at $z\in B$. $\Ww, \Ww^{\pm}$ are
canonically related to the spaces $\pi_{*}(\Ef | \Pp),
\pi_{*}(\Ef^{\pm}|\Pp^{\pm})$ of interior sections by the exact
sequences
\begin{equation}\label{e:exactW}
0 \too \pi_{*}(\Ef | \Pp) \too \pi_{*}(\Ef) \stackrel{\Pp\circ \g
}{\too} \Ww \too 0
\end{equation}
and
\begin{equation*}
0 \too \pi_{*}(\Ef^{\pm} | \Pp^{\pm}) \too \pi_{*}(\Ef^{\pm})
\stackrel{\Pp^{\pm}\circ \g }{\too} \Ww^{\pm} \too 0 \ .
\end{equation*}
\vskip 1mm \noi
\begin{prop}\label{p:Pbundles}
$\pi_{*}(\Ef|\Pp)$, $\pi_{*}(\Ef^{\pm}|\Pp^{\pm})$ and $\Ww,
\Ww^{\pm}$ are smooth infinite-dimensional Frechet vector bundles
on $B$.
\end{prop}
\begin{proof}
We restrict our comments to $\pi_{*}(\Ef|\Pp)$ and $\Ww$. Because
$$(P_{z_0},P_z) = P_z\circ P_{z_0} : W_{z_0} \too W_z$$ is Fredholm
for $z$ close to $z_0\in B$ and is the identity operator for
$z=z_0$, then, since invertibility is an open condition in
$\textsf{Fred}$, the operator $(P_{z_0},P_z)$ is invertible for
$z$ close to $z_0$. Since $P_z-P(D_z)\in\Psi^{\o}(Y_z,E_z)$, the
Fredholm property for $(P_{z_0},P_z)$ is equivalent to that of
$(P(D_{z_0}),P(D_z))$ and this is proved in
\propref{p:infinitebundle}. By its local invertibility,
$(P_{z_0},P_z)$  defines around $z_0$ a local trivialization of
$\Ww$  modelled on $W_{z_0}$ bounded in each $C^l$ norm. Since
$z_0$ was arbitrary $\Ww$ is hence a smooth subbundle of the
Frechet bundle $\pi_{*}^{N}(\Ef^{'})$.

The bundle structure on $\pi_{*}(\Ef | \Pp)$ is defined by a
Fredholm splitting of the exact sequence \eqref{e:exactW}. Let
$\chi$ be a $\Ci$ cut-off function: $\chi(t) = 1$ when $0\leq t <
1/4$ and $\chi(t) = 0$ when $t > 3/4$. Then we define a vertical
Poisson operator extending boundary sections into interior
sections with support in the collar neighborhood $\Uu$
\begin{equation*}
    \Ks \in \Aa(B,\Hom(\pi_{*}^N(\Ef^{'}),\pi_{*}(\Ef))) \
\end{equation*}
by defining $\Ks : \G(N,\Ef^{'}) \too \G(M,\Ef)$ by
\begin{equation}\label{e:Kk}
  \Ks(x) =
  \begin{cases}
    \chi(u). e^{-u\dd_y^2} &  \ x = (u,y)\in \Uu, \\
    0 & \ x \in M\backslash \Uu \ .
  \end{cases}
\end{equation}
Here $e^{-u\dd^2}\in \Aa(B,\Psi^{-\o}(\Ef^{'}))$ is the heat
operator defined by the family of boundary Laplacians $\dd$, and
as such we have
\begin{equation}\label{e:gK}
  \g\circ \Ks  = \lim_{u\too 0}  \chi(u). e^{-u\dd^2}(y) = \If_N \ ,
\end{equation}
where $\If_N$ is the vertical identity operator on
$\pi^N_{*}(\Ef^{'})$. We use  $\Ks$ to define the operator
\begin{equation}\label{e:P}
\Ps = \Ps(\Pp) := \If -  \Ks\cdot\Pp\cdot\g
 \in\Aa(B,\End(\pi_{*}(\Ef)))\ ,
\end{equation}
with $\If$ the vertical identity operator on $\pi_{*}(\Ef)$. Then
$$\Ps(\Pp) = \{\Ps(P_z) := I_z -  \Ks_z\cdot P_z \cdot\g \ | \ z\in B\}$$
is a smooth family of projection operators with range $\pi_{*}(\Ef
| \Pp)$:
\begin{equation*}
    \Ps(\Pp)^2 = \Ps(\Pp) \ ,
\end{equation*}
\begin{equation*}
 \Pp\g \Ps(\Pp) = 0   \ ,
\end{equation*}
\begin{equation*}
 \Ps(\Pp) = \If  \ \ \ {\rm on} \ \ \  \pi_{*}(\Ef | \Pp) \ .
\end{equation*}
The first two identities are readily verified from \eqref{e:gK};
note that the second says that $\ran(\Ps(P_z)) = \G(X_z, E_z |
P_z)$.

It is enough to prove the bundle structure of $\pi_{*}(\Ef|\Pp)$
over an open ball $U$ in $B$ endowed with a local trivialization
$M_{|U} = X_{z_0} \times U$, some $z_0 \in U$. This defines a
trivialization (in the weak sense) $ \pi_{*}(\Ef)_{|U} =
\G(X_{z_0},E_{z_0}) \times U$ with respect to which the local
family of APS-type projections $\Pp = \{P_z \in
\Gr(\G(Y_{z_0},E^0_{z_0}) \ | \ z\in U\}$ may be taken in a fixed
Grassmannian. In a neighborhood of $z_0\in B$ in $U$ we can
therefore define the local map between fibres of $\pi_{*}(\Ef |
\Pp)_{|U}$
\begin{equation}\label{e:P4}
\Ps(P_z)\cdot\Ps(P_{z_0}) : \G(X_{z_0}, E_{z_0} | P_{z_0}) \too
\G(X_{z_0}, E_{z_0} | P_z)\ .
\end{equation}
From the first part of the proof, a local trivialization of $\Ww$
is defined by the smooth family of invertible operators
$P_{z_0}\circ P_z : W_{z} \too W_{z_0}$ in a small neighborhood
$U_0$ of $z_0$. An easy check then shows that the family of
operators
$$
I_z -  \Ks_z P_z(P_{z_0}\circ P_z)\ii P_{z_0} \cdot\g
$$
defines a two-sided inverse to $\Ps(P_z)\cdot\Ps(P_{z_0})$ in
$U_0$. Smoothness is obvious as compositions of smooth families of
operators. The isomorphisms \eqref{e:P4} hence define local weak
trivializations of $\pi_{*}(\Ef | \Pp)$ and this completes the
proof.
\end{proof}

By \propref{p:Pbundles}, $\Pp_i, \Pp_j\in \Gr(\pi^{N}_{*}(\Ef^0))$
define the canonical smooth Hermitian family of odd-parity
Fredholm operators
\begin{equation}\label{e:pipj}
L_{i,j} = \begin{bmatrix}
  0 &   \Pp_i \cdot \Pp_j \\
  \Pp_j \cdot \Pp_i & 0 \\
\end{bmatrix} \in \Aa^0 (B, \End(\Ww_i \oplus \Ww_j)) \ ,
\end{equation}
where $\Pp_j \cdot \Pp_i $ parameterizes the operators $P_{j,z}
\circ P_{i,z}: W_{i,z} \too W_{j,z}$. We may also write
\begin{equation*}
    L_{i,j} = (\Pp_i,\Pp_j) \ .
\end{equation*}
In particular, associated to $\Ds_{\Pp}$ is the vertical
`scattering' operator
\begin{equation*}
\Sss(\Pp) = \begin{bmatrix}
  0 &   \Sss(\Pp^+)^* \\
  \Sss(\Pp^+)  & 0 \\
\end{bmatrix} \in \Aa^0 (B, \End(\Ww(\Ds^+) \oplus \Ww^+)) \ ,
\end{equation*}
where $\Sss(\Pp^+) = \Pp^+ \cdot P(\Ds^+) : \Ww(\Ds^+) \too
\Ww^+$.

\vskip 2mm

 The choice of a Grassmann section $\Pp$ restricts $\Ds$ to a
family of Fredholm operators. It also has the consequence that the
kernels of the restricted operators no longer define a vector
bundle. They do, however, still define a virtual bundle:
\begin{prop}\label{p:indexbundle1}
$\Ds_{\Pp}$ defines a smooth family of Fredholm operators with
kernel and cokernel consisting of smooth sections, and hence an
index bundle
\begin{equation*}
\Ind(\Ds_{\Pp}) \in K(B) \ .
\end{equation*}
Similarly, $L_{1,2} = (\Pp_1, \Pp_2)$ is a smooth Toeplitz family
of Fredholm operators with kernel and cokernel consisting of
smooth sections defining an index bundle
\begin{equation*}
\Ind(\Pp_1,\Pp_2) \in K(B) \ .
\end{equation*}
\end{prop}

\noi The principal relations between these K-theory elements are
as follows.

\begin{prop}\label{p:indexbundle2}
As elements of $K(B)$
\begin{equation}\label{e:indexbundle3}
\Ind(\Ds_{\Pp}) = \Ind(\Ss(\Pp))
\end{equation}
and
\begin{equation*}
\Ind(\Pp_1,\Pp_3) - \Ind(\Pp_1,\Pp_2) = \Ind(\Pp_2,\Pp_3)\ ,
\end{equation*}
and hence
\begin{equation}\label{e:indexbundle5}
\Ind(\Ds_{\Pp_1}) - \Ind(\Ds_{\Pp_2}) = \Ind(\Pp_2,\Pp_1) \ .
\end{equation}
\end{prop}

\noi Proofs for spectral sections can be accessed in
\cite{MePi97,DaZh97}. The general case is proved in \cite{Sc03} by
a quite different method, extending to the corresponding
determinant bundle isomorphisms. We omit further comment as the
details are not relevant here.

\vskip 2mm

The objective is to construct canonical differential form
representatives for the K-class $\Ind(\Ds_{\Pp})$ by computing
Chern character forms induced from the superconnections on
$\pi_*(\Ef)$ and $\pi_*^{N}(\Ef^0)$ defined as follows. A
superconnection on $\pi_*(\Ef)$ adapted to $\Ds$ is a differential
operator on $\Aa(B,\pi_*(\Ef))$ of odd-parity with respect to the
induced $\Zf_2$-grading such that
\begin{equation}\label{e:superconnection}
 \Af (\omega\psi) = d\omega\psi + (-1)^{|\omega|}\omega\Af(\psi)
\end{equation}
for $\omega\in\Aa(B), \psi\in \Aa(B,\pi_*(\Ef))$, with $\Af_{[0]}
= \Ds$. Notice \eqref{e:superconnection} implies $\Af_{[1]}$ is a
connection on $\pi_*(\Ef) = \pi_*(\Ef^+)\oplus \pi_*(\Ef^-)$
preserving the grading.

The induced superconnection on $\Aa(B,\End(\pi_*(\Ef)))$ is
defined by $\Af a = [\Af,a]$, where $[ \ , \ ]$ is the
supercommutator. The curvature of $\Af$ is the vertical
differential operator with differential form coefficients $\Fff =
\Af^2 \in\Aa(B,\Dd(\Ef))$ with degree zero component the
Dirac-Laplacian
\begin{equation*}
\Fff_{[0]} = \Ds^2 =
\begin{bmatrix}
  \Ds^-\Ds^+ & 0 \\
  0 & \Ds^+\Ds^- \\
\end{bmatrix}
 \ .
\end{equation*}
For $t> 0$ we consider the scaled superconnection
\begin{equation*}
\Af_t = t^{1/2} \d_t (\Af) = t^{1/2}\Ds + \Af_{[1]} +
t^{-1/2}\Af_{[2]} + \ldots
\end{equation*}
with curvature operator $\Fff_t = t^{1/2} \d_t (\Ff) = t\Ds^2 +
F_t$ where $F_t$ raises exterior degree.

A superconnection on $\pi_*(\Ef^{'})$ defines by restriction a
superconnection $ \Af^{'}$ on $\pi_*^{N}(\Ef^{'})$ adapted to
$\dd$. We assume throughout the compatibility condition
\begin{equation}\label{e:pullback2}
 \Af_{|\Uu} = p^*(\Af^{'}) \ ,
\end{equation}
where $p:\Uu\too N$ is the canonical projection map. By
functoriality $\Fff_{|\Uu} = p^* (\Fff^{'}),$ where $\Fff^{'} =
(\Af^{'})^2\in\Aa(B,\Dd(\Ef^0))$ is the boundary curvature.

\vskip 2mm

Though much of what follows holds for a general superconnection,
Theorems (I) and (II) are proved with the canonical
superconnection induced from a choice of connection on
$\pi_*(\Ef)$.  This is defined by a compatible connection
$\na^{\Ef}$ on $\Ef$ over $M$ such that
\begin{equation}\label{e:pullback4}
\na^{\Ef}_{|\Uu} = p^* (\na^{\Ef^{'}}) \ ,
\end{equation}
with $\na^{\Ef^{'}}$ the induced boundary connection on $\Ef^{'}$,
along with a connection on the fibration  $\pi : M \too B$,
meaning a choice of splitting
\begin{equation}\label{e:fibrationsplit}
T^* M = T^* (M/B) \oplus T_{H}^*  M \ ,
\end{equation}
and hence an isomorphism $\tau:T_{H}^*  M  \too \pi^*(T^*B)$. The
connection $$\nE: \Aa^0 (B,\pi_*(\Ef)) \too \Aa^1 (B,\pi_*(\Ef)) \
$$  is then defined by the composition
\begin{equation}\label{e:connection}
 \G(M,\Ef) \stackrel{\na^{\Ef}}{\too} \G(M,T^* M \otimes\Ef)
\stackrel{\tau}{\too} \G(M,\pi^*(T^*B) \otimes\Ef) \ .
\end{equation}
Hence
\begin{equation}\label{e:canonicalsuperconnection}
\Af_t = \nE + t^{1/2}\Ds = \begin{bmatrix}
\na^{\pi_{*}(\Ef^+)}   &  t^{1/2}\Ds^- \\
  t^{1/2}\Ds^+  & \na^{\pi_{*}(\Ef^-)}\\
\end{bmatrix}
\end{equation}
maps $\G(M,\Ef^{\pm}) \too \A^{\pm}(M,\Ef)$, where  $\pm$ refer to
the $\Zf_2$-gradings,  and thus extends to a superconnection on
all of $\Aa(B,\pi_*(\Ef))$ by setting $\Af (\omega s ) = d\omega s
+ (-1)^{|\omega|}\omega\Af s$, for $\omega\in \pi^*(\Aa(B)), s\in
\G(M,\Ef)$. Assumption \eqref{e:pullback2} is assured by
\eqref{e:product-type} and \eqref{e:pullback4}.

In this case we have
\begin{equation}\label{e:canonicalcurv}
    \Fff_t = \Rs^{\pi_*(\Ef)} + t^{1/2}\nE \Ds + t\Ds^2 \ ,
\end{equation}
where the curvature $\Rs^{\pi_*(\Ef)} = (\nE)^2
\in\Aa(B,\Dd(\Ef))$ of $\nE$ is the vertical first-order
differential operator valued 2-form on $B$
\begin{equation}\label{e:canonicalRcurv}
    \Rs^{\pi_*(\Ef)}  =  -\nE_{\sigma} + R^{(M/B)} \ .
\end{equation}
Here $\sigma(\xi_1,\xi_2) = -P_H [\xi^{H}_1,\xi^{H}_2]$, with
$P_H$ the projection onto $T_H M$ arising from the dual splitting
to \eqref{e:fibrationsplit} and $\xi^H$ the corresponding lifting
of a vector field $\xi$ on $B$. $R^{(M/B)}$ is the curvature
tensor of $P_{M/B}\cdot \na^{TM} \cdot P_{M/B}$, where $\na^{TM}$
is the Levi-Civita connection for $g^M + \pi^*(g^B)$ for any
metric on the base, and $P_{M/B}$ the projection onto $T(M/B)$
complementary to $P_H$.
  Note that $\nE \Ds := [\nE ,\Ds]$ in \eqref{e:canonicalcurv}.

The compatibility of metric, connections all other geometrical
structures in the collar neighborhood $\Uu$ of $N$ in $M$ means
that the homolog of each of these objects is induced on the
boundary fibration of closed manifolds. Thus we have a splitting
\begin{equation*}
T^* N = T^* (N/B) \oplus T_{H}^*  N
\end{equation*}
which combined with  \eqref{e:pullback2} defines the connection
$$\na^{\pi_*^N(\Ef^{'})} = \na^{\pi_*^N(\Ef^{0})}\oplus \na^{\pi_*^N(\Ef^{1})}
 : \Aa^i (B,\pi_*(\Ef^{'})) \too \Aa^{i+1}
(B,\pi_*(\Ef^{'}))$$ on $\pi_*^N(\Ef^{'}),$ and hence the
superconnection $\Af^{'}_t = \na^{\pi_*^N(\Ef^{'})} + t^{1/2}\dd$
which the Bismut superconnection extends.

\vskip 3mm

\section{Chern Character Forms from the Boundary}

The K-theory identities of \propref{p:indexbundle2} mean that to
construct a canonical Chern character form representative for
$\ch(\Ind(\Ds_{\Pp}))$ from the geometry of the boundary fibration
of closed manifolds, it is enough to compute the Chern character
form  of $L_{i,j}$ \eqref{e:pipj} considered as a vertical $\pdo$
of order $0$ on the super bundle $\pi_{*}^{N} (\Ef^0)\oplus
\pi_{*}^{N} (\Ef^0)$.

\vskip 2mm

\subsection{Construction of a Relative Heat Operator}
We consider a superconnection $\Bf$ on the superbundle
$\pi_{*}^{N} (\Ef^0)\oplus \pi_{*}^{N} (\Ef^0)$ adapted to an odd
family of self-adjoint pseudodifferential operators $\Tf  =
\begin{bmatrix}
  0 & \Tf_- \\
  \Tf_+ & 0 \\
\end{bmatrix}\in \Aa^0 (B, \Psi(\pi_{*}^{N} (\Ef^0)\oplus \pi_{*}^{N}
(\Ef^0)))$. Thus $\Bf$ is a differential operator of odd-parity
$\Bf = \sum_{i=0}^{\dim B}\Bf_{[i]} $ on $\Aa(B, \pi_{*}^{N}
(\Ef^0)\oplus \pi_{*}^{N} (\Ef^0))$ satisfying $$ \Bf(\omega.\phi)
= d_B \omega . \phi + (-1)^{|\omega|}\omega.\Bf\phi$$ for
$\omega\in\Aa(B), \phi \in\Aa(B,\pi_{*}^{N} (\Ef^0)\oplus
\pi_{*}^{N} (\Ef^0))$, with $\Bf_{[0]} = \Tf$.

Given $\Pp_i, \Pp_j \in \Gr(\pi_{*}^{N}(\Ef^0))$ we define the
even element
$$\Pp_{i,j} = \Pp_i \oplus \Pp_j \in  \Aa^0 (B, \End(\p_{*}^{N}
(\Ef^0\oplus \Ef^0))) \ .$$ We then have the graded Hermitian
subbundle $\Ww_i \oplus \Ww_j$ of $\pi_{*}^{N} (\Ef^0\oplus
\Ef^0)$ and canonical maps
$$\Aa (B, \pi_{*}^{N} (\Ef^0\oplus \Ef^0))\too  \Aa (B,
\Ww_i \oplus \Ww_j) \ , \hskip 10mm \omega\otimes s \mtoo
\omega\otimes \Pp_{i,j}s  \ ,$$ and
$$\Aa (B, \End(\pi_{*}^{N} (\Ef^0\oplus \Ef^0)))\too  \Aa (B,
\End(\Ww_i \oplus \Ww_j)) \ , \hskip 10mm \omega\otimes A \mtoo
\omega\otimes \Pp_{i,j} A \Pp_{i,j}  \ .$$ The latter extends to
define the induced superconnection on $\Ww_i \oplus \Ww_j$
$$\Bf_{i,j} : = \Pp_{i,j}\cdot \Bf \cdot \Pp_{i,j} :
\Aa(B, \Ww_i \oplus \Ww_j)\too \Aa(B, \Ww_i \oplus \Ww_j)  \ .$$
The degree 0 term is the odd-parity operator
\begin{equation}\label{e:Bij0}
\Bf^{i,j}_{[0]} = \Pp_{i,j}\, \Tf \, \Pp_{i,j} =
\begin{bmatrix}
  0 & \Pp_i \,\Tf_- \, \Pp_j \\
  \Pp_j\,\Tf_+ \, \Pp_i & 0 \\
\end{bmatrix} \ ,
\end{equation}
while the degree 1 term is the unitary connection of even parity
$\Bf^{i,j}_{[1]} = \Pp_{i,j}\cdot \Bf_{[1]} \cdot \Pp_{i,j} $ on
$\Ww_i \oplus \Ww_j$.

We shall consider the asymptotic behaviour of the scaled
superconnections
$$ \Bf_t = t^{1/2} \d_t \cdot\Bf \cdot\d_t\ii \ , \hskip 15mm
\Bf^{i,j}_t = t^{1/2} \d_t \cdot\Bf^{i,j}\cdot \d_t\ii \ .$$

For our purposes here it will be enough to consider just the
canonical superconnection defined by setting $\Tf_{\pm} =
\If_{N},$ and $\Bf_{[1]} = \na^{\pi_{*}^N(\Ef^0)} \oplus
\na^{\pi_{*}^N(\Ef^0)},$ where $\If_N$ is the vertical identity
operator on $\pi_{*}^{N} (\Ef^0)$. Thus
\begin{equation*}
\Bf_t =
\begin{bmatrix}
  \na^{\pi_{*}^N(\Ef^0)}& 0\\
  0 & \na^{\pi_{*}^N(\Ef^0)}\\
\end{bmatrix} +
t^{1/2}
\begin{bmatrix}
  0 & \If_N\\
  \If_N & 0 \\
\end{bmatrix}
\end{equation*}
and with $\na^i = \Pp_i\cdot\na^{\pi_{*}^N(\Ef^0)}\cdot\Pp_i$
\begin{equation}\label{e:bftij}
\Bf^{i,j}_t =
\begin{bmatrix}
  \na^{i}& 0\\
  0 & \na^{j}\\
\end{bmatrix} +
t^{1/2} L_{i,j}
\end{equation}
is a superconnection on $\Ww_i \oplus \Ww_j$ adapted to the family
$t^{1/2}L_{i,j}$. The curvature form $\Fff^{i,j}_t =
(\Bf^{i,j}_t)^2 \in \Aa^2(B,\End(\Ww_i \oplus \Ww_j))$ is then
\begin{equation}\label{e:Fij}
\Fff^{i,j}_t  =  \Rr_{i,j} + t^{1/2}\nabla L_{i,j} + tL_{i,j}^2 =
 \begin{bmatrix}
  \Rr_i + t \Pp_i \Pp_j \Pp_i  &  t^{1/2}\nabla^{i,j} (\Pp_i \Pp_j )\\
  t^{1/2}\nabla^{j,i} (\Pp_j \Pp_i ) &  \Rr_j + t \Pp_j \Pp_i \Pp_j\\
\end{bmatrix} \ ,
\end{equation}
where $\Rr_i = (\na^i)^2, \, \Rr_{i,j} = \Rr_i\oplus \Rr_j$ are
the curvature 2-forms on $\Ww_i$, $\Ww_i\oplus\Ww_j$, and
$\nabla^{i,j}\a = \nabla^j\a - \a\nabla^i$ the induced connection
on $\Hom(\Ww_i,\Ww_j)$. Note with  $\Pp_i = \Pp(\Ds), \Pp_j = \Pp$
\eqref{e:bftij} and \eqref{e:Fij} become $$\begin{bmatrix}
  \na^{\Ww(\Ds)}& t^{1/2}\Sss(\Pp)^{*}\\
  t^{1/2}\Sss(\Pp) & \na^{\Ww}\\
\end{bmatrix} \ , \hskip 10mm
\begin{bmatrix}
  \Rr_{\Ww(\Ds)} + t \Sss(\Pp)^* \Sss(\Pp)  &  t^{1/2}\nabla \Sss(\Pp)^*\\
  t^{1/2}\nabla \Sss(\Pp) &  \Rr_{\Ww} + t \Sss(\Pp)\Sss(\Pp)^*\\
\end{bmatrix} \ . $$ Since $(\Fff^{i,j}_t - \la \If_{i,j})_{[0]} = tL_{i,j}^2 - \la
\If_{i,j},$ where $\If_{i,j} = \If_i \oplus \If_j$ is the identity
operator on on $\Aa(B,\Ww_i \oplus \Ww_j)$, and $tL_{i,j}^2 $ is
positive self-adjoint, then $\Fff^{i,j}_t$ has a resolvent for
$\la\in\Cf\backslash\Rf_+$
\begin{equation}\label{e:Fijresolvent}
(\Fff^{i,j}_t - \la \If_{i,j})\ii = \sum_{k=0}^{\dim B}(-1)^k
(tL_{i,j}^2 - \la \If_{i,j})\ii \left\{(\Rr_{i,j} + t^{1/2}\nabla
L_{i,j})(tL_{i,j}^2 - \la \If_{i,j})\ii\right\}^k
 \ ,
\end{equation}
using \eqref{e:resolvent}.
Since $L_{i,j}^2 $ is bounded in each $C^l$ norm we can expand
this further via the infinite Neumann series for $t/|\la| <
\|L_{i,j}^2 \|_{l}\ii$
\begin{equation}\label{e:Lijresolvent}
(tL_{i,j}^2 - \la \If_{i,j})\ii = -\sum_{l\geq 0} t^l \la^{-l-1}
L_{i,j}^{2l} \ .
\end{equation}
It follows from \eqref{e:Fijresolvent}, which shows $(\Fff^{i,j}_t
- \la \If_{i,j})\ii $ is (the restriction of) an unbounded
differential operator, and from \eqref{e:Lijresolvent}, that the
heat operator $e^{-\Fff^{i,j}_t}$ can be defined. This is not,
however, the object of interest--in particular it is not trace
class and so has no Chern character. Rather we consider the
relative Chern character, defined via the choice of an
intermediary Grassmann section, as follows.

Let $\Pp_1, \Pp_2, \Pp_3 \in \Gr(\pi_{*}^{N}(\Ef^0))$. Then the
subbundles
\begin{equation*}
\Ww_1 \oplus \Ww_2 \ \ \ \ {\rm  and} \ \ \ \ \Ww_2 \oplus \Ww_3
\end{equation*}
of $\pi_{*}^{N} (\Ef^0\oplus \Ef^0)$ are endowed with the
canonical superconnections
\begin{equation*}
\Bf^{1,2}_t = \na^{1,2} + t^{1/2} L_{1,2} \ \ \ \ {\rm  and} \ \ \
\ \Bf^{2,3}_t = \na^{2,3} + t^{1/2} L_{2,3}
\end{equation*}
with resolvent operators
\begin{equation*}
(\Fff^{1,2}_t - \la \If_{1,2})\ii\in \Aa(B,\End(\Ww_1 \oplus
\Ww_2)) \ \ \ \ {\rm  and} \ \ \ \ (\Fff^{2,3}_t - \la
\If_{2,3})\ii\in \Aa(B,\End(\Ww_2 \oplus \Ww_3)) \ .
\end{equation*}

We construct the relative heat operator on $\pi_{*}^{N}
(\Ef^0\oplus \Ef^0)$  by extending the resolvents by zero on the
bundles $\Ww_i\pp \oplus \Ww_j\pp$. By inclusion $(\Fff^{i,j}_t -
\la \If_{i,j})\ii$ is identified with the $\pdo$
$$(\Fff^{i,j}_t - \la \If_{i,j})\ii_{|\Ef^0} := \Pp_{i,j}\cdot
(\Fff^{i,j}_t - \la \If_{i,j})\ii \cdot \Pp_{i,j}
$$ on $\pi_{*}^{N} (\Ef^0\oplus \Ef^0)$ and we have
\begin{equation}\label{e:inclusion}
(\Fff^{i,j}_t - \la \If_{i,j})\ii_{|\Ef^0} =
 (\Fff^{i,j}_t - \la \If)\ii - (-\la)\ii \Pp_{i,j}\pp \ \ \ \
 {\rm in} \ \ \
\Aa(B,\Psi(\Ef^0\oplus \Ef^0))\ ,
\end{equation}
where on the right-side $\Fff^{i,j}_t$ is the operator
$\Pp_{i,j}\cdot \Fff^{i,j}_t\cdot \Pp_{i,j}$ and $\If$ is the
vertical identity operator on $\pi_{*}^{N} (\Ef^0\oplus \Ef^0)$.
It follows that the relative resolvent $$(\Fff^{1,2}_t - \la
\If)_{|\Ef^0}\ii - (\Fff^{2,3}_t - \la \If)_{|\Ef^0}\ii $$ is
$O(|\la|\ii)$ as $\la\to\o$ in all $C^l$ norms (see
\propref{p:converges} below) and hence that the relative heat
operator can be defined by
\begin{equation}\label{e:relboundaryheat}
e^{-\Fff^{1,2}_t} - e^{-\Fff^{2,3}_t} =
\frac{i}{2\pi}\int_{\Cc}e^{-\la}\left((\Fff^{1,2}_t - \la
\If)_{|\Ef^0}\ii - (\Fff^{2,3}_t - \la \If)_{|\Ef^0}\ii\right) \
d\la \ .
\end{equation}

Let $\Psi^{-\o}(\Ef^0\oplus \Ef^0)$ be the vertical bundle of
smoothing operators, as in \cite{BeGeVe92}, whose sections are
smooth families of smoothing operators.
\begin{prop}\label{p:converges}
\noi [1] For each $\la\in \Cf\backslash\Rf_+$
\begin{equation*}
(\Fff^{1,2}_t - \la \If)_{|\Ef^0}\ii - (\Fff^{2,3}_t - \la
\If)_{|\Ef^0}\ii \in \Aa(B,\Psi^{-\o}(\Ef^0\oplus \Ef^0)) \ .
\end{equation*}
\noi [2]  The integral \eqref{e:relboundaryheat} converges in the
$\Ci$-topology to an element
\begin{equation*}
e^{-\Fff^{1,2}_t} - e^{-\Fff^{2,3}_t} \in
\Aa(B,\Psi^{-\o}(\Ef^0\oplus \Ef^0)) \ .
\end{equation*}
\end{prop}
\begin{proof}
From \eqref{e:Fijresolvent} and \eqref{e:inclusion} we have
\begin{equation*}
(\Fff^{1,2}_t - \la \If)_{|\Ef^0}\ii - (\Fff^{2,3}_t - \la
\If)_{|\Ef^0}\ii =  (\Fff^{1,2}_t - \la \If)\ii - (\Fff^{2,3}_t -
\la \If)\ii + (-\la)\ii (\Pp_{2,3}\pp  - \Pp_{1,2}\pp)
\end{equation*}
\begin{equation*}
 = (-\la)\ii (\Pp_{1,2}  - \Pp_{2,3}) +
 \sum_{k=0}^{\dim B}(-1)^k (tL_{1,2}^2 - \la \If)\ii
\left\{(\Rr_{1,2} + t^{1/2}\nabla L_{1,2})(tL_{1,2}^2 - \la
\If)\ii\right\}^k
\end{equation*}
$$- (-1)^k (tL_{2,3}^2 - \la \If)\ii
\left\{(\Rr_{2,3} + t^{1/2}\nabla L_{2,3})(tL_{2,3}^2 - \la
\If)\ii\right\}^k$$

$$ = (-\la)\ii (\Pp_{1,2}  - \Pp_{2,3})$$
\begin{equation*}
  + \sum_{k=0}^{\dim B}(-1)^k \sum_{i=0}^k  \left\{(tL_{2,3}^2 - \la \If)\ii
(\Rr_{2,3} + t^{1/2}\nabla L_{2,3})\right\}^i \left\{(tL_{1,2}^2 -
\la \If)\ii - (tL_{2,3}^2 - \la \If)\ii\right\}
\end{equation*}
$$\times \left\{(\Rr_{1,2}
+ t^{1/2}\nabla L_{1,2})(tL_{1,2}^2 - \la \If)\ii \right\}^{k-i}$$
\begin{equation}\label{e:relresboundary}
  + \sum_{k=1}^{\dim B}(-1)^k \sum_{j=0}^{k-1}  \left\{(tL_{2,3}^2 - \la \If)\ii
(\Rr_{2,3} + t^{1/2}\nabla L_{2,3})\right\}^j (tL_{2,3}^2 - \la
\If)\ii
\end{equation}
\begin{equation*}
\times \left\{(\Rr_{1,2} - \Rr_{2,3}) + t^{1/2}(\nabla L_{1,2} -
\nabla L_{2,3})\right\}\left\{(tL_{1,2}^2 - \la \If)\ii (\Rr_{1,2}
+ t^{1/2}\nabla L_{1,2})(tL_{1,2}^2 - \la \If)\ii \right\}^{k-j-1}
\ .
\end{equation*}

\noi From \eqref{e:relresboundary} and \eqref{e:Lijresolvent} we
hence obtain for any vector fields $\xi_1,\ldots,\xi_l$ on $B$
$$
\|\na^{\pi_*(\Ef^0\oplus\Ef^0)}_{\xi_1}\ldots\na^{\pi_*(\Ef^0\oplus\Ef^0)}_{\xi_l}
\left((\Fff^{1,2}_t - \la \If)_{|\Ef^0}\ii - (\Fff^{2,3}_t - \la
\If)_{|\Ef^0}\ii\right)\| \leq C_l t^l |\la|^{-l-1}$$ and so
\eqref{e:relboundaryheat} is convergent in each $C^l$ norm.

The vertical operator $D_{i,j} = \Rr_{i,j} + t^{1/2}\nabla
L_{i,j}$ is a differential form with differential operator
coefficients. Since these are families of operators over a closed
manifold we have by \cite{BeGeVe92} Cor(2.40) that if $K\in
\Aa(B,\Psi^{-\o}(\Ef^0\oplus \Ef^0))$, then $D_{i,j}K$ and
$KD_{i,j}$ are also smooth families of smoothing operators. Hence
since $(tL_{i,j}^2 - \la \If)\ii$ is a vertical $\pdo$ of order
$0$, then [1] is immediate from \eqref{e:relresboundary} once we
have the following sublemma.
\begin{lem}\label{lem:relativelyboundarysmooth}
\begin{equation}\label{e:bsmooth0}
\Pp_{1,2} - \Pp_{2,3}  \in \Aa(B,\Psi^{-\o}(\Ef^0\oplus \Ef^0)) \
.
\end{equation}
\begin{equation}\label{e:bsmooth2}
(tL_{1,2}^2 - \la \If)\ii - (tL_{2,3}^2 - \la \If)\ii\in
\Aa(B,\Psi^{-\o}(\Ef^0\oplus \Ef^0)) \ .
\end{equation}
\begin{equation}\label{e:bsmooth3}
\Rr_{1,2}^{k} - \Rr_{2,3}^{k}\in \Aa(B,\Psi^{-\o}(\Ef^0\oplus
\Ef^0)) \ .
\end{equation}
\begin{equation}\label{e:bsmooth4}
\nabla L_{1,2}-\nabla L_{2,3}\in \Aa(B,\Psi^{-\o}(\Ef^0\oplus
\Ef^0)) \ .
\end{equation}
\end{lem}
\begin{proof}
\eqref{e:bsmooth0} is obvious. For \eqref{e:bsmooth2} we have
\begin{equation*}
(tL_{1,2}^2 - \la \If)\ii - (tL_{2,3}^2 - \la \If)\ii = t
(tL_{1,2}^2 - \la \If)\ii(L_{2,3}^2 - L_{1,2}^2)(tL_{2,3}^2 - \la
\If)\ii
\end{equation*}
while $(tL_{1,2}^2 - \la \If)\ii\in \Aa(B,\Psi^0(\Ef^0\oplus
\Ef^0))$ and
\begin{equation}\label{e:relLsmooth}
L_{2,3}^2 - L_{1,2}^2 = (\Pp_2\Pp_3\Pp_2 - \Pp_1\Pp_2\Pp_1)\oplus
(\Pp_3\Pp_2\Pp_3 - \Pp_2\Pp_1\Pp_2)\in
\Aa(B,\Psi^{-\o}(\Ef^0\oplus \Ef^0))
\end{equation}
since $\Pp_i - \Pp_j\in \Aa(B,\Psi^{-\o}(\Ef^0))$.

It is enough to prove \eqref{e:bsmooth3} over a  local
trivialization $N_{|U} = Y_{z_0} \times U$ of the boundary
fibration in a neighborhood $U$ of $z_0 \in B$. This defines an
identification $\Ef^0_{|U} = E^0_{z_0} \times U$, where
$E^0_{z_0}$ is the superbundle over $Y_{z_0} = \dd X_{z_0}$, and
hence a trivialization (in the weak sense) of $\pi^{N}_*
(\Ef^0\oplus\Ef^0)$ over $U$ as $ \pi^{N}_*
(\Ef^0\oplus\Ef^0)_{|U} = \G(Y_{z_0},E^0_{z_0}\oplus E^0_{z_0})
\times U$ and
$$\Aa(U,\pi^{N}_* (\Ef^0\oplus\Ef^0))
= \Aa(U)\otimes \G(Y_{z_0},E^0_{z_0}\oplus E^0_{z_0}) \ .$$ Any
superconnection $\Bf$ on $\pi^{N}_* (\Ef^0\oplus\Ef^0)$ then takes
the local form
\begin{equation}\label{e:Bflocal}
\Bf = d_U + \sum_I D_{I}^z dz_I \ ,
\end{equation}
where $dz_I \in \Aa^{|I|}(U)$, and $D_I^z $ are differential
operators on $\G(Y_{z_0},E^0_{z_0}\oplus E^0_{z_0})$ depending
smoothly on $z$. We have over $U$
\begin{equation}\label{e:APlocalb}
\Bf^{i,j}_{|U}  = P^{z}_{i,j}\cdot d_U \cdot P^{z}_{i,j} +
\sum_{I} P^{z}_{i,j}\cdot D_{I}^z \cdot P^{z}_{i,j} dz_I \ ,
\end{equation}
where the $P^{z}_{i,j}$ are elements of the fixed Grassmannian
$\Gr_{\o}(\dd_{z_0}\oplus\dd_{z_0})$ depending smoothly on $z$ and
defining $P^{z}_{i,j} : U \too \Aa(U)\otimes
\G(Y_{z_0},\End(E^0_{z_0}\oplus E^0_{z_0}))$. Fixing
$P_0\in\Gr_{\o}(\dd_{z_0}\oplus\dd_{z_0})$ we have $P^{z}_{i,j} -
P_0 \in \Aa^0(U,\Psi^{-\o}(E_{z_0}^0\oplus E_{z_0}^0))$ and hence
that
\begin{equation}\label{e:dPsmooth}
d P^{z}_{i,j} = d (P^{z}_{i,j} - P_0)  \in
\Aa^1(U,\Psi^{-\o}(E_{z_0}^0\oplus E_{z_0}^0)) \ ,
\end{equation}
is a smooth family of smoothing operators.

From \eqref{e:APlocalb}
\begin{eqnarray}
& &  \hskip 5mm (\Rr_{i,j})_{|U} =  P^{z}_{i,j} dP^{z}_{i,j}
\wedge dP^{z}_{i,j} + \sum_I
(P^{z}_{i,j}\cdot D_I^z \cdot P^{z}_{i,j}) dz_I\wedge dP^{z}_{i,j}  \label{e:curvlocala}\\
& + & P^{z}_{i,j}\sum_I\sum_{l=0}^{\dim B} \frac{\dd}{\dd
z_l}(P^{z}_{i,j}\cdot D_I^z \cdot P^{z}_{i,j} ) dz_l\wedge dz_I +
\left(\sum_I (P^{z}_{i,j}\cdot D_I^z \cdot P^{z}_{i,j})
dz_I\right)^2 \ . \label{e:curvlocalN}
\end{eqnarray}
From \eqref{e:dPsmooth} the first two terms on the right-side are
in $\Aa(U,\Psi^{-\o}(E_{z_0}^0\oplus E_{z_0}^0))$. The  terms in
\eqref{e:curvlocalN} are non-smoothing, but their difference in
$\Rr_{1,2} - \Rr_{2,3}$ are smooth families of smoothing
operators. For the first term in \eqref{e:curvlocalN} we have
\begin{eqnarray*}
& & \sum_I\sum_{l=0}^{\dim B} \left(P^{z}_1 \frac{\dd D_I^z }{\dd
z_l} P^{z}_1  - P^{z}_2 \frac{\dd D_I^z }{\dd z_l} P^{z}_2\right)
dz_l\wedge dz_I  + F(\Pp_{1,2},\Pp_{2,3}) \\
& = & \sum_I\sum_{l=0}^{\dim B} (P^{z}_{1,2} - P^{z}_{2,3})
\frac{\dd D_I^z }{\dd z_l} P^{z}_{1,2}  - P^{z}_{2,3} \frac{\dd
D_I^z }{\dd z_l} (P^{z}_{1,2}- P^{z}_{2,3}) dz_l\wedge dz_I +
F(\Pp_{1,2},\Pp_{2,3}) \ ,
\end{eqnarray*}
where $F(\Pp_{1,2},\Pp_{2,3})$ is a sum of  terms involving
derivatives of $\Pp_{i,j}$ and hence smoothing. Since the
composition of a differential operator and a smoothing operator is
smoothing, this is in $\Aa(U,\Psi^{-\o}(E_{z_0}^0\oplus
E_{z_0}^0))$. Similarly the second term in \eqref{e:curvlocalN} is
in $\Aa(U,\Psi^{-\o}(E_{z_0}^0\oplus E_{z_0}^0))$. This proves
\eqref{e:bsmooth3} for $k=1$. The general case now follows from
the identity
\begin{equation}\label{e:R1k-R2k}
\Rr_{1,2}^{k} - \Rr_{2,3}^{k} = \sum_{j=0}^{k-1} \Rr_{2,3}^j \,
(\Rr_{1,2} - \Rr_{2,3}) \Rr_{1,2}^{k-j-1} \ .
\end{equation}
The final assertion \eqref{e:bsmooth4} follows by a similar
argument.
\end{proof}

Since smoothing operators form an ideal in the bundle of $\pdo$s
$\Psi(\Ef^0\oplus \Ef^0)$, \eqref{e:relresboundary} and the lemma
 therefore imply [1]. Equivalently,
$(\Fff^{1,2}_t - \la \If)_{|\Ef^0}\ii - (\Fff^{2,3}_t - \la
\If)_{|\Ef^0}\ii $ is a vertical integral operator with $\Ci$
kernel $K_{t,\la}^{1,2,3}\in\Aa(N\times_{\pi^N} N,
(\Ef^0\oplus\Ef^0)\boxtimes (\Ef^0\oplus\Ef^0)^*)$. Hence
$e^{-\Fff^{1,2}_t} - e^{-\Fff^{2,3}_t}$ has
 kernel $$E_{t}^{1,2,3}(x,y) = \frac{i}{2\pi}\int_{\Cc} e^{-\la}
 K_{t,\la}^{1,2,3}(x,y) \, d\la$$ which by a standard argument
 is therefore also in $\Aa(N\times_{\pi^N} N, (\Ef^0\oplus\Ef^0)\boxtimes
(\Ef^0\oplus\Ef^0)^*)$. This completes
 the proof.
\end{proof}

\vskip 2mm

\subsection{Chern Character Forms and Transgression}
Defined on $\Aa(B,\Psi^{-\o}(\Ef^0\oplus \Ef^0))$ is the fibrewise
supertrace
\begin{equation}\label{e:Nstr}
\str : \Aa(B,\Psi^{-\o}(\Ef^0\oplus \Ef^0)) \too \Aa(B) \ , \hskip
10mm \str(A)(z) = \int_{N_z} \str(a_z(y,y)) \ dy \ ,
\end{equation}
where $A$ has fibrewise kernel $a_z\in \Aa(N_z\times_{\pi^N} N_z,
(E_z^0\oplus E_z^0)\boxtimes (E_z^0\oplus E_z^0)^*)$, and hence by
\propref{p:converges} we have the relative Chern character form
associated to the superconnections $\Bf^{1,2}_t, \Bf^{2,3}_t$ as
the even degree element of $\Aa(B) $ defined by
\begin{equation}\label{e:Nch}
    \ch(\Bf^{1,2}_t, \Bf^{2,3}_t) = \str(e^{-\Fff^{1,2}_t} -
    e^{-\Fff^{2,3}_t}) \ .
\end{equation}
By \propref{p:converges}[1]
\begin{equation}\label{e:Nch1}
    \ch(\Bf^{1,2}_t, \Bf^{2,3}_t) = \frac{i}{2\pi}\int_{\Cc} e^{-\la}
\,\str\left((\Fff^{1,2}_t - \la \If)_{|\Ef^0}\ii - (\Fff^{2,3}_t -
\la \If)_{|\Ef^0}\ii\right)\, d\la \ .
\end{equation}
The relative Chern character form defines a deRham cohomology
class in $H^{2\bullet}(B)$:
\begin{prop}\label{p:chclosed}
The differential form $\ch(\Bf^{1,2}_t, \Bf^{2,3}_t)$ is closed.
\end{prop}
\begin{proof}
We require the following lemma.
\begin{lem}\label{lem:dstr}
Let $\Pp_{\a},\Pp_{\b}\in\Gr(\pi_{*}^N(\Ef^0\oplus\Ef^0))$ be of
even parity, and let $$a_{\a} = \Pp_{\a} a_{\a} \Pp_{\a}, \ a_{\b}
= \Pp_{\b} a_{\b} \Pp_{\b} \in
\Aa(B,\End(\pi_{*}^N(\Ef^0\oplus\Ef^0)))$$ be even parity such
that $a_{\a}- a_{\b} \in \Aa(B,\Psi(\Ef^0\oplus\Ef^0))$. Let $\Bf$
be a superconnection on $\pi_{*}^N(\Ef^0\oplus\Ef^0)$ and let
$\Bf_{\a} = \Pp_{\a}\cdot \Bf\cdot \Pp_{\a}, \Bf_{\b}=
\Pp_{\b}\cdot \Bf\cdot \Pp_{\b}$. Then in $\Aa(B)$
\begin{equation}\label{e:Ndstr}
d_{B} \cdot \str(a_{\a}- a_{\b}) =
    \str(\Bf_{\a} a_{\a} - \Bf_{\b} a_{\b}) \ ,
\end{equation}
where $\Bf_{\a} a_{\a} = [\Bf_{\a} ,a_{\a}]$, $\Bf_{\b} a_{\b} =
 [\Bf_{\b} ,a_{\b}]$.
\end{lem}
\begin{proof}
$a_{\a}- a_{\b}  \in \Aa(B,\Psi^{-\o}(\Ef^0\oplus\Ef^0))$ is a
consequence of $\Pp_{\a}- \Pp_{\b}\in
\Aa(B,\Psi^{-\o}(\Ef^0\oplus\Ef^0))$. Hence from \cite{BeGeVe92}
Lem(9.14)
\begin{equation*}
 d_{B} \cdot \str(a_{\a}- a_{\b} ) = \str(\Bf a_{\a} - \Bf a_{\b}) \ .
\end{equation*}
For $\sigma_1, \sigma_2 \in
\Aa(B,\End(\pi_{*}^N(\Ef^0\oplus\Ef^0)))$ one has $\Bf (\sigma_1
\sigma_2) = (\Bf \sigma_1)\cdot \sigma_2 + (-1)^{|\sigma_1 |}
\sigma_1\cdot \Bf(\sigma_2).$ Hence $\Bf a_{\a} = \Bf
(\Pp_{\a}\cdot a_{\a})  = (\Bf \Pp_{\a})\cdot a_{\a} + \Pp_{\a}\Bf
a_{\a}  =   (\Bf \Pp_{\a})\cdot a_{\a} + \Bf_{\a} a_{\a}$ since
$|\Pp_{\a}| = 0$. Thus
\begin{equation}\label{e:str00}
\str(\Bf a_{\a} - \Bf a_{\b}) = \str((\Bf \Pp_{\a})\cdot a_{\a} -
(\Bf \Pp_{\b})\cdot a_{\b}) + \str(\Bf_{\a} a_{\a} - \Bf_{\b}
a_{\b}) \ ,
\end{equation}
since $(\Bf \Pp_{\a})\cdot a_{\a} - (\Bf \Pp_{\b})\cdot a_{\b}$ is
clearly a smooth family of smoothing operators.

The proof is thus to show that the first term on the right-side of
\eqref{e:str00} vanishes. We may work locally over an open set $U$
in $B$, where $\Bf = d_U + \omega$ and $\omega = \sum_I D_I dz_I \
.$ The term in question is
\begin{equation*}
\str((d\Pp_{\a})\cdot a_{\a}- (d\Pp_{\a})\cdot a_{\b} ) +
\str([\omega, \Pp_{\a}]\cdot a_{\a} - [\omega, \Pp_{\b}]\cdot
a_{\b})
\end{equation*}
\begin{equation}\label{e:strlocal}
= \str((d\Pp_{\a})\cdot a_{\a}\cdot \Pp_{\a}) -
\str((d\Pp_{\a})\cdot a_{\b}\cdot \Pp_{\b} ) + \str(\Pp_{\a}\pp
\omega a_{\a} - \Pp_{\b}\pp \omega a_{\b}) \ .
\end{equation}
Here we use the fact that locally the operators $(\Pp_{\a})_{|U}
=(P_{\a,z} \ | \ z\in U)$ can be taken in a fixed Grassmannian
$\Gr_{\o}(\dd_{z_0}\oplus\dd_{z_0})$ and hence relative to the
trivialization $d\Pp_{\a}\in \Aa(U,\Psi^{-\o}(\Ef^0\oplus\Ef^0))$.
We have
\begin{equation*}
 \str((d\Pp_{\a})\cdot a_{\a}\cdot \Pp_{\a}) =
  \str((d\Pp_{\a})\cdot \Pp_{\a} a_{\a}\cdot \Pp_{\a}^2) =
   \str((\Pp_{\a}(d\Pp_{\a}) \Pp_{\a}) a_{\a}\cdot \Pp_{\a}) = 0
   \ ,
\end{equation*}
since $\Pp_{\a}^2 = \Pp_{\a}$ implies $\Pp_{\a}(d\Pp_{\a})
\Pp_{\a} = 0$. Since $\Pp_{\a}\pp -\Pp_{\b}\pp, a_{\a}- a_{\b} \in
\Aa(B,\Psi^{-\o}(\Ef^0\oplus\Ef^0))$, the remaining term in
\eqref{e:strlocal} is
\begin{eqnarray*}
\str(\Pp_{\a}\pp \omega a_{\a} - \Pp_{\b}\pp \omega a_{\b}) & = &
\str((\Pp_{\a}\pp -\Pp_{\b}\pp)\omega \Pp_{\a} a_{\a} \Pp_{\a}) +
\str(\Pp_{\b}\pp \omega (\Pp_{\a}a_{\a}\Pp_{\a} -
\Pp_{\b}a_{\b}\Pp_{\b})) \\
& = & - \str(\Pp_{\a}\Pp_{\b}\pp\omega \Pp_{\a} a_{\a} \Pp_{\a}) +
\str(\Pp_{\b}\pp\omega \Pp_{\a} a_{\a} \Pp_{\a}\Pp_{\b}\pp) \\
& = & \str([\Pp_{\b}\pp\omega \Pp_{\a} a_{\a}
\Pp_{\a},\Pp_{\a}\Pp_{\b}\pp]) = 0  \ ,
\end{eqnarray*}
where the second equality follows by cycling the trace, the third
since $\Pp_{\b}\pp\omega \Pp_{\a} a_{\a} \Pp_{\a}$ and
$\Pp_{\a}\Pp_{\b}\pp$ have odd and even parity respectively, and
the final equality since $\str$ vanishes on supercommutators.
\end{proof}
Applying \lemref{lem:dstr} with $\Pp_{\a} = \Pp_{1,2}, \Pp_{\b} =
\Pp_{2,3}$, so that $\Bf_{\a} = \Bf_{t}^{1,2}, \Bf_{\b} =
\Bf_{t}^{2,3}$,  to \eqref{e:Nch1} we have
\begin{equation*}
    d_B \ch(\Bf^{1,2}_t, \Bf^{2,3}_t)  =  \frac{i}{2\pi}\int_{\Cc} e^{-\la}
\str\left([\Bf_{t}^{1,2}, (\Fff^{1,2}_t - \la \If)_{|\Ef^0}\ii ] -
[\Bf_{t}^{2,3}, (\Fff^{2,3}_t - \la \If)_{|\Ef^0}\ii ]\right)\, d\la \\
\end{equation*}
and this vanishes since $[\Bf_{t}^{i,j}, (\Fff^{i,j}_t - \la
\If)_{|\Ef^0}\ii ] = [\Bf_{t}^{i,j}, ((\Bf_{t}^{i,j})^2 - \la
\If)_{|\Ef^0}\ii ]$ vanishes.
\end{proof}
Next, we have the transgression formula:
\begin{prop}\label{p:Ntrans}
One has
\begin{equation}\label{e:Ntrans}
\frac{d}{dt} \ch(\Bf^{1,2}_t, \Bf^{2,3}_t) = - d_B
\str(\dot{\Bf}_{t}^{1,2} e^{-\Fff^{1,2}_t} - \dot{\Bf}_{t}^{2,3}
e^{-\Fff^{2,3}_t}) \ .
\end{equation}
Equivalently, for $0 < t < T < \o$
\begin{equation}\label{e:Ntrans1}
\ch(\Bf^{1,2}_T, \Bf^{2,3}_T)  -  \ch(\Bf^{1,2}_t, \Bf^{2,3}_t) =
- d_B \int_{t}^T\str(\dot{\Bf}_{t}^{1,2} e^{-\Fff^{1,2}_t} -
\dot{\Bf}_{t}^{2,3} e^{-\Fff^{2,3}_t}) \, dt \ .
\end{equation}
\end{prop}
\begin{proof}
From \eqref{e:inclusion}
\begin{equation*}
\dd_t \ch(\Bf^{1,2}_t, \Bf^{2,3}_t) = \frac{i}{2\pi}\int_{\Cc}
e^{-\la} \, \str\left(\dd_t(\Fff^{1,2}_t - \la \If)\ii -
\dd_t(\Fff^{2,3}_t - \la \If)\ii\right)\, d\la \ ,
\end{equation*}
while since $\Bf^{i,j}_t$ commutes with  $(\Fff^{i,j}_t - \la
\If)\ii$
\begin{eqnarray*}
\dd_t(\Fff^{i,j}_t - \la \If)\ii & = &  - (\Fff^{i,j}_t - \la
\If)\ii \cdot (\Bf^{i,j}_t \dot{\Bf}_{t}^{i,j}) \cdot(\Fff^{i,j}_t
- \la \If)\ii - (\Fff^{i,j}_t - \la \If)\ii\cdot
(\dot{\Bf}_{t}^{i,j}\Bf^{i,j}_t )\cdot (\Fff^{i,j}_t - \la \If)\ii
\\ & = & -[\Bf^{i,j}_t,  (\Fff^{i,j}_t - \la \If)\ii\cdot
\dot{\Bf}_{t}^{i,j}\cdot (\Fff^{i,j}_t - \la \If)\ii ] \ .
\end{eqnarray*}
By \lemref{lem:dstr} we conclude
\begin{equation*}
\dd_t \ch(\Bf^{1,2}_t, \Bf^{2,3}_t)  = - d_B \tau_t \ ,
\end{equation*}
where
\begin{equation*}
\tau_t = \frac{i}{2\pi}\int_{\Cc} e^{-\la} \, \str\left(
(\Fff^{1,2}_t - \la \If)\ii  \dot{\Bf}_{t}^{1,2} (\Fff^{1,2}_t -
\la \If)\ii - (\Fff^{2,3}_t - \la \If)\ii  \dot{\Bf}_{t}^{2,3}
(\Fff^{2,3}_t - \la \If)\ii\right)\, d\la \ .
\end{equation*}
Writing the integrand as
\begin{equation*}
\str\left( \left((\Fff^{1,2}_t - \la \If)\ii  - (\Fff^{2,3}_t -
\la \If)\ii\right)\dot{\Bf}_{t}^{1,2} (\Fff^{1,2}_t - \la \If)\ii
+ (\Fff^{1,2}_t - \la \If)\ii \left(\dot{\Bf}_{t}^{1,2}-
\dot{\Bf}_{t}^{2,3}\right) (\Fff^{2,3}_t - \la \If)\ii\right.
\end{equation*}
$$\left.+
(\Fff^{1,2}_t - \la \If)\ii \dot{\Bf}_{t}^{2,3}\left((\Fff^{1,2}_t
- \la \If)\ii  - (\Fff^{2,3}_t - \la \If)\ii\right)\right)$$ each
of the differences is an element of
$\Aa(B,\Psi^{-\o}(\Ef^0\oplus\Ef^0))$ and hence we can cycle the
trace on each term to get
\begin{equation*}
\tau_t = \frac{i}{2\pi}\int_{\Cc} e^{-\la} \, \str\left(
\dot{\Bf}_{t}^{1,2} (\Fff^{1,2}_t - \la \If)^{-2} -
\dot{\Bf}_{t}^{2,3} (\Fff^{2,3}_t - \la \If)^{-2}\right)\, d\la \
,
\end{equation*}
and integrating by parts we obtain \eqref{e:Ntrans}.
\end{proof}

\vskip 2mm

\subsection{Computing {\bf $\lim_{t\to 0} \ch(\Bf^{1,2}_t,
\Bf^{2,3}_t) $} }

The curvature $\Rr_{i}$ of the bundle $\Ww_i$ is a pure 2-form
with restricted differential operator coefficients. It is
therefore nilpotent. The Chern character form
\begin{equation}\label{e:eR}
e^{-\Rr_i}  =  \frac{i}{2\pi}\int_{c_0}e^{-\la} \left((\Rr_i -
\la\If)\ii  - \Pp_i\pp(-\la)\ii \right) \ d\la \ ,
\end{equation}
where $c_0$ is any simple closed contour around the origin, is
therefore well-defined and is equal to the sum with only finitely
many non-zero terms $ \Pp_i + \sum_{k\geq
1}(-1)^k\frac{1}{k!}\,\Rr_i^k$.

For two choices of Grassmann section $\Pp_i ,\Pp_j $ we have
$e^{-\Rr_i} - e^{-\Rr_j} \in \Aa(B,\Psi^{-\o}(\Ef^0))$ by the same
argument as \lemref{lem:relativelyboundarysmooth}. We therefore
have the well-defined relative Eta-form
\begin{equation}\label{e:genetaform1}
    \eta(\Pp_i ,\Pp_j ) = \Tr(e^{-\Rr_i} - e^{-\Rr_j}) \in \Aa^{2\bullet}(B) \
\end{equation}
and since $\frac{i}{2\pi}\int_{c_0} e^{-\la} \, (-\la)^{-k-1} \
d\la = (k!)\ii$,
\begin{equation}\label{e:genetaformseries1}
\eta(\Pp_i ,\Pp_j )  = \Tr(\Pp_i - \Pp_j) +
\sum_{k=1}^{\o}(-1)^k\frac{1}{k!}\,\Tr\left(\Rr_i^k -
\Rr_i^k\right) \ .
\end{equation}
In particular, the degree 0 term is the pointwise relative index:
\begin{equation}\label{e:etadeg0}
\eta(\Pp_i ,\Pp_j )_{[0]}  = \Tr(\Pp_i - \Pp_j) = \ind(D_{P_j}) -
\ind(D_{P_i}) \ ,
\end{equation}
where $D_{P_i}$ is any global boundary problem in the family
$\Ds_{\Pp_i}$.\\[1mm]

\noi We have:

\begin{prop}\label{p:relindexthm1}
The differential form $\ch(\Bf^{1,2}_t, \Bf^{2,3}_t)$ has a limit
as $t\too 0$ in $\Aa(B)$ given by
\begin{equation}\label{e:tlimitatzero1}
\lim_{t\to 0} \ch(\Bf^{1,2}_t, \Bf^{2,3}_t)  = \eta(\Pp_1 ,\Pp_3 )
\ .
\end{equation}
\end{prop}
\begin{proof}
From \eqref{e:relresboundary} we have
$$\ch(\Bf^{1,2}_t, \Bf^{2,3}_t)$$
\begin{equation}\label{e:ch1}
 = \frac{i}{2\pi}\int_{\Cc}e^{-\la}
(-\la)\ii \ d\la \ . \ \str(\Pp_{1,2}  - \Pp_{2,3})
\end{equation}
\begin{equation}\label{e:ch2}
  + \sum_{k=0}^{\dim B}(-1)^k \sum_{i=0}^k  \frac{i}{2\pi}\int_{\Cc}e^{-\la}
  \ \str\left[\left\{(tL_{2,3}^2 - \la \If)\ii
(\Rr_{2,3} + t^{1/2}\nabla L_{2,3})\right\}^i \right.
\end{equation}
$$\times \left.\left((tL_{1,2}^2 -
\la \If)\ii - (tL_{2,3}^2 - \la \If)\ii\right)\left\{(\Rr_{1,2} +
t^{1/2}\nabla L_{1,2})(tL_{1,2}^2 - \la \If)\ii
\right\}^{k-i}\right] \ d\la $$
\begin{equation}\label{e:ch3}
  + \sum_{k=1}^{\dim B}(-1)^k \sum_{j=0}^{k-1}  \frac{i}{2\pi}\int_{\Cc}e^{-\la}
  \ \str \left[\left\{(tL_{2,3}^2 - \la \If)\ii
(\Rr_{2,3} + t^{1/2}\nabla L_{2,3})\right\}^j (tL_{2,3}^2 - \la
\If)\ii\right.
\end{equation}
$$\times \left((\Rr_{1,2} - \Rr_{2,3}) + t^{1/2}(\nabla L_{1,2} -
\nabla L_{2,3})\right)$$
$$
\times \left.\left\{(tL_{1,2}^2 - \la \If)\ii (\Rr_{1,2} +
t^{1/2}\nabla L_{1,2})(tL_{1,2}^2 - \la \If)\ii
\right\}^{k-j-1}\right] \ d\la \ . $$ Clearly, the first term
\eqref{e:ch1} is the degree zero term \eqref{e:etadeg0}. Next,
\eqref{e:relLsmooth} implies
$$(tL_{1,2}^2 - \la \If)\ii
-(tL_{2,3}^2 - \la \If)\ii = t(tL_{1,2}^2 - \la \If)\ii(L_{2,3}^2
- L_{1,2}^2)(tL_{2,3}^2 - \la \If)\ii$$ is $O(|\la|^{-2}).t$ in
any $C^l$ norm as $t\too 0$. It follows using
\eqref{e:Lijresolvent} that the second term \eqref{e:ch2} is
$o(1)$ as $t\too 0 $ and hence makes no contribution to
\eqref{e:tlimitatzero1}.

Similarly, asymptotically as $t\too 0$ the third term
\eqref{e:ch3} is
\begin{equation*}
\sum_{k\geq 1}\sum_{j=0}^{k-1} (-1)^k
\frac{i}{2\pi}\int_{\Cc}e^{-\la}
  \ \left( \str (\Rr^j_{2,3}(\Rr_{1,2} -
  \Rr_{2,3})\Rr_{1,2}^{k-j-1})
  (-\la)^{-k-1} t^0 + O(|\la|^{-k-2} ).t\right) \ d\la
\end{equation*}
\begin{equation*}
= \sum_{k\geq 1}\sum_{j=0}^{k-1} \frac{(-1)^k}{k!}
  \ \str (\Rr^j_{2,3}(\Rr_{1,2} - \Rr_{2,3})\Rr_{1,2}^{k-j-1})
   + o(1)
\end{equation*}
\begin{equation*}
= \sum_{k\geq 1} \frac{(-1)^k}{k!}
  \ \Tr (\Rr^k_1 - \Rr^k_3)  + o(1 ) \ .
\end{equation*}
This completes the proof.
\end{proof}

One can prove a somewhat stronger property of the transgression
form
\begin{equation*}
\tau(\Bf_{t}^{1,2},\Bf_{t}^{2,3}) = \str(\dot{\Bf}_{t}^{1,2}
e^{-\Fff^{1,2}_t} - \dot{\Bf}_{t}^{2,3} e^{-\Fff^{2,3}_t}) \ .
\end{equation*}
Namely:
\begin{prop}\label{e:tform0}
As $t\too 0$ there is an asymptotic expansion
\begin{equation}\label{e:tformasymp}
\tau(\Bf_{t}^{1,2},\Bf_{t}^{2,3}) \sim \sum_{j\geq 0} \tau_j
t^{j-1/2} \ ,
\end{equation}
where $\tau_j\in\Aa(B)$.
\end{prop}

It follows that the integral in \eqref{e:Ntrans1} converges as
$t\too 0$. Combined with \propref{p:relindexthm1} this implies:
\begin{cor}
For any $t > 0$ one has in $\Aa(B)$
\begin{equation*}
\ch(\Bf^{1,2}_t, \Bf^{2,3}_t) = \eta(\Pp_1 ,\Pp_3 ) - d_B\Tt_t \ ,
\end{equation*}
where $\Tt_t = \lim_{s\to 0}\int_{s}^t
\tau(\Bf_{s}^{1,2},\Bf_{s}^{2,3}) \ ds$.
\end{cor}

\vskip 2mm

\subsection{Computing {\bf $\lim_{t\to \o} \ch(\Bf^{1,2}_t,
\Bf^{2,3}_t) $} }

To  complete the proof of Theorem(I) first observe by the
linearity of the supertrace
\begin{equation*}
\eta(\Pp_1 ,\Pp_2) + \eta(\Pp_2 ,\Pp_3) = \eta(\Pp_1 ,\Pp_3)\ .
\end{equation*}
Hence from \propref{p:indexbundle2}, we see that \eqref{e:tI1}
follows from \eqref{e:tI2} on setting $\Pp_3 = \Pp(\Ds), \Pp_1 =
\Pp$, since $\ch(\Ind(\Ds_{\Pp(\Ds)})) = 0$ in cohomology. From
the transgression formula \eqref{e:Ntrans1} and and the small $t$
limit \eqref{e:tlimitatzero1}, and by \propref{p:indexbundle2},
the remaining component in the proof is the following
identification.
\begin{prop}
The cohomology class of $\ch(\Bf^{1,2}_t, \Bf^{2,3}_t)$ is equal
to the Chern character of the index bundle $\Ind(\Pp_3,\Pp_1)$.
\end{prop}
\begin{proof}
Broadly the proof follows that of \cite{BeGeVe92}, but generalized
to the relative context. We begin by assuming that the dimensions
of the kernels of the operators $L_{i,j}^z$ on the fibres
$W_{i,z}\oplus W_{j,z}$ of the bundles $\Ww_i\oplus\Ww_j$ are
constant as $z$ varies in $B$. Let $\Pi^z_{i,j}$ be the orthogonal
projection onto $\Ker(L_{i,j}^z)$ considered as a subspace of
$\pi_{*}^N(\Ef^0\oplus\Ef^0)$. Then $\Pi^0_{i,j} = \{\Pi^z_{i,j} \
| \ z\in B\}\in\Aa(B,\Psi^{-\o}(\Ef^0\oplus\Ef^0))$ and it follows
that the kernels of the $L_{i,j}^z$ form a smooth finite-rank
super vector bundle $\Ker(L_{i,j}) = \ran(\Pi_{i,j})$ endowed with
the connection $\na_{0}^{i,j} =
\Pi_{i,j}^0\cdot\na^{\pi_{*}^N(\Ef^0\oplus\Ef^0)}\cdot\Pi^0_{i,j}$.
Clearly,
$$ \Ker(L_{i,j})  = \Ker(\Pp_i\Pp_j)\oplus\Ker((\Pp_i\Pp_j)^*)$$
relative to the grading of
$\pi_{*}^N(\Ef^0)\oplus\pi_{*}^N(\Ef^0)$, while
\begin{equation}\label{e:chker}
\ch(\na_{0}^{i,j}) = \ch(\Ind(\Pp_j,\Pp_i)) \ ,
\end{equation}
since $\Ind(\Pp_j,\Pp_i) = \Ker(\Pp_i\Pp_j) - \Ker(\Pp_j\Pp_i)$,
in this case, and $\ch(\na_{0}^{i,j})=\str((\na_{0}^{i,j})^2)$ is
defined since the bundle is finite-rank.
\begin{lem}\label{lem:toinfinity}
For each $C^l$ norm on $B$ one has
\begin{equation}\label{e:toinfinity}
    \lim_{t\too \o} \ch(\Bf^{1,2}_t, \Bf^{2,3}_t) = \ch(\na_{0}^{1,3})
\end{equation}
and a constant $C(l)$ such that
\begin{equation}\label{e:estimateinfinity}
   \|\tau(\Bf_{t}^{1,2},\Bf_{t}^{2,3})\|_l \leq C(l) t^{-3/2} \ .
\end{equation}
\end{lem}
\begin{proof}
Consider the spaces of $\pdo$ valued forms
\begin{equation*}
\Mm_{i,j} = \Aa(B,\Psi^*(\Ww_i\oplus\Ww_j)) \hskip 10mm {\rm
and}\hskip 10mm \Nn_{i,j} = \Aa(B,\Psi^{-\o}(\Ww_i\oplus\Ww_j))
\end{equation*}
where sections of $\Psi^*(\Ww_i\oplus\Ww_j)$ are families
$\Pp_{i,j}\cdot\Tf\cdot\Pp_{i,j}$, where $\Tf$ is a smooth family
of $\pdo$s on $\pi_{*}^N(\Ef^0\oplus\Ef^0)$, that is, a section of
$\Psi^*(\Ef^0\oplus\Ef^0)$, and similarly for
$\Psi^{-\o}(\Ww_i\oplus\Ww_j)$ but with $\Tf$ a smooth family of
smoothing operators.

Let $\Pi_{i,j}^0, \Pi_{i,j}^1\in\Aa(B,\Psi^{0}(\Ww_i\oplus\Ww_j))
$ be the smooth families of orthogonal projections onto
$\Ker(L_{i,j})$, $\Ker(L_{i,j})\pp$ as subspaces of
$\Ef^0\oplus\Ef^0$, or of $\Ww_i\oplus\Ww_j$. Since $L_{i,j}$ is
the restriction of the smooth family of elliptic operators $\Pp_i
\cdot(\Pp_i\cdot\Pp_j + \Pp_i\pp\cdot\Pp_j\pp)\ii \cdot\Pp_j,$
then by elliptic regularity on closed manifolds we have
$\Pi_{i,j}^0\in\Nn_{i,j}$. Relative to the decomposition
$\Ker(L_{i,j})\oplus\Ker(L_{i,j})\pp$ of $\Ww_i\oplus\Ww_j$ we
thus have, with $\Fff^{i,j} = \Fff^{i,j}_{t=1}$,
\begin{equation}\label{e:Fijker}
\Fff^{i,j} = \begin{bmatrix}
  \Pi_{i,j}^0\Fff^{i,j}\Pi_{i,j}^0 & \Pi_{i,j}^0\Fff^{i,j}\Pi_{i,j}^1 \\
  \Pi_{i,j}^1\Fff^{i,j}\Pi_{i,j}^0& \Pi_{i,j}^1\Fff^{i,j}\Pi_{i,j}^1 \\
\end{bmatrix} \in
\begin{bmatrix}
  \Nn^2_{i,j} & \Nn^1_{i,j} \\
  \Nn^1_{i,j} & \Mm^0_{i,j} \\
\end{bmatrix} \ ,
\end{equation}
where $\Mm_{i,j}, \Nn_{i,j}$ are filtered by
\begin{equation*}
\Mm_{i,j}^k = \sum_{r\geq k}\Aa^r(B,\Psi^*(\Ww_i\oplus\Ww_j)),
\hskip 10mm \Nn_{i,j}^k = \sum_{r\geq
k}\Aa^r(B,\Psi^{-\o}(\Ww_i\oplus\Ww_j)) \ .
\end{equation*}

Write $R_{i,j}^0 := (\na^{i,j})^2$. It is not quite true that the
curvature form $R_{i,j}^0$ coincides with the 2-form part of
$\Pi_{i,j}^0\Fff^{i,j}\Pi_{i,j}^0$, due to the Second Fundamental
Form defined by the off-diagonal terms in \eqref{e:Fijker}, rather
$$R_{i,j}^0 = (\Pi_{i,j}^0\Fff^{i,j}\Pi_{i,j}^0)_{[2]} -
(\Pi_{i,j}^0\Fff^{i,j}\Pi_{i,j}^1)_{[1]} G_{i,j}
(\Pi_{i,j}^1\Fff^{i,j}\Pi_{i,j}^0)_{[1]} \ , $$ where $G_{i,j}$ is
the inverse of $L_{i,j}$ restricted to $\Ker(L_{i,j})\pp$, and
zero elsewhere. However, by the Diagonalization Lemma of
\cite{BeVe90} there is an invertible $$g_{i,j}\in\Mm =
\Aa(B,\Psi^*(\Ef^0\oplus\Ef^0))$$ with, since
$\Pi_{i,j}^0\in\Nn_{i,j}$, $g_{i,j} - \If\in \Nn_1$, where $\Nn_k$
is the standard filtration on $\Nn =
\Aa(B,\Psi^{-\o}(\Ef^0\oplus\Ef^0))$, such that
\begin{equation*}
\Fff^{i,j} = g_{i,j}\ii\begin{bmatrix}
  U_{i,j} & 0 \\
  0 & V_{i,j} \\
\end{bmatrix}g_{i,j} =
g_{i,j}\ii\cdot \Pi_{i,j}^0\cdot U_{i,j}\cdot\Pi_{i,j}^0 \cdot
g_{i,j} + g_{i,j}\ii\cdot \Pi_{i,j}^1\cdot V_{i,j}\cdot\Pi_{i,j}^1
\cdot g_{i,j}
\end{equation*}
where
\begin{equation*}
U_{i,j} =  R_{i,j}^0 \ \ \mod\Nn_{i,j}^3 \ , \hskip 15mm V_{i,j} =
\Pi_{i,j}^1\cdot L_{i,j}\cdot\Pi_{i,j}^1 \ \ \mod\Nn_1 \ .
\end{equation*}
Moreover, since $\Fff^{i,j}_t = t\d_t(\Fff^{i,j})$, we have
$$\Fff_t^{i,j} = \d_t(g_{i,j})\ii\begin{bmatrix}
  t\d_t(U_{i,j}) & 0 \\
  0 & t\d_t(V_{i,j}) \\
\end{bmatrix}\d_t(g_{i,j})$$
\begin{equation}\label{e:deltatFt}
= \d_t(g_{i,j})\ii\cdot \Pi_{i,j}^0\cdot
t\d_t(U_{i,j})\cdot\Pi_{i,j}^0 \cdot \d_t(g_{i,j}) +
\d_t(g_{i,j})\ii\cdot \Pi_{i,j}^1\cdot
t\d_t(V_{i,j})\cdot\Pi_{i,j}^1 \cdot \d_t(g_{i,j})
\end{equation}
with
\begin{equation}\label{e:dtg}
\d_t(g_{i,j})^{\pm 1} = \If_{i,j} + O(t^{-1/2})
\end{equation}
\begin{equation}\label{e:dtu}
t\d_t(U_{i,j}) = R_{i,j}^0 + O(t^{-1/2})
\end{equation}
\begin{equation}\label{e:dtv}
t\d_t(V_{i,j}) = tL_{i,j}^2 + t^{1/2}\Mm_{i,j}^1 \ ,
\end{equation}
$\If_{i,j}$ the identity operator on $\Ww_i\oplus\Ww_j$. Here, a
1-parameter family of $A_t\in\Nn$, $t > 0$, is $O(f(t))$ for some
positive function $f(t)$ if for any $\e>0, l\in\Nf$ and
$\phi\in\Ci_c (B)$, one has $\|(\pi^{N})^*(\phi)(x)A_t (x,y)\|_l
\leq C_{l,\e,\phi} f(t)$ for $t>\e$, while $B_t\in
t^{\d}\Mm_{i,j}^k$ if $B_t = \sum_{r\geq 0}
t^{\d_r}\omega_{[k+r]}$ with $\d_0 = \d > \d_r \searrow -\o$ and
$\omega_{[k+r]}\in\Mm_{i,j}^k$.


It follows from \eqref{e:deltatFt} and \eqref{e:inclusion} that
\begin{equation}\label{e:relttoinfty}
    (\Fff_t^{1,2} - \la\If)_{|\Ef^0}\ii -  (\Fff_t^{2,3} - \la\If)_{|\Ef^0}\ii
\end{equation}
\begin{equation*}
= \d_t(g_{1,2})\ii\cdot \Pi_{1,2}^0\cdot (t\d_t(U_{1,2}-
\la\If^0_{1,2})\ii\cdot\Pi_{1,2}^0 \cdot \d_t(g_{1,2})  -
\d_t(g_{2,3})\ii\cdot \Pi_{2,3}^0\cdot (t\d_t(U_{2,3}-
\la\If^0_{2,3})\ii\cdot\Pi_{2,3}^0 \cdot \d_t(g_{2,3})
\end{equation*}
\begin{equation*}
+ \d_t(g_{1,2})\ii\cdot \Pi_{1,2}^1\cdot (t\d_t(V_{1,2}-
\la\If^1_{1,2})\ii\cdot\Pi_{1,2}^1 \cdot \d_t(g_{1,2}) -
\d_t(g_{2,3})\ii\cdot \Pi_{2,3}^1\cdot (t\d_t(V_{2,3}-
\la\If^1_{2,3})\ii\cdot\Pi_{2,3}^1 \cdot \d_t(g_{2,3}) \ .
\end{equation*}
where $\If_{i,j}^0, \If_{i,j}^1$ are the vertical identity
operators on $\Ker(L_{i,j})$ and $\Ker(L_{i,j})\pp$, respectively,
and zero elsewhere, and $\d_t(g_{i,j}) =\d_t(g_{i,j})_{|\Ef^0}$ is
included as an element of $\Mm$. The first two terms, involving
the resolvents of $U_{i,j}$, are equal to
\begin{equation}\label{e:1st2terms}
\d_t(g_{1,2})\ii\cdot \Pi_{1,2}^0\cdot (R^0_{1,2} -
\la\If_{1,2}^0)\ii\cdot\Pi_{1,2}^0 \cdot \d_t(g_{1,2})  -
\d_t(g_{2,3})\ii\cdot \Pi_{2,3}^0\cdot (R^0_{2,3}-
\la\If_{2,3}^0)\ii\cdot\Pi_{2,3}^0 \cdot \d_t(g_{2,3})
\end{equation}
\begin{equation*}
+ \d_t(g_{1,2})\ii\cdot \Pi_{1,2}^0\cdot
F^{[3]}_{1,2}(t^{1/2},\la)\cdot\Pi_{1,2}^0 \cdot \d_t(g_{1,2})  -
\d_t(g_{2,3})\ii\cdot \Pi_{2,3}^0\cdot
F^{[3]}_{2,3}(t^{1/2},\la)\cdot\Pi_{2,3}^0 \cdot \d_t(g_{2,3}) \ ,
\end{equation*}
where the notation $F^{[m]}_{i,j}(t^{\e},\la)$ means a function of
the form
$$F^{[m]}_{i,j}(t^{\e},\la) = \sum_{k=1}^{\dim B} t^{-\e k}
(A_{i,j} - \la\If_{i,j}^0)\ii w_{i,j}[t]$$ for some differential
form valued bundle endomorphism $A_{i,j}$, and $w_{i,j}[t] \in t^0
\Nn_{i,j}^m$ polynomial in $t$. From \eqref{e:dtg}
\begin{equation*}
\d_t(g_{i,j})^{\pm 1}\cdot \Pi_{i,j}^0 = (\Pi_{i,j}^0 +
\Pi_{i,j}^1)\d_t(g_{i,j})^{\pm 1}\cdot \Pi_{i,j}^0  = \Pi_{i,j}^0
+ \sum_{k=0}^1\Pi_{i,j}^k \cdot O(t^{-1/2})\cdot \Pi_{i,j}^0 \ ,
\end{equation*}
\begin{equation*}
\Pi_{i,j}^0\cdot \d_t(g_{i,j})^{\pm 1} = \Pi_{i,j}^0 +
\sum_{k=0}^1\Pi_{i,j}^0 \cdot O(t^{-1/2})\cdot \Pi_{i,j}^k \ .
\end{equation*}
Hence \eqref{e:1st2terms} has the form
\begin{equation}\label{e:1st2terms2}
\Pi_{1,2}^0\cdot (R^0_{1,2} - \la\If_{1,2}^0)\ii\cdot\Pi_{1,2}^0 -
\Pi_{2,3}^0\cdot (R^0_{2,3} - \la\If_{2,3}^0)\ii\cdot\Pi_{2,3}^0
\end{equation}
\begin{equation*}
+ \sum_{i=1}^2\sum_{k=0}^1 \Pi_{i,i+1}^k\cdot
F^{[0]}_{i,i+1}(t^{-1/2},\la)\cdot \Pi_{i,i+1}^{1-k}
+\Pi_{1,2}^0\cdot B^{[0]}_{1,2}(t^{-1},\la)\cdot\Pi_{1,2}^0
-\Pi_{2,3}^0\cdot B^{[0]}_{2,3}(t^{-1},\la)\cdot\Pi_{2,3}^0
\end{equation*}
From \eqref{e:dtv}
\begin{equation}\label{e:2nd2terms}
\Pi_{i,j}^1\cdot (t\d_t(V_{i,j}-\la\If_{i,j})\ii\cdot\Pi_{i,j}^1
 = \Pi_{i,j}^1\cdot (tL^2_{i,j}-\la\If_{i,j})\ii\cdot\Pi_{i,j}^1
+ \Pi_{i,j}^1\cdot C^{[0]}_{i,j}(t^{-1/2},\la)\cdot\Pi_{i,j}^1 \ .
\end{equation}
Hence we obtain from \eqref{e:relttoinfty}, \eqref{e:1st2terms},
\eqref{e:1st2terms2}, \eqref{e:2nd2terms}
\begin{equation}\label{e:relatinfty}
    e^{-\Fff^{1,2}_t} - e^{-\Fff^{2,3}_t}
\end{equation}
\begin{equation*}
=\Pi_{1,2}^0\cdot \frac{i}{2\pi}\int_{\Cc} e^{-\la} \,(R_{1,2}^0 -
\la \If_{1,2})) \, d\la \cdot \Pi_{1,2}^0 - \Pi_{2,3}^0\cdot
\frac{i}{2\pi}\int_{\Cc} e^{-\la} \,(R_{2,3}^0 - \la \If_{2,3}))
\, d\la \cdot \Pi_{2,3}^0
\end{equation*}
\begin{equation*}
+ \sum_{i=1}^2\sum_{k=0}^1 \Pi_{i,i+1}^k\cdot
\frac{i}{2\pi}\int_{\Cc} e^{-\la} \, F^{[0]}_{i,i+1}(t^{-1/2},\la)
\, d\la \cdot \Pi_{i,i+1}^{1-k}
\end{equation*}
\begin{equation*}
+ \sum_{i=1}^2 (-1)^{i+1}\d_t(g_{i,i+1})\ii\left(
\Pi_{i,i+1}^1\cdot\frac{i}{2\pi}\int_{\Cc} e^{-\la} \,
(tL^2_{i,i+1}-\la\If_{i,i+1})\ii
 \, d\la \cdot \Pi_{i,i+1}^{1}
 \right.
\end{equation*}
\begin{equation*}
\left. +
 \Pi_{i,i+1}^1\cdot\frac{i}{2\pi}\int_{\Cc} e^{-\la} \,
A_{i,i+1}(t^{-1/2},\la)
 \, d\la \cdot \Pi_{i,i+1}^{1}\right)\d_t(g_{i,i+1})
\end{equation*}
\begin{equation}\label{e:estimatesatinfty}
=\Pi_{1,2}^0\cdot e^{-R_{1,2}^0} \Pi_{1,2}^0 - \Pi_{2,3}^0\cdot
e^{-R_{2,3}^0} \Pi_{2,3}^0
\end{equation}
\begin{equation*}
+ \sum_{i=1}^2\sum_{k=0}^1 \Pi_{i,i+1}^k\cdot O(t^{-1/2}) \cdot
\Pi_{i,i+1}^{1-k}
\end{equation*}
\begin{equation*}
+ \d_t(g_{1,2})\ii \cdot\Pi_{1,2}^1\cdot e^{-tL^2_{1,2}} \cdot
\Pi_{1,2}^{1}\cdot \d_t(g_{1,2}) - \d_t(g_{2,3})\ii
\cdot\Pi_{2,3}^1\cdot e^{-tL^2_{2,3}} \cdot \Pi_{2,3}^{1}\cdot
\d_t(g_{2,3})
\end{equation*}
\begin{equation*}
+ \d_t(g_{1,2})\ii \cdot\Pi_{1,2}^1\cdot O(t^{-1/2})
\Pi_{1,2}^{1}\cdot \d_t(g_{1,2}) - \d_t(g_{2,3})\ii
\cdot\Pi_{2,3}^1\cdot O(t^{-1/2}) \cdot \Pi_{2,3}^{1}\cdot
\d_t(g_{2,3})
 \end{equation*}
The operators in the first line of \eqref{e:estimatesatinfty} are
finite-rank with supertraces $\str(e^{-R^0_{i,j}})$. The terms in
the second line are also finite-rank but with vanishing traces by
symmetry and $\Pi_{i,j}^{0}\Pi_{i,j}^{1} = 0$. The terms in the
third and fourth line are not trace class, but the differences are
trace class and of order $O(e^{-ct})$, some $c>0$, and
$O(t^{-1/2})$ as $t\too \o$.

Totally as an element of $\Aa(B)$ we therefore obtain
\begin{equation*}
\ch(\Bf^{1,2}_t, \Bf^{2,3}_t)  = \str(e^{-R^0_{1,2}}) -
\str(e^{-R^0_{2,3}}) + O(t^{-1/2}) \ ,
\end{equation*}
proving equation \eqref{e:toinfinity}.

For \eqref{e:estimateinfinity} we can revert for simplicity to the
case of the canonical superconnections (the proof of
\eqref{e:toinfinity} holds for any induced superconnection). Then
clearly
$$\dd_t \Bf^{1,2}_t =\left\{%
\begin{array}{ll}
    0 \ , & {\rm on} \ \ \Ker(L_{i,j}) = \ran(\Pi_{i,j}^0) \\
    t^{-1/2}L_{i,j} \ , & {\rm on} \ \ \Ker(L_{i,j})\pp \\
\end{array}%
\right.$$ It is then immediate from \eqref{e:estimatesatinfty}
that the estimate \eqref{e:estimateinfinity} holds.
\end{proof}

Notice that it follows from \propref{e:tform0} and
\eqref{e:estimateinfinity} that the integral $\int_{0}^{\o}
\tau(\Bf_{t}^{1,2},\Bf_{t}^{2,3}) \ dt$ exists, and we have then
the transgression formula:
\begin{cor}
In $\Aa(B)$
\begin{equation*}
\ch(\Ind(\Pp_1,\Pp_3)) = \eta(\Pp_3 ,\Pp_1 ) - d_B\int_{0}^{\o}
\tau(\Bf_{t}^{1,2},\Bf_{t}^{2,3}) \ dt \ .
\end{equation*}
\end{cor}

\vskip 2mm

The extension  to the general case is obtained in the manner of
\cite{BeGeVe92}, applying the above argument to a perturbed family
$L_{i,j} \oplus Z_{i,j} : \Ww_i\oplus\Ww_j\oplus\Cf^N\too
\Ww_i\oplus\Ww_j$ with constant kernel dimension and whose index
bundle defines the same class in $K(B)$ as $\Ind(\Pp_i,\Pp_j)$.
\end{proof}

\vskip 5mm

\section{Chern Character Forms from the Interior}

In this Section we construct a superconnection directly on the
bundle $\pi_{*}(\Ef|\Pp)$ whose Chern character form represents
$\ch(\Ind(\Ds_{\Pp}))$. We then prove Theorem (II) by
transgressing this form using a generalized smoothing operator
calculus adapted to manifolds with boundary.

\subsection{Construction of the Heat Operator}
Let $\Af$ be the canonical superconnection
\eqref{e:canonicalsuperconnection} on $\pi_{*}(\Ef)$. Let $\Pp$ be
a Grassmann section for $\Ds$, and let $\Ps =
\Ps(\Pp)\in\Aa^0(B,\End(\pi_*(\Ef)))$ be the corresponding family
of projections \eqref{e:P} onto $\pi_{*}(\Ef|\Pp)$ The induced
superconnection $\Ps\cdot\Af_t \cdot\Ps $ on the subbundle
$\pi_{*}(\Ef|\Pp)$ of $\pi_{*}(\Ef)$ takes the form
\begin{equation}\label{e:inducedcanonicalsuper}
\Ps\cdot\Af_t \cdot\Ps  =  \na^{\pi_*(\Ef | \Pp)} +
t^{1/2}\Ps\cdot \Ds\cdot\Ps
\end{equation}
where $\na^{\pi_*(\Ef | \Pp)} = \Ps\cdot \nE \cdot\Ps$ is the
induced connection on $\pi_{*}(\Ef|\Pp)$, in the usual sense. The
curvature of $\Ps\cdot\Af_t \cdot\Ps $ is the 2-form in
$\Aa^2(B,\End(\pi_{*}(\Ef|\Pp)))$
\begin{equation}\label{e:inducedcurv}
(\Ps\cdot\Af_t \cdot\Ps )^2 = \Rs^{\pi_{*}(\Ef|\Pp)} +
t^{1/2}\na^{\pi_*(\Ef | \Pp)} (\Ps\cdot \Ds\cdot\Ps) + t (\Ps\cdot
\Ds\cdot\Ps)^2 \ .
\end{equation}
Here $\Rs^{\pi_{*}(\Ef|\Pp)}  = (\na^{\pi_*(\Ef | \Pp)})^2$ is the
curvature 2-form of $\na^{\pi_*(\Ef | \Pp)}$, given by
\begin{equation}\label{e:inducedcurv2}
\Rs^{\pi_{*}(\Ef|\Pp)} =  -\Ps\cdot \nE_{\sigma}\cdot\Ps +
\Ps\cdot R^{M/B}\cdot\Ps - II^{*}_{\Pp}\wedge II_{\Pp} \ ,
\end{equation}
where $II_{\Pp} = \nE_{|\pi_{*}(\Ef|\Pp)} - \na^{\pi_*(\Ef |
\Pp)}$ is the second fundamental form.

\vskip 2mm

The formal construction of a heat operator on $\pi_{*}(\Ef|\Pp)$
from this superconnection is problematic as $(\Ps\cdot
\Ds\cdot\Ps)^2$ is neither a differential operator nor formally
self-adjoint. However, its restriction to the subbundle
$\pi_{*}(\Ef|\Pp^3)$ coincides with the Dirac Laplacian
\begin{equation}\label{e:equaltolaplacian}
(\Ps\cdot \Ds\cdot\Ps)_{|\pi_{*}(\Ef|\Pp^3)}^2 =
(\Ds^2_{\Pp})_{|_{\pi_{*}(\Ef|\Pp^3)}}  \ ,
\end{equation}
which is a self-adjoint differential operator on
$\pi_{*}(\Ef|\Pp^2)$. Here, for each integer $k\geq 1$ we have an
$L^2$-dense subbundle of $\pi_{*}(\Ef)$
\begin{equation}\label{e:EPK}
\pi_{*}(\Ef|\Pp^k) = \{\psi\in\pi_{*}(\Ef)\ | \ \Pp\g \Ds^i\psi =
0 \ , \ 0\leq i \leq k-1\} \ .
\end{equation}
Thus for $k=1$ this is $\pi_{*}(\Ef | \Pp)$, the domain of
$\Ds_{\Pp}$, while $k=2$ is the domain of the family of the
Laplacian boundary problems
\begin{equation*}
\Ds^2_{\Pp } : \Aa(B,\pi_{*}(\Ef|\Pp^2)) \too \Aa(B,\pi_{*}(\Ef))
\ .
\end{equation*}
The restriction of the curvature of $\Ps\cdot\Af_t \cdot\Ps$ to
$\pi_{*}(\Ef|\Pp^3)$ has range in $\pi_{*}(\Ef|\Pp)$ equal there
to the curvature of the superconnection with domain
$\pi_{*}(\Ef|\Pp)$
\begin{equation}\label{e:superP}
\Af_{t,\Pp} := \na^{\pi_*(\Ef | \Pp)} + t^{1/2}\Ds_{\Pp}
\end{equation}
 adapted to the family of global Dirac boundary
problems $\Ds_{\Pp}$. One has $$\Af_{t,\Pp} = (\Ps\cdot\Af_t
\cdot\Ps)_{|\pi_{*}(\Ef|\Pp^2)} \ .$$ Thus $\Af_{t,\Pp}$ extends
$(\Ps\cdot\Af_t \cdot\Ps)_{|\pi_{*}(\Ef|\Pp^2)}$ to
$\pi_{*}(\Ef|\Pp)$ with range $\pi_{*}(\Ef)$ and self-adjoint zero
order part. The curvature of $\Af_{t,\Pp}$,
 defined by $$\Ff_{t,\Pp} = (\Af^2_{t,\Pp})_{|\pi_{*}(\Ef|\Pp^2)}
\in \Aa(B,\Hom(\pi_{*}(\Ef|\Pp^2),\pi_{*}(\Ef)) ) \ $$ is the
vertical differential operator, supercommuting with action of
$\Aa(B)$,
\begin{equation}\label{e:superPcurv}
\Ff_{t,\Pp} = \Rt + t^{1/2}\na^{\pi_*(\Ef | \Pp)}\Ds_{\Pp} +
t\Ds_{\Pp}^2 : \Aa(B,\pi_{*}(\Ef|\Pp^2)) \too \Aa(B,\pi_{*}(\Ef))
\ ,
\end{equation}
(Recall here that $\Ds_{\Pp} = \Ds$ as an operator, the subscript
just refers to the domain restriction.)

\begin{prop}
For any $\la \in \Cf\backslash\Rf_+$, $\Ff_{t,\Pp}$ has a
resolvent
\begin{equation}\label{e:superresolvent}
(\Ff_{t,\Pp} - \la\If)\ii \in \Aa(B,\End(\pi_{*}(\Ef))) \ ,
\end{equation}
with range in $\Aa(B,\pi_{*}(\Ef|\Pp^2))$, that is,
\begin{equation}\label{e:rangeinP}
\Pp\g (\Ff_{t,\Pp} - \la\If)\ii = 0 \ , \hskip 10mm
\Pp\g\Ds(\Ff_{t,\Pp} - \la\If)\ii = 0 \ ,
\end{equation}
where $\If$ is the vertical identity operator on
$\Aa(B,\pi_{*}(\Ef))$. One has
\begin{equation}\label{e:nilpotent}
(\Ff_{t,\Pp} - \la\If)\ii  = \sum_{k=0}^{\dim B}(-1)^k
(t\Ds_{\Pp}^2-\la\If)\ii \left((\Rs + t^{1/2}\na^{\pi_*(\Ef |
\Pp)}\Ds_{\Pp})(t\Ds_{\Pp}^2 - \la\If)\ii\right)^k \  .
\end{equation}
\end{prop}
\begin{proof}
The first statement is just the statement that
$(\Ff_{t,\Pp})_{[0]} = \Ds_{\Pp}^2$ is a family of self-adjoint
generalized Laplacians and hence has a resolvent
$(t\Ds_{\Pp}^2-\la\If)\ii$ for $\la \in \Cf\backslash\Rf_+$. Then
\begin{eqnarray*}
(\Ff_{t,\Pp} - \la\If)\ii & = & (\Rs + t^{1/2}\na^{\pi_*(\Ef |
\Pp)}\Ds_{\Pp} +
t\Ds_{\Pp}^2 - \la\If)\ii \\
& = & (t\Ds_{\Pp}^2-\la\If)\ii \left((\Rs + t^{1/2}\na^{\pi_*(\Ef
| \Pp)}\Ds_{\Pp})(t\Ds_{\Pp}^2 - \la\If)\ii + \If\right)\ii \ ,
\end{eqnarray*}
and since $(t\Ds_{\Pp}^2-\la\If)\ii$ has range in
$\Aa(B,\pi_{*}(\Ef|\Pp^2))$, \eqref{e:rangeinP} follows. Since
$$(\Rs^{\pi_*(\Ef | \Pp)}+ t^{1/2}\na^{\pi_*(\Ef | \Pp)}\Ds_{\Pp})(t\Ds_{\Pp}^2 - \la\If)\ii \in
\Aa^1(B,\End(\pi_{*}(\Ef)))$$ it is nilpotent of degree $\dim B$,
consisting of terms which raise form degree by 1 or 2, and
\eqref{e:nilpotent} is the corresponding finite Neumann expansion.
\end{proof}

\noi Now $(\Rs^{\pi_*(\Ef | \Pp)}+ t^{1/2}\na^{\pi_*(\Ef |
\Pp)}\Ds_{\Pp})(t\Ds_{\Pp}^2 - \la\If)\ii$ is bounded in the
operator norm , since the first factor is a first-order vertical
differential operator, while $(t\Ds_{\Pp}^2 - \la\If)\ii$ has
Sobolev order -2, then \eqref{e:nilpotent} implies the estimate
$$\|(\Ff_{t,\Pp} - \la\If)\ii\| = O(|\la|^{-1/2})$$ as $|\la|\to \o$
(for details, see the proof of \propref{p:limitexists} where a
similar fact is used), as an operator on $\Aa(M,\pi^*(\wedge
T*B)\otimes \Ef)$. Hence we can define the heat operator:

\begin{defn}
For $t>0$ the heat operator is the convergent contour integral
\begin{equation}\label{e:heatP} e^{-\Ff_{t,\Pp}} :=
\frac{i}{2\pi}\int_{\Cc}e^{-\la} (\Ff_{t,\Pp} - \la\If)\ii \ d\la
\end{equation}
where $\Cc$ is a contour surrounding the positive real axis, as in
the Introduction.
\end{defn}
$e^{-\Ff_{t,\Pp}}$ is convergent as a vertical operator in
$C^0(B,\pi_*(\End( \Ef)))$. More precisely, on a vertical family
of bounded operators in $C^0(B,\pi_*(\End( \Ef)))$ the $sup$-norm
is defined via compact subsets $C$ of $B$ by $\|T\|_C =
\sup\{\|T_z\| \ | \ z\in B\}$ . Using identities such as
$$\na_{\xi}^{\pi_*(\Ef|\Pp)}(t\Ds_{\Pp}^2 - \la\If)\ii =
-t(t\Ds_{\Pp}^2 - \la\If)\ii
(\na_{\xi}^{\pi_*(\Ef|\Pp)}\Ds_{\Pp}^2)(t\Ds_{\Pp}^2 -
\la\If)\ii$$ we obtain from \eqref{e:nilpotent} the estimate as
$\la\to\o$
$$\|\na_{\xi_1}^{\pi_*(\Ef|\Pp)}\ldots\na_{\xi_k}^{\pi_*(\Ef|\Pp)}(
\Rs + t^{1/2}\na^{\pi_*(\Ef | \Pp)}\Ds_{\Pp} + t\Ds_{\Pp}^2 -
\la\If)\ii \|_C \leq a(k,C) t^k |\la|^{-k} \ ,$$ for any
$\xi_1,\dots\xi_k\in {\rm Vect}(B)$ and a constant $a(k,C)$.  This
proves convergence in each $C^l$ norm:
\begin{prop}
$$e^{-\Ff_{t,\Pp}}\in \Aa(B,\End(\pi_{*}(\Ef))) \ .$$
\end{prop}

To see $e^{-\Ff_{t,\Pp}}$ has a well defined trace in $\Aa(B)$ we
show it is a vertical `smoothing' operator, in the following
generalized sense.

\subsection{Generalized Smoothing Operators}
A smooth family of (generalized) smoothing operators is an element
of the space $\Aa(B,\Psi^{-\o}(\Ef,\Ef^{'}))$ of vertical
smoothing operators acting on $\Aa(M,\Ef)$. Vertical means
supercommuting with the action of $\Aa(B)$. A section of
$\G(M,\Psi^{-\o}(\Ef,\Ef^{'}))$ is an operator of the form
\begin{equation}\label{e:gensmoperator}
    \Ks = \Ks_{X} + \Gs^{-\o} \ ,
\end{equation}
where $\Ks_{X}$ is a smooth family of interior smoothing
operators, a family of $\pdo$ in the usual sense, and $\Gs^{-\o}$
is a smooth family of singular Green's operators (sgo) of order
$-\o$. Thus, $\Ks_{X}$ is defined by a kernel
\begin{equation}\label{e:intsmoperator}
    <x_1 | \Ks_{X} | x_2 > \in \G(M\times_{\pi}M,\pi^*(\wedge T^*
    B) \otimes \Ef \boxtimes_{\pi}\Ef),
\end{equation}
$\Ci$ up to the boundary $\dd M$--equivalently, it is the
restriction of a vertical smoothing operator from the closed
double of $M$.

A vertical sgo is a vertical operator of the form
\begin{equation}\label{e:sgo}
    \Gs^{-\o} = \Kk\Tt \Gg_X \in \Aa(B,\End(\pi_{*}(\Ef)))
\end{equation}
where
\begin{equation}\label{e:poisson}
    \Kk \in \Aa(B,\Hom(\pi^{N}_{*}(\Ef^{'}), \pi_{*}(\Ef)))
\end{equation}
\begin{equation}\label{e:trace}
    \Tt  = \sum_{0\leq j \leq r} \Ss_j \g_j + \widetilde{\Tt}
    \in \Aa(B,\Hom(\pi_{*}(\Ef)), \pi^{N}_{*}(\Ef^{'}))
\end{equation}
\begin{equation}\label{e:interiorsgo}
     \Gg_X \in \Aa(B,\End(\pi_{*}(\Ef)))
\end{equation}
are, respectively, smooth vertical families of {\it Poisson,
trace}, and {\it interior singular Green's operators} of order
$m_1, m_2, m_3$ and of total order $-\o$. This means the
following. For any local trivialization of $\pi_{*}(\Ef)$ we have
an identification
\begin{equation}\label{e:Etriv3}
\Aa(U,\End(\pi_* (\Ef))) = \Aa(U)\otimes \End(\G(X_{z_0},E_{z_0}))
\ .
\end{equation}
Then $\Gs^{-\o}$ is a smooth family of sgo if with respect to each
such trivialization
\begin{equation}\label{e:localsgo}
     \Gs^{-\o}_{|U} = \sum_{I}G_{I}^z dz_I
\end{equation}
where $G_{I}^z\in \End(\G(X_{z_0},E_{z_0}))$ is a sgo of order
$-\o$ in the classical single operator sense \cite{Gr96} depending
smoothly on $z\in U$. Similarly, $\Kk$ is a smooth family of
Poisson operators if locally $\Kk_{|U} = \sum_{I}K_{I}^z dz_I$ ,
where $K_{I}^z : \G(\dd X_{z_0},E^{0}_{z_0}) \too
\G(X_{z_0},E_{z_0})$ is a classical Poisson operator of order $m_1
\leq 0$, defined by a symbol kernel satisfying the estimates in
\cite{Gr96,Gr00} and depending smoothly on $z$. The trace operator
$\Tt$ takes the local form
\begin{equation}\label{e:glocatrace}
\Tt_{|U} = \sum_{I}\sum_{0}^{r}(S^z_{j,I} \g^{z_0}_j +
\widetilde{T}_{I}^z) dz_I \ ,
\end{equation}
where $\g^{z_0}_j s= \frac{\dd^j}{\dd r^j}s_{|\dd X_{z_0}}$ is the
restriction map assigning the $j^{{\rm th}}$ normal derivative at
$\dd X_{z_0}$, $S^z_{j,I}$ is the local form of a smooth family of
boundary pseudodifferential operators
$\Ss_j\in\Aa(B,\Psi(\Ef^{'}))$, and $\widetilde{T}_{I}^z :
\G(X_{z_0},E_{z_0})\too \G(\dd X_{z_0},E^{0}_{z_0})$ a trace
operator of order $m_2\leq 0$ defined by a classical symbol kernel
 smooth in $z$. Finally, $\Gg_X = \sum_{I}G_{X,I}^z
dz_I$ with  $G_{X,I}^z$ a sgo operator of order $m_3\leq 0$
defined by a classical symbol kernel depending smoothly on $z$
\cite{Gr96,Gr00}.

If the orders of $\Kk, \Ss_j, \widetilde{\Tt}, \Gg_X$ are
sufficiently negative they are given by integral operators with
continuous  kernels; thus if $m_1 \leq -\dim(X) - k$ then $\Kk$ is
given by a kernel
\begin{equation*}
    <x| \Kk | y > \in \G^k (M\times_{\pi}\dd M,\pi^*(\wedge T^*
    B) \otimes \Ef \boxtimes_{\pi}\Ef^{'}) \ ,
\end{equation*}
where $M\times_{\pi}\dd M = \{(x,y)\in M\times\dd M \ | \ \pi(x) =
\pi^N (y)\}$, and  so on. Since $\Gs^{-\o}$ is of order $-\o$ it
extends to a bounded operator $\Aa^s (B,\pi_*(\Ef))\to \Aa^{s
+d}(B,\pi_*(\Ef))$ on the Sobolev completions for any $d>0$. It
follows that $\Gs^{-\o}$ is a vertical trace class operator. Hence
there is a trace on $\Psi^{-\o}(\Ef,\Ef^{'})$ over each fibre
$X_z$ of $M$ defining a supertrace
\begin{equation}\label{e:supertrace}
\str : \Aa(B, \Psi^{-\o}(\Ef,\Ef^{'})) \too \Aa(B) \ .
\end{equation}
More precisely, by the boundary pseudodifferential operator
calculus of \cite{Gr96,Gr00} the cycled operator
\begin{equation}\label{e:cycled}
\Ks_{\dd X} = \widetilde{\Tt}\Gg_X\Kk\in
\Aa(B,\Psi^{-\o}(\Ef^{'}))
\end{equation}
 is a smooth family of smoothing operators over the family of closed
manifolds $N$   in the usual sense \cite{BeGeVe92} given by a
$\Ci$ kernel
 \begin{equation*}
<y_1| \Ks_{\dd X} | y_2 > \in \G^k (N\times_{\pi^N}
N,(\pi^N)^*(\wedge T^*B) \otimes \Ef^{'} \boxtimes_{\pi^N}\Ef^{'})
\ .
\end{equation*}
Using the invariance of the trace with respect to cyclic
permutation of products of bounded and trace-class operators, and
Lidskii's Theorem, we therefore obtain:
\begin{lem}\label{lem:str}
The supertrace  \eqref{e:supertrace} of an operator $\Ks = \Ks_{X}
+ \Gs^{-\o}\in \Aa(B, \Psi^{-\o}(\Ef,\Ef^{'}))$ is the
differential form on $B$
\begin{equation}\label{e:supertracekernel}
z\mtoo \int_{X_z} \str <x | \Ks^z_{X} | x >  + \int_{\dd X_z} \str
<y | \Ks^z_{\dd X} | y >  \ ,
\end{equation}
where $\Ks_{\dd X} $ is the cycled boundary operator
\eqref{e:cycled}.
\end{lem}

\vskip 2mm The boundary trace of the singular Green's operator
term in \eqref{e:supertracekernel} has the consequence that  the
trace of the commutator of a differential operator with a
smoothing operator is non-zero, in contrast to closed manifolds
\cite{BeGeVe92}:
\begin{lem}\label{lem:commutatortrace}
Let $\Df\in\Aa(B,\Dd(\Ef))$ be a vertical differential operator of
order $r$ and let $K\in\Aa(B,\Psi^{-\o}(\Ef,\Ef^{'}))$. Then\\[1mm]
{\rm [1]} $\Df K$ and $K\Df$ are elements of
$\Aa(B,\Psi^{-\o}(\Ef,\Ef^{'}))$.\\[1mm]
{\rm [2]} $[\Df, K] = A\Kk\rho,$ where $\Kk$ is a Poisson operator
of order $-\o$,  $\rho = (\g_0,\g_1,\ldots,\g_{r-1})$, and $A$ is
a boundary Green's form.
 $\Tr([\Df, K]) $ is the differential form on $B$
\begin{equation}\label{e:commutatortrace}
z\mtoo \int_{\dd X_z} \Tr( A_z <y | \Kk | y> )\   dy \ .
\end{equation}
\end{lem}

For certain preferred classes of vertical smoothing operators the
vanishing of the trace on such commutators can nevertheless be
restored by imposing APS-boundary conditions $\Pp$ on $\Ds$ . The
trangression formula for the Chern character, below, depends on
this fact. The proof of \lemref{lem:commutatortrace} of this lemma
is a straightforward application of Green's formula -- see
\cite{DaFr97} for the case of a Dirac operator.

\subsection{The Chern Character}

\begin{prop}\label{p:heat}
The integral \eqref{e:heatP} converges in the $\Ci$ topology to an
element of $\Aa(B,\Psi^{-\o}(\Ef,\Ef^{'}))$. In particular,
$e^{-\Ff_{t,\Pp}} $ is trace class for each $t>0$. Moreover, with
$\Ff_{\Pp} := \Ff_{1,\Pp}$, $\Hf(t) = e^{-t\Ff_{\Pp}} $ is the
unique solution to the heat equation along the fibre $X_z$:
\begin{equation}\label{e:heatequation}
    \left(\frac{d}{dt} +( \Af_{P}^z)^2 \right)\Hf(t)^z = 0
\end{equation}
subject to $\lim_{t\to 0}\Hf(t)^z\psi^z = \psi^z$ for $\psi^z\in
\pi^*(\wedge T^{*}_{z}B)\otimes\G(X_z,E_z)$. One has
\begin{equation}\label{e:heatdeltat}
e^{-\Ff_{t,\Pp}} = \d_t(e^{-t\Ff_{\Pp}}) \ .
\end{equation}
\end{prop}

First, we require the following Lemma. For $k=1,2$ let $\Ps_k =
\Ps(\Pp_k)$ be the families of projection operators associated to
Grassmann sections $\Pp_k \in\Aa(B,\End(\pi_{*}^{N}(\Ef^{'})))$
defining the subbundles $\pi_*(\Ef|\Pp_k)$ of $\pi_* (\Ef)$. Let
$\Rs_k = \Rs^{\pi_*(\Ef|\Pp_k)}$ be the vertical curvature forms.

\begin{lem}\label{lem:relativelysmooth}
 \begin{equation}\label{e:smooth0}
  \Ps_1 - \Ps_2 \in \Aa(B,\Psi^{-\o}(\Ef,\Ef^{'}))
\end{equation}
with supertrace
\begin{equation}\label{e:smooth1}
    \str(\Ps_1 - \Ps_2) = 2\,\ind(\Pp_2\Pp_1) \ .
\end{equation}
\begin{equation}\label{e:smooth2}
    (t\Ds_{\Pp_1}^2-\la\If)\ii - (t\Ds_{\Pp_2}^2-\la\If)\ii \in
    \Aa(B,\Psi^{-\o}(\Ef,\Ef^{'}))\ .
\end{equation}
For $k \in\Nf$
\begin{equation}\label{e:smooth3}
\Rs_1^{k} - \Rs_2^{k}\in \Aa(B,\Psi^{-\o}(\Ef,\Ef^{'}))\ .
\end{equation}
\begin{equation}\label{e:smooth4}
\na^{\pi_*(\Ef|\Pp_1)}\cdot \Ds - \na^{\pi_*(\Ef|\Pp_2)}\cdot
\Ds\in \Aa(B,\Psi^{-\o}(\Ef,\Ef^{'}))\ .
\end{equation}
\end{lem}
\begin{proof}
From \eqref{e:P}
\begin{equation}\label{e:smooth5}
    \Ps_1 - \Ps_2 = \Ks_{\chi}(\Pp_2 - \Pp_1)\g \ ,
\end{equation}
where $\Ks_{\chi}$ is a vertical Poisson operator, while $\Pp_2 -
\Pp_1 \in \Aa(B,\Psi^{-\o}(\Ef^{'}))$. Hence $\Ps_1 - \Ps_2$ is a
sgo operator of order $-\o$. Cycling the trace we have
\begin{eqnarray*}
   \str(\Ps_1 - \Ps_2) & = & \str\left(\Ks_{\chi}
   \begin{bmatrix}
     \Pp_2 - \Pp_1 & 0 \\
      & \Psi\ii(\Pp_2 - \Pp_1)\Psi\\
   \end{bmatrix}\g \right) \\
& = & \Tr(\Pp_2 - \Pp_1) - \Tr(\Psi\ii(\Pp_2 - \Pp_1)\Psi) =
2\Tr(\Pp_2 -
 \Pp_1)\ .
\end{eqnarray*}

To see \eqref{e:smooth2} we have from \cite{Gr99}
\begin{equation}\label{e:smooth6}
(t\Ds_{\Pp_i}^2-\la\If)\ii - (t\Ds^2-\la\If)_{+}\ii =
\Kk_{\la}S_{i,\la}\Pp_i \g_1 (t\Ds^2-\la\If)_{+}\ii \ ,
\end{equation}
where $(t\Ds^2-\la\If)_{+}\ii$ is the restriction to $X$ of the
inverse of the invertible double operator of $t\Ds^2-\la\If$,
$\Kk_{\la}$ is the classical vertical Poisson operator for
$t\Ds^2-\la\If$,  and $S_{i,\la}\in \Aa(B,\Psi^0(\Ef^{'}))$ is a
vertical right inverse to $\Pp\circ P(\Ds^2)$. Hence the relative
resolvent is a vertical sgo. Since $\Pp\circ P(\Ds^2)^{\perp}\in
\Aa(B,\Psi^{-\o}(\Ef^{'}))$, and
$$\Pp\g_1 (t\Ds^2-\la\If)_{+}\ii   =
\Pp(\Ds^2))^{\perp}\g_1 (t\Ds^2-\la\If)_{+}\ii $$ it is of order
$-\o$, and \eqref{e:smooth2} follows.

It is enough to prove \eqref{e:smooth3} locally over $U\subset B$.
Since $\na^{\pi_*(\Ef|\Pp_i)} = \Ps_i\cdot \na^{\pi_*(\Ef)}\cdot
\Ps_i$ we have
\begin{equation*}
\na^{\pi_*(\Ef|\Pp_i)}_{|U}  = \Ps^{z}_i\cdot d_U \cdot \Ps^{z}_i
+ \sum_{j=1}^{\dim B} \Ps^{z}_i\cdot \Dd_{j}^z \cdot \Ps^{z}_i
dz_j \ .
\end{equation*}
$\Dd_{j}^z $ are first-order differential operators on
$\G(X_{z_0},E_{z_0})$. The $\Ps^{z}_i$ are defined with respect to
elements of the fixed Grassmannian $\Gr_{\o}(\dd_{z_0})$ depending
smoothly on $z$ and defining $\Ps^{z}_i : U \too \Aa(U)\otimes
\G(\dd X_{z_0},\End(E^0_{z_0}))$. Since $\Ps^{z}_i = I -\Ks^{z}_i
\Pp^{z}_i\g$
\begin{equation*}
d \Ps^{z}_i  = -(d\Ks^{z}_i) \Pp^{z}_i\g -\Ks^{z}_i (d\Pp^{z}_i)\g
\in \Aa^1(U,\Psi^{-\o}(E_{z_0},E_{z_0}^0)) \ ,
\end{equation*}
is a smooth family of (generalized) smoothing operators. For,
$$d\Ks^{z}_i = \chi(u)\sum_{j=0}^{\dim B} \frac{\dd}{\dd
z_j}e^{-u(\dd^z)^2} dz_i$$ and by Duhamel's formula $\dd_j
e^{-u(\dd^z)^2}\in\G(U,\Psi^{-\o}(E^0_{z_0}\oplus E^0_{z_0})))$,
while for fixed $\Ps_0 = \Ps(\Pp_0)$ we have $$d\Ps^{z}_i =
d(\Ps^{z}_i - \Ps_0) \in\Aa(U,\Psi^{-\o}(\Ef,\Ef^{'}))\ , $$ since
$\Ps^{z}_i - \Ps_0\in\Aa(U, \Psi^{-\o}(\Ef,\Ef^{'}))$.

The curvature can be computed by a local formula similar to that
in \lemref{lem:relativelyboundarysmooth}, consisting of elements
in $\Aa(U,\Psi^{-\o}(\Ef,\Ef^{'}))$ and non-smoothing terms but
whose difference in $\Rs_1 - \Rs_2$ are in
$\Aa(U,\Psi^{-\o}(\Ef,\Ef^{'}))$. The details are similar to
\lemref{lem:relativelyboundarysmooth}. Equation \eqref{e:smooth4}
follows by an analogous argument.
\end{proof}

\noi Proof of \propref{p:heat}:
\begin{proof}
Integrating by parts in \eqref{e:heatP} we have
\begin{equation}\label{e:heatbyparts}
 e^{-\Ff_{t,\Pp}} :=
\frac{i}{2\pi}\int_{\Cc}e^{-\la} \dd_{\la}^m(\Ff_{t,\Pp} -
\la\If)\ii \ d\la \ .
\end{equation}
By \eqref{e:nilpotent}, $\dd_{\la}^m(\Ff_{t,\Pp} - \la\If)\ii$ is
a sum of terms of the form
$$\dd_{\la}^{j_1}G_{t,\la}\ii F  \dd_{\la}^{j_2}G_{t,\la}
 \ldots F \dd_{\la}^{j_k}G_{t,\la} \ ,$$
where $G_{t,\la}\ii = (t\Ds_{\Pp}^2-\la\If)$, and $F\in
\Aa^2(B,\Dd(\Ef))$ is a 2-form with first-order differential
operator coefficients. From sgo theory on manifolds with boundary
\cite{Gr96,Gr99} and equations \eqref{e:nilpotent} and
\eqref{e:smooth6}, and \lemref{lem:commutatortrace} [1],
\lemref{lem:relativelysmooth}, we find that
$\dd_{\la}^m(\Ff_{t,\Pp} - \la\If)\ii$ is a vertical sgo
consisting of terms of order $-2m-2-k$ with $k\geq 0$. Since this
is so for any $m$, from \eqref{e:heatbyparts} $ e^{-\Ff_{t,\Pp}}$
is a vertical sgo of order $-\o$. In particular it is trace class
with a generalized smoothing kernel. A detailed study of form
degree zero term is given in $\S$4.2 of \cite{Gr96}.

Since $(t\Ff_{\Pp} - \la\If)\ii$ has range in
$\Aa(B,\pi_*(\Ef|\Pp^2))$, then  $\Ff(t\Ff_{\Pp} - \la\If)\ii =
\Ff_{\Pp}(t\Ff_{\Pp} - \la\If)\ii$ and hence
\begin{eqnarray*}\label{e:heateqn}
 \dd_t(t\Ff_{\Pp} - \la\If)\ii & = & - t (t\Ff_{\Pp} - \la\If)\ii
 \Ff_{\Pp} (t\Ff_{\Pp} - \la\If)\ii \\
 &= & -\Ff_{\Pp} (t\Ff_{\Pp} - \la\If)^{-2} \ .
\end{eqnarray*}
Hence
$$\dd_t  e^{-t\Ff_{\Pp}} = -\int_{C}e^{-\la}\Ff_{\Pp}
(t\Ff_{\Pp} - \la\If)^{-2} \, d\la = -\Ff_{\Pp} e^{-t\Ff_{\Pp}} $$
integrating by parts for the second equality. Since $\Ff_{\Pp} =
\Af_{\Pp}^2$ on $\Aa(B,\pi_*(\Ef|\Pp^2))$, this proves
\eqref{e:heatequation}.  A similar argument shows that
$\|(e^{-t\Ff_{\Pp}} - e^{-t})\psi \|_s \to 0$ as $t\to 0$ in each
Sobolev norm and hence the boundary condition at $t=0$ holds. This
completes the proof.
\end{proof}

We can now define the {\em Chern character form} of the
superconnection $\Af_{t,\Pp}$ adapted to the family of APS-type
boundary problems $\Ds_{\Pp}$ to be the differential form in
$\Aa^{{\rm even}}(B)$

\begin{equation}\label{e:cherncharacter}
 \ch(\Af_{t,\Pp}) = \str \left(e^{-\Ff_{t,\Pp}}\right) \ .
\end{equation}
\vskip 2mm \noi For later use, observe that for large $m$,
$\dd_{\la}^m(\Ff_{t,\Pp} - \la\If)\ii$ is a trace class vertical
operator and hence from \eqref{e:nilpotent}, \eqref{e:heatbyparts}
\begin{eqnarray}
&  & \hskip 20mm \ch(\Af_{t,\Pp})  =
\frac{i}{2\pi}\int_{\Cc}e^{-\la}
 \str(\dd_{\la}^m(\Ff_{t,\Pp} - \la\If)\ii) \ d\la \label{e:cherncharacter2}
 \\  & = & \sum_{k=0}^{\dim B }(-1)^k \frac{i}{2\pi}
 \int_{\Cc}  e^{-\la}  \str\left[\dd_{\la}^m\left(
(t\Ds_{\Pp}^2-\la\If)\ii \left((\Rs^{\pi_*(\Ef | \Pp)}+
t^{1/2}\na^{\pi_*(\Ef | \Pp)}\Ds_{\Pp})(t\Ds_{\Pp}^2 -
\la\If)\ii\right)^k\right) \right]\ d\la \ . \nonumber
\end{eqnarray}

\vskip 5mm

\subsection{The Transgression Formula}
Because we are considering Dirac operators over a manifold with
boundary with domain restricted by the Grassmann section the
supertrace still vanishes on commutators of APS-type boundary
problems with smoothing operators which have the correct range.
This gives us a direct analogue of the transgression formula for
closed manifolds \cite{BeVe90,BeGeVe92,Bi86}.

\vskip 1mm

\begin{thm}\label{t:transgression}  \noi \newline \vskip 1mm
\noi {\rm [1]} The differential form $\ch(\Af_{t,\Pp})$
on $B$ is closed.\\[1mm]

\noi {\rm [2]} The cohomology class defined by  $\ch(\Af_{t,\Pp})$
is independent of the parameter $t$. One has
\begin{equation}\label{e:transgression1}
\frac{d}{dt}\ch(\Af_{t,\Pp}) = -d_B \str\left(\dot{\Af}_{t,\Pp}\,
e^{-\Ff_{t,\Pp}}\right) \ ,
\end{equation}
where $\dot{\Af}_{t,\Pp} = (d/dt)\Af_{t,\Pp}$. Equivalently, for
$0 < t < T < \o$,
\begin{equation}\label{e:transgression2}
\ch(\Af_{T,\Pp}) - \ch(\Af_{t,\Pp}) = -d_B\int_{t}^{T}
\str\left(\dot{\Af}_{s,\Pp}\, e^{-\Ff_{s,\Pp}}\right) ds \ .
\end{equation}
More generally, $\ch(\Af_{\sigma,\Pp})$ is invariant for a general
family $\Af_{\sigma,\Pp}$ for variations $\sigma$ that leave
unchanged the symbol of the boundary family $\Ds_Y$ of Dirac
operators, but it is not a general homotopy invariant.
\end{thm}

First, we need the following property.
\begin{lem}\label{lem:tracecommutes}
Let $X$ be a compact manifold with boundary and let $D$ be an
elliptic operator of Dirac type acting on $\G(X,E)$. Let $P\in
\Gr_{\o}(D_{|\dd X})$ be an APS-type boundary condition for $D$.
Then if there exists a generalized smoothing operator $\Lambda_P$
on $\G(X,E)$ such that for all non-negative integers
$k=0,1,\ldots$
\begin{equation}\label{e:Prange}
P\g D^k \Lambda_P = 0 \ ,
\end{equation}
then for any differential operator $A$ acting on $\G(X,E)$
\begin{equation}\label{e:commutes}
\str\left([A,\Lambda_P]\right) = 0 \ .
\end{equation}
\end{lem}
\begin{proof}
The reason for this property is that $D_P$ has both a left and a
right parametrix. Let us briefly explain this further. Suppose to
begin with that $D_P$ is invertible. Then one has
\begin{equation}\label{e:parametrix}
D \cdot D_P\ii = I, \hskip 10mm   D_P\ii\cdot D = I - \verb"K"
S(P)\ii P\g \ ,
\end{equation}
where $\verb"K"$ is the Poisson operator for $D$. The second
identity means that $D$ has no left inverse, which is the source
of Green's theorem and the non-commutativity of the trace in
general for compositions of differential and smoothing operators.
On the other hand, $I = D\cdot D_P\ii = D_P \cdot D_P\ii$, so
$D_P\ii$ is a right inverse for $D_P$. Further, since $\verb"K"
S(P)\ii P\g = 0$ on $\dom(D_P)$ then $D_P\ii\cdot D_P =
I_{|\dom(D_P)}$, and so it is a two-sided inverse for $D_P$. The
sufficient condition \eqref{e:Prange} allows this to be repeated
this with all higher powers $D^k_P$ (with domain $\dom(D^k_P) =
\{\psi \in \G(X,E) \ | \ P\g D^i\psi = 0, \, 0\leq i \leq k-1 \}$)
and hence sufficiently bound the Sobolev norm. Precisely, choosing
$N > \ord(A)$ we have that $AD_{P}^{-N}$ is bounded in the
operator norm on $\G(X,E)$ defined by metrics on $E$ and $TX$ (--
more is true: by taking $N$ sufficiently large we can ensure a
$C^r$ kernel). Now we have
\begin{equation}\label{e:parametrix2}
I_{|\dom(D_P^{N})} = D_{P}^{-N} D^N_{|\dom(D_P^{N})} = D_{P}^{-N}
D_P^N \ ,
\end{equation}
while taking $k=N-1$ in \eqref{e:Prange} we have
$\ran(\Lambda_P)\subset \dom(D_P^{N})$, and hence
$$D_{P}^{-N} D^N \Lambda_P = \Lambda_P \ . $$ By
\lemref{lem:commutatortrace} [1], $D^N \Lambda_P$ is a generalized
smoothing operator, and so we can follow the argument of
\cite{BeGeVe92} Prop(2.48) to obtain
$$\str(A\Lambda_P) = \str\left((AD_{P}^{-N}) D^N \Lambda_P\right)
= \str\left((D^N \Lambda_P A)D_{P}^{-N} \right) $$ $$ =
\str\left((D_{P}^{-N}D^N) \Lambda_P A \right) = \str(\Lambda_P
A),$$ using the commutativity of the trace for combinations of
bounded and smoothing operators. If $D_P$ is not invertible, then
the argument obviously extends by replacing the inverse by a
parametrix $Q_P$ (so $\ran(Q_P)\subset \dom(D_P^{N})$ and $Q_P D_P
- I$ and $I - DQ_P$ are generalized smoothing operators).
\end{proof}

\noi We prove \thmref{t:transgression} for {\em any}
superconnection adapted to $\Ds_{\Pp}$:
\begin{proof}
The essential point is to show that for $k = 0,1,\ldots$
\begin{equation}\label{e:Prange2}
\Pp\g \Ds^{2k}e^{-\Ff_{t,\Pp}}  = 0 \ .
\end{equation}
We have $t\Ds^2 (t\Ds^2 - \la\If)\ii = \If + \la (t\Ds^2 -
\la\If)\ii$ on $\Aa(B,\pi_*(\Ef))$ and hence
\begin{equation*}
\Ds^2 e^{-t\Ds_{\Pp}^2} = \frac{i}{2\pi}\int_{\Cc}t\ii \la
e^{-\la} (t\Ds_{\Pp}^2 - \la\If)\ii  \ d\la \ ,
\end{equation*}
and iteratively
\begin{equation*}
\Ds^{2k} e^{-t\Ds_{\Pp}^2} = \frac{i}{2\pi}\int_{\Cc}t^{-k}\la^k
e^{-\la} (t\Ds_{\Pp}^2 - \la\If)\ii  \ d\la \ .
\end{equation*}
Hence, by \eqref{e:rangeinP}, for any $t>0$ and non-negative
integer $k$
\begin{equation}\label{e:Prange3}
\Pp\g \Ds^{2k}e^{-t\Ds_{\Pp}^2} =
\frac{i}{2\pi}\int_{\Cc}t^{-k}\la^k e^{-\la} \Pp\g  (t\Ds_{\Pp}^2
- \la\If)\ii \ d\la  = 0 \ .
\end{equation}
On the other hand Duhamel's formula gives
\begin{equation}\label{e:Prange4}
e^{-t\Ff_{\Pp}}  = e^{-t\Ds_{\Pp}^2} + \sum_{k=1}^{\dim B}(-t)^k
I_{k,t} \ ,
\end{equation}
where
\begin{equation*}
I_{k,t}  = \int_{\sigma_k} e^{-s_0 \Ds_{\Pp}^2}(\Ff_{\Pp}-
\Ds_{\Pp}^2)e^{-s_1 \Ds_{\Pp}^2}(\Ff_{\Pp}- \Ds_{\Pp}^2) \ldots
e^{-s_{k-1} \Ds_{\Pp}^2}(\Ff_{\Pp}- \Ds_{\Pp}^2)e^{-s_k
\Ds_{\Pp}^2} \ d\sigma_k \ ,
\end{equation*}
and $\sigma_k$ is the standard $k$-simplex. From
\eqref{e:Prange3}, \eqref{e:Prange4} we therefore have $\Pp\g
\Ds^{2k}e^{-t\Ff_{\Pp}}  = 0$, and \eqref{e:Prange2} then follows
from \eqref{e:heatdeltat}. (Equivalently, one may use
\eqref{e:Prange3} and \eqref{e:nilpotent}.)

Now let $\Af$ be any superconnection on $\pi_*(\Ef)$. Locally,
$\Af_{|U} = d_U + \sum_I D_{I}^z dz_I$, while $e^{-\Ff_{t,\Pp}}  =
\sum_J K(\Pp)_{J}^z dz_J$, with $K(\Pp)_{J}^z$ a generalized
smoothing operator satisfying
\begin{equation}\label{e:Prange5}
\Pp\g D^{2k} K(\Pp)_J^z = 0 \ .
\end{equation}
Hence by \lemref{lem:tracecommutes}
\begin{equation}\label{e:Prange6}
\str([\sum_I D_{I}^z dz_I,\sum_J K(\Pp)_{J}^z dz_J]) = \sum_{I,J}
\str([D_{I}^z, K(\Pp)_{J}^z])\, dz_I\wedge dz_J = 0  \ .
\end{equation}
Thus $\str([\Af_{|U},e^{-\Ff_{t,\Pp}} ]) =
\str([d_U,e^{-\Ff_{t,\Pp}} ]) = d_U \str(e^{-\Ff_{t,\Pp}})$,
implying the global formula
\begin{equation}\label{e:dstr}
d_B \str(e^{-\Ff_{t,\Pp}}) = \str([\Af,e^{-\Ff_{t,\Pp}} ])  \ .
\end{equation}
To see the right-side vanishes we show it is equal to
$\str([\Af_{t,\Pp},e^{-\Ff_{t,\Pp}} ])$ -- the potential
obstruction to this comes from the second fundamental form
$II_{\Pp} = \nE_{|\pi_{*}(\Ef|\Pp)} - \na^{\pi_*(\Ef | \Pp)}$. It
is enough to establish this locally, and by the above discussion
it is sufficient to prove:
\begin{lem}\label{lem:sff}
$$\str([\Ps\cdot d \cdot \Ps, e^{-\Ff_{t,\Pp}}]) = \str([d,
e^{-\Ff_{t,\Pp}} ]) \ .$$
\end{lem}
\begin{proof}
For any $\a\in\Aa(B,\Psi^{-\o}(\Ef,\Ef^{'}))$ one has $\str(\a) =
\str(\Ps\a\Ps) + \str(\Ps\pp\a\Ps\pp).$ Hence, since
$e^{-\Ff_{t,\Pp}} = \Ps e^{-\Ff_{t,\Pp}} $ is even, with
$d^{\pi_{*}(\Ef|\Pp)}_U := \Ps\cdot d_U\cdot\Ps$,
\begin{eqnarray*}
\str([d_U,e^{-\Ff_{t,\Pp}} ]) & = & \str(d_U\cdot\Ps\cdot
e^{-\Ff_{t,\Pp}} - \Ps\cdot e^{-\Ff_{t,\Pp}}\cdot d_U ) \\
& = & \str(\Ps\cdot d_U\cdot\Ps e^{-\Ff_{t,\Pp}}\cdot\Ps -
\Ps\cdot e^{-\Ff_{t,\Pp}}(P+P\pp)d_U\cdot \Ps  \\
&  &  + \, \Ps\pp\cdot d_U\cdot\Ps e^{-\Ff_{t,\Pp}}\cdot\Ps\pp
- \Ps\pp\Ps\cdot e^{-\Ff_{t,\Pp}}(P+P\pp)d_U\cdot \Ps\pp) \\
& = & \str([d^{\pi_{*}(\Ef|\Pp)}_U, e^{-\Ff_{t,\Pp}}])   +
\str([\Ps\pp d\Ps,\Ps e^{-\Ff_{t,\Pp}}\Ps\pp])     \\
& = & \str([d^{\pi_{*}(\Ef|\Pp)}_U, e^{-\Ff_{t,\Pp}}]) \ ,
\end{eqnarray*}
since $\str$ vanishes on supercommutators in
$\Aa(B,\End(\pi_*(\Ef)))$.
\end{proof}
\noi We conclude
\begin{equation}\label{e:dstrP}
d_B \str(e^{-\Ff_{t,\Pp}}) = \str([\Af_{t,\Pp},e^{-\Ff_{t,\Pp}} ])
= \frac{i}{2\pi}\int_{\Cc} e^{-\la} \, \str\left((\Ff_{t,\Pp} -
\la\If)\ii [\Af_{t,\Pp}, \Ff_{t,\Pp}] (\Ff_{t,\Pp} -
\la\If)\ii\right) \ d\la
\end{equation}
which vanishes since $\Ff_{t,\Pp}= \Af_{t,\Pp}^2$ on the
$L^2$-dense subbundle $\pi_*(\Ef|\Pp^2)$ of $\pi_*(\Ef)$.

To prove [2], from \eqref{e:heatbyparts} we have that
\begin{eqnarray*}
\frac{d}{dt} e^{-\Ff_{t,\Pp}} & = & \frac{i}{2\pi}\int_{\Cc}
e^{-\la} \,  \dd_{t}\dd_{\la}^m(\Ff_{t,\Pp} - \la\If)\ii  \ d\la \\
& = & -\frac{i}{2\pi} m! \int_{\Cc} e^{-\la} \, (\Ff_{t,\Pp} -
\la\If)^{-m-1}
 \dd_{t}((\Ff_{t,\Pp} - \la\If)^{m+1}) (\Ff_{t,\Pp} - \la\If)^{-m-1} \ d\la
\ ,
\end{eqnarray*}
noting that the domains match-up; that is, $\Ff_{t}^m (\Ff_{t,\Pp}
- \la\If)^{-m-1} = \Ff^m_{t,\Pp}(\Ff_{t,\Pp} - \la\If)^{-m-1}$.
Similarly, since the variation is interior $\Pp\g\dot{\Ff}_{t} =
0$, so $\dot{\Ff}_{t} = \dot{\Ff}_{t,\Pp}$. Hence
\begin{equation*}
\dd_{t}(\Ff_{t,\Pp} - \la\If)^{m+1} = \sum_{k=0}^m (\Ff_{t,\Pp} -
\la\If)^{m-k} \dot{\Ff}_{t,\Pp} (\Ff_{t,\Pp} - \la\If)^k \ ,
\end{equation*}
and so
\begin{equation*}
\frac{d}{dt} e^{-\Ff_{t,\Pp}}  = -\frac{i}{2\pi} m! \sum_{k=0}^m
\int_{\Cc} e^{-\la} \, (\Ff_{t,\Pp} - \la\If)^{-k-1}
\dot{\Ff}_{t,\Pp} (\Ff_{t,\Pp} - \la\If)^{k-m-1} \ d\la \ .
\end{equation*}
For $m > 4n$, at least one of the vertical operators $(\Ff_{t,\Pp}
- \la\If)^{-k-1}$ or $\dot{\Ff}_{t,\Pp} (\Ff_{t,\Pp} -
\la\If)^{k-m-1}$ is trace class, while the other is bounded
 in the operator norm.  Hence we can swap the order of the
 operators in the trace to obtain
\begin{eqnarray*}
\frac{d}{dt} \str(e^{-\Ff_{t,\Pp}}) & = & -\frac{i}{2\pi} m!
\sum_{k=0}^m \int_{\Cc} e^{-\la} \, \str(\dot{\Ff}_{t,\Pp}
(\Ff_{t,\Pp} - \la\If)^{-m-2}) \ d\la \\
& = & -\frac{i}{2\pi} m! \sum_{k=0}^m \int_{\Cc} e^{-\la} \,
\str(\dot{\Ff}_{t,\Pp} \, \dd_{\la}^{m+1}(\Ff_{t,\Pp} -
\la\If)\ii) \
d\la \\
& = & - \,\str(\dot{\Ff}_{t,\Pp}\, e^{-\Ff_{t,\Pp}}) \ .
\end{eqnarray*}
Since $\Ff_{t,\Pp}= \Af_{t,\Pp}^2$ on $\ran(e^{-\Ff_{t,\Pp}})$, we
have
\begin{equation}\label{e:chvariation}
\dot{\Ff}_{t,\Pp} \, e^{-\Ff_{t,\Pp}} =
\Af_{t,\Pp}\dot{\Af}_{t,\Pp} \, e^{-\Ff_{t,\Pp}} +
\dot{\Af}_{t,\Pp} \Af_{t,\Pp} \, e^{-\Ff_{t,\Pp}} \ .
\end{equation}
Since $\Af_{t,\Pp}(\Ff_{t,\Pp} - \la\If )\ii = (\Ff_{t,\Pp} -
\la\If )\ii \Af_{t,\Pp}$ in the $L^2$-dense subbundle
$\pi_*(\Ef|\Pp)$ we can commute $\Af_{t,\Pp}$ and
$e^{-\Ff_{t,\Pp}}$ in the supertrace, and hence from
\eqref{e:chvariation}
\begin{equation}\label{e:d-strace}
\str(\dot{\Ff}_{t,\Pp}\, e^{-\Ff_{t,\Pp}}) =
\str([\Af_t,\dot{\Af}_{t,\Pp}\, e^{-\Ff_{t,\Pp}}]) \ .
\end{equation}
But $\dot{\Af}_{t,\Pp}\, e^{-\Ff_{t,\Pp}} = \frac{1}{2}t^{-1/2}\Ds
\, e^{-\Ff_{t,\Pp}}$ and hence, by the argument of part [1]
satisfies
$$\Pp\g \Ds^{k}\dot{\Af}_{t,\Pp}\, e^{-\Ff_{t,\Pp}}  = 0 \ . $$ We
therefore conclude, as in \eqref{e:dstr}, that \eqref{e:d-strace}
is equal to $d_B \str(\dot{\Af}_{t,\Pp} \, e^{-\Ff_{t,\Pp}})$.

The above argument depends on the variation leaving the boundary
condition invariant, but otherwise holds for a general 1-parameter
family of superconnections $\Af_t$. In particular, given two
superconnections $\Af_1, \Af_2$ it holds for $\Af_t = t\Af_1 +
(1-t)\Af_2$ adapted to $\Ds_{\Pp}$. Hence the statement on
interior homotopy invariance follows. This completes the proof.
\end{proof}

\vskip 2mm

\subsection{Computing the Limiting Values of {\bf $\ch(\Af_{t,\Pp_1}) -
\ch(\Af_{t,\Pp_2})$.}}\vskip 1mm

Next we study the relative form $\ch(\Af_{t,\Pp_1}) -
\ch(\Af_{t,\Pp_2})$ as $ t\to 0 $ and $t \to \o$.

Notice that $\ch(\Af_{t,\Pp_i})$ do not have limits as $t\to 0$;
for closed manifolds the content of the local families index
theorem is that such a limit does exist for the Bismut
superconnection.

First, let us point out that the limit
\begin{equation}\label{e:relchlimit}
\lim_{t\to 0}(\ch(\Af_{t,\Pp_1}) - \ch(\Af_{t,\Pp_2}))
\end{equation}
is not (quite) the generalized relative Eta-form
\begin{equation}\label{e:genetaform}
    \eta^{[M]}(\Ps_1,\Ps_2) = \str(e^{-\Rs_1} - e^{-\Rs_2}) \ ,
\end{equation}
but also includes a regularized trace term. The right-side of
\eqref{e:genetaform} exists for the same reason as in
finite-dimensions explained in the introduction. More precisely:
\begin{prop}\label{p:genetaform}
$$e^{-\Rs_1} - e^{-\Rs_2} \in \Aa(B,\Psi^{-\o}(\Ef,\Ef^{'})) \ .$$
One has
\begin{equation}\label{e:genetaformseries}
\eta^{[M]}(\Ps_1,\Ps_2) = \sum_{k\geq
1}(-1)^k\frac{1}{k!}\,\str\left(\Rs_{1}^k - \Rs_{2}^k\right)
\end{equation}
\end{prop}
\begin{proof}
The essential point is that since the curvature $\Rs_i =
(\Ps_i\cdot\nE\cdot\Ps_i)^2$ is a 2-form, with restricted
differential operator coefficients, it is nilpotent. Thus
\begin{equation}\label{e:relRresolvents}
(\Rs_1 - \la\If)\ii - (\Rs_2 - \la\If)\ii = \sum_{k\geq
1}(-\la)^{-k-1}(\Rs_{1}^k - \Rs_{2}^k)
\end{equation}
is a finite sum, and hence so is
\begin{equation*} e^{-\Rs_1} - e^{-\Rs_2}  =
\frac{i}{2\pi}\int_{c}e^{-\la} ((\Rs_1 - \la\If)\ii - (\Rs_2 -
\la\If)\ii) \ d\la \ ,
\end{equation*}
where $c$ can be any simple closed contour around the origin. The
first assertion now follows from \eqref{e:smooth3}, and the
formula \eqref{e:genetaformseries} from \eqref{e:relRresolvents}.
\end{proof}

Notice that there is no order zero term in
\eqref{e:genetaformseries}, \eqref{e:relRresolvents}, in contrast
to \eqref{e:inclusion}, since here we work on the full space of
sections -- now the order zero contribution (the pointwise index)
comes from the regularized zeta term.

Consider the limit \eqref{e:relchlimit}.

\begin{lem}\label{lem:relsmooth}
\begin{equation}\label{e:relsmooth}
(\Ff_{t,\Pp_1} - \la\If )\ii -(\Ff_{t,\Pp_2} - \la\If )\ii
\in\Aa(B,\Psi^{-\o}(\Ef,\Ef^{'})) \ .
\end{equation}
\end{lem}
\begin{proof}
Let $Z_{t,i} = t\Rs_i + t^{1/2}\na^{\pi_*(\Ef|\Pp)}D_i$. From
\eqref{e:nilpotent} we have
\begin{equation*}
(\Ff_{t,\Pp_1} - \la\If )\ii -(\Ff_{t,\Pp_2} - \la\If )\ii  =
\end{equation*}
\begin{equation}\label{e:relexpand1}
\left\{(t\Ds^2_{\Pp_1} - \la\If )\ii -(t\Ds^2_{\Pp_2} - \la\If
)\ii\right\} \sum_{k=0}^{\dim B}(-1)^k (Z_{t,1}(t\Ds_{\Pp_1}^2 -
\la\If)\ii)^k
\end{equation}
\begin{equation}\label{e:relexpand2}
+ \ (t\Ds_{\Pp}^2 - \la\If)\ii \sum_{k=0}^{\dim B}(-1)^k
\left\{Z_{t,1}(t\Ds_{\Pp_1}^2 - \la\If)\ii)^k -
(Z_{t,2}(t\Ds_{\Pp_2}^2 - \la\If)\ii)^k \right\} \ .
\end{equation}
From \lemref{lem:relativelysmooth} the first factor in
\eqref{e:relexpand1} is smooth family of smoothing operators.
Since the second factor consists of terms which are smooth
families of first-order differential operator valued forms ,
coming from $Z_{t,1}$, and sgo operators of order $-1$, coming
from $(t\Ds_{\Pp_1}^2 - \la\If)\ii$ then the composition is in
$\Aa(B,\Psi^{-\o}(\Ef,\Ef^{'}))$ by the sgo calculus and
\lemref{lem:commutatortrace}[1].

Similarly, \eqref{e:relexpand2} is in
$\Aa(B,\Psi^{-\o}(\Ef,\Ef^{'}))$ since the first factor is a
vertical family sgo operators of order $-1$, while it follows from
\lemref{lem:relativelysmooth} and \lemref{lem:commutatortrace}[1]
that the second factor is a smooth family of smoothing operators.
This can be checked term by term. When $k=0$ this is trivial, the
$k=1$ term can be expanded
$$(Z_{t,1} -Z_{t,2})(t\Ds^2_{\Pp_1} - \la\If )\ii  +
Z_{t,2}((t\Ds^2_{\Pp_1} - \la\If )\ii -(t\Ds^2_{\Pp_2} - \la\If
)\ii) \ ,$$ which is in $\Aa(B,\Psi^{-\o}(\Ef,\Ef^{'}))$, as
above. The general term is given by a similar formula, and we
hence reach the conclusion.
\end{proof}

\noi By the linearity of the supertrace we therefore have:

\begin{cor}\label{c:relch}
For $t>0$
\begin{equation}\label{e:relchcharacter3128}
\ch(\Af_{t,\Pp_1}) - \ch(\Af_{t,\Pp_2})
\end{equation}
$$ = \frac{i}{2\pi}\int_{\Cc}e^{-\la}\, \str\left((\Rs_1 +
t^{1/2}\na^1 \Ds_{\Pp_1}  + t\Ds_{\Pp_1}^2 - \la\If)\ii  - (\Rs_2
+ t^{1/2}\na^2 \Ds_{\Pp_2}  + t\Ds_{\Pp_2}^2 - \la\If)\ii\right) \
d\la $$
$$ = \sum_{k=0}^{\dim B} \frac{i}{2\pi}\int_{\Cc}e^{-\la}\,
\str\left((t\Ds^2_{\Pp_1} - \la\If )\ii (Z_{t,1}(t\Ds_{\Pp_1}^2 -
\la\If)\ii)^k - (t\Ds^2_{\Pp_2} - \la\If )\ii
(Z_{t,2}(t\Ds_{\Pp_2}^2 - \la\If)\ii)^k \right) \ d\la \ .$$
\end{cor}

The crucial fact in the proof of Theorem (II) is that the limit
\eqref{e:relchlimit} has no contribution coming from the term
$t^{1/2}\na^i \Ds_{\Pp_i}$. This eliminates complicated Wodzicki
residue terms in the regularized limit $\LIM_{t\to 0}
\ch(\Af_{t})$--see \cite{Sc02} for precise formulas on closed
manifolds.

\begin{prop}\label{p:limitexists} The limit
 $\lim_{t\to 0}(\ch(\Af_{t,\Pp_1}) - \ch(\Af_{t,\Pp_2}))$
exists. One has
\begin{equation}\label{e:simplerlimit}
\lim_{t\to 0}(\ch(\Af_{t,\Pp_1}) - \ch(\Af_{t,\Pp_2}))
\end{equation}
$$ = \lim_{t\to 0} \frac{i}{2\pi}\int_{\Cc}e^{-\la}\, \str\left((\Rs_1 +
 t\Ds_{\Pp_1}^2 - \la\If)\ii  - (\Rs_2
 + t\Ds_{\Pp_2}^2 - \la\If)\ii\right)
\ d\la $$ in each $C^l$ norm.
\end{prop}
\begin{proof} Notice that the right-side of \eqref{e:relchcharacter3128}
at $t=0$ is just $\eta(\Ps_1,\Ps_2)$, but this is not the limit
\eqref{e:simplerlimit} --- $\ch(\Af_{t,\Pp_i})$ do not have limits
as $t\to 0$.

The supertrace of the relative vertical resolvent $(\Ff_{t,\Pp_1}
- \la\If )\ii -(\Ff_{t,\Pp_2} - \la\If )\ii$ is equal to the sum
of the supertraces of the expressions in \eqref{e:relexpand1} and
\eqref{e:relexpand2}. For the first of these we have using
\eqref{e:smooth6}, with $\| .\|, \|.\|_1$ respectively the
operator and trace norms on $\Aa(B,\End(\pi_*(\Ef)))$,
\begin{equation}\label{e:str1}
|\str\left(\{(t\Ds^2_{\Pp_1} - \la\If )\ii -(t\Ds^2_{\Pp_2} -
\la\If )\ii\}\right. \sum_{k=0}^{\dim B}\left.(-1)^k
(Z_{t,1}(t\Ds_{\Pp_1}^2 - \la\If)\ii)^k\right) |
\end{equation}
\begin{equation*}
 \leq  t\ii \| (\Ds^2_{\Pp_1} - t\ii\la\If )\ii
-(\Ds^2_{\Pp_2} - t\ii\la\If )\ii\|_1 \sum_{k=0}^{\dim B} t\ii
\|(Z_{t,1}(\Ds_{\Pp_1}^2 - t\ii\la\If)\ii)\|^k
\end{equation*}
\begin{equation*}
 \leq  t\ii (\|\Kk_{t\ii\la}S_{1,t\ii\la}\|.\|\Pp_1
 P(\Ds^2)^{\perp}\g_1\|_1
 + \|\Kk_{t\ii\la}S_{2,t\ii\la}\|.\|\Pp_2P (\Ds^2))^{\perp}\g_1\|_1)
  \|(\Ds^2-t\ii\la\If)_{+}\ii \|
\end{equation*}
\begin{equation}\label{e:firstestimate}
\times \sum_{k=0}^{\dim B} t\ii \|(Z_{t,1}(\Ds_{\Pp_1}^2 -
t\ii\la\If)\ii)\|^k \ .
\end{equation}
From standard estimates in elliptic theory on manifolds with
boundary \cite{Gr96} and closed manifolds respectively, we have as
$t\too 0$
\begin{equation}\label{e:term1}
\|\Kk_{t^{-1}\la}S_{1,t\ii\la}\| = O(1), \hskip 10mm
\|(\Ds^2-t\ii\la\If)_{+}\ii \| = O(t|\la|\ii) \ .
\end{equation}
To estimate the sum, let $Q_{\Pp_1}\in\Aa(B,\End(\pi_*(\Ef)))$ be
a vertical parametrix for $\Ds_{\Pp_1}$, so that
$Q_{\Pp_1}\Ds_{\Pp_1} = \If + S$ with
$S\in\Aa(B,\Psi^{-\o}(\Ef,\Ef^{'}))$. Then
\begin{equation*}
 t\ii Z_{t,1}(\Ds_{\Pp_1}^2 -
t\ii\la\If)\ii  = t\ii Z_{t,1}Q_{\Pp_1}\Ds(\Ds_{\Pp_1}^2 -
t\ii\la\If)\ii  - t\ii Z_{t,1}S(\Ds_{\Pp_1}^2 - t\ii\la\If)\ii
\end{equation*}
\begin{equation*}
= \frac{t}{2}\ii Z_{t,1}Q_{\Pp_1}\left((\Ds_{\Pp_1} -
t^{-1/2}\la^{1/2}\If)\ii + (\Ds_{\Pp_1} +
t^{-1/2}\la^{1/2}\If)\ii\right) - t\ii Z_{t,1}S(\Ds_{\Pp_1}^2 -
t\ii\la\If)\ii \ ,
\end{equation*}
since $\Dd_{\Pp_1}$ is a family of self-adjoint operators. Since
$Z_{t,1}$ is a vertical family of first-order differential
operators it follows that $Z_{t,1}Q_{\Pp_1}$ is bounded in the
operator norm. Hence
\begin{equation}\label{e:suminequality}
\|(Z_{t,1}(\Ds_{\Pp_1}^2 - t\ii\la\If)\ii)\| \leq \frac{t}{2}\ii
\|Z_{t,1}Q_{\Pp_1}\|(\|(\Ds_{\Pp_1} - t^{-1/2}\la^{1/2}\If)\ii)\|
+ \|(\Ds_{\Pp_1} + t^{-1/2}\la^{1/2}\If)\ii\|)
\end{equation}
$$- t\ii
\|Z_{t,1}S\|\|(\Ds_{\Pp_1}^2 - t\ii\la\If)\ii)\| \ .$$ We have
\begin{equation*}
\|Z_{t,1}Q_{\Pp_1}\| = O(t^{1/2}) \ ,
\end{equation*}
\begin{equation*}
\|(\Ds_{\Pp_1} \pm t^{-1/2}\la^{1/2}\If)\ii\| =
O(|\la|^{-1/2}).t^{1/2} \hskip 10mm {\rm as} \ t\too 0 \ ,
\end{equation*}
\begin{equation*}
\|Z_{t,1}S\|  = O(t^{1/2}) \ ,
\end{equation*}
\begin{equation*}
\|(\Ds_{\Pp_1}^2 - t\ii\la\If)\ii\|  = O(t |\la|\ii) \hskip 10mm
{\rm as} \ t\too 0 \ ,
\end{equation*}
and hence that
\begin{equation}\label{e:term2}
\|(Z_{t,1}(\Ds_{\Pp_1}^2 - t\ii\la\If)\ii)\| = O(t^0 |\la|\ii)
\hskip 10mm {\rm as} \ t\too 0 \ .
\end{equation}
From \eqref{e:firstestimate}, \eqref{e:term1}, \eqref{e:term2} we
have that \eqref{e:str1} is  $O(t^0 |\la|\ii)$ as $t\too 0$. An
entirely similar analysis gives the same result for the supertrace
of \eqref{e:relexpand2}. Hence we have
\begin{equation*}
|\str((\Ff_{t,\Pp_1} - \la\If )\ii -(\Ff_{t,\Pp_2} - \la\If )\ii)|
= O(t^0 |\la|\ii)
\end{equation*}
as $t\too 0$, which by \corref{c:relch} proves the first sentence
of the Proposition.

To prove \eqref{e:simplerlimit}, we expand the resolvent
\begin{equation*}
(\Rs_i + t^{1/2}\na^i \Ds_{\Pp_i}  + t\Ds_{\Pp_i}^2 - \la\If)\ii
 = (\Rs_i  + t\Ds_{\Pp_i}^2 - \la\If)\ii
\end{equation*}
$$+ \sum_{k=1}^{\dim B}t^{k/2}
(\Rs_i + t\Ds_{\Pp_i}^2 - \la\If)\ii\left((\na^i \Ds_{\Pp_i})
(\Rs_i + t\Ds_{\Pp_i}^2 - \la\If)\ii\right)^k \ .$$ Hence
\begin{equation*} \left\{(\Ff_{t,1} - \la\If)\ii - (\Ff_{t,1} -
\la\If)\ii\right\}- \left\{(\Rs_1  + t\Ds_{\Pp_1}^2 - \la\If)\ii -
(\Rs_2  + t\Ds_{\Pp_2}^2 - \la\If)\ii\right\}
\end{equation*}
$$ =  \sum_{k=1}^{\dim B}t^{k/2}\left\{
(\Rs_1 + t\Ds_{\Pp_1}^2 - \la\If)\ii\left((\na^1\Ds_{\Pp_1})
(\Rs_1 + t\Ds_{\Pp_1}^2 - \la\If)\ii\right)^k  \right.\ .$$
\begin{equation}\label{e:difference}
 \left. - (\Rs_2 + t\Ds_{\Pp_2}^2 - \la\If)\ii\left((\na^2
\Ds_{\Pp_2}) (\Rs_2 + t\Ds_{\Pp_2}^2 - \la\If)\ii \right)^k
\right\} \ . \end{equation} Since
$$(\Rs_i + t\Ds_{\Pp_i}^2 - \la\If)\ii(\na^i \Ds_{\Pp_i})
(\Rs_i + t\Ds_{\Pp_i}^2 - \la\If)\ii  $$
$$ = (t\Ds_{\Pp_i}^2 - \la\If)\ii .\sum_{j\geq 0}(\Rs_i(t\Ds_{\Pp_i}^2 -
\la\If)\ii)^j . (\na^i \Ds_{\Pp_i}) .\sum_{j\geq
0}(\Rs_i(t\Ds_{\Pp_i}^2 - \la\If)\ii)^j $$ we obtain using
\lemref{lem:relativelysmooth} that the terms in the sum in
\eqref{e:difference} are smooth families of smoothing operators,
and by a exactly similar argument to the first part that $$
|\str\left( (\Rs_1 + t\Ds_{\Pp_1}^2 -
\la\If)\ii\left\{(\na^1\Ds_{\Pp_1}) (\Rs_1 + t\Ds_{\Pp_1}^2 -
\la\If)\ii\right\}^k  \right. $$
$$
 \left. - (\Rs_2 + t\Ds_{\Pp_2}^2 - \la\If)\ii\left\{(\na^2
\Ds_{\Pp_2}) (\Rs_2 + t\Ds_{\Pp_2}^2 - \la\If)\ii \right\}^k
\right)| = O(1)
$$
as $t\too 0$. It follows that the supertrace of
\eqref{e:difference} is $O(t^{1/2})$ as $t\too 0$, and
\eqref{e:simplerlimit} follows.

Taking sums of derivatives with respect to the base $B$ produces
similar terms and hence the assertions hold in each $C^l$ norm.
\end{proof}

It follows, then, that the object of interest is the asymptotic
behavior as $t\too 0$ of
\begin{equation}\label{e:relparametrix}
\str\left\{(\Rs_1 + t\Ds_{\Pp_1}^2 - \la\If)\ii - (\Rs_2 +
t\Ds_{\Pp_2}^2 - \la\If)\ii \right\}
\end{equation}
$$ = \sum_{k=0}^{\o} (-1)^k \str\left(\Gs_{t,\la}\ii(\Rs_1
\Gs_{t,\la}\ii)^k -  \Hs_{t,\la}\ii(\Rs_1 \Hs_{t,\la}\ii)^k
\right)$$ where
\begin{equation*}
\Gs_{t,\la}\ii = (t\Ds_{\Pp_1}^2 - \la\If)\ii \ , \hskip 15mm
\Hs_{t,\la}\ii = (t\Ds_{\Pp_2}^2 - \la\If)\ii \ .
\end{equation*}
\vskip 2mm
\begin{prop}\label{p:lim0}
As $\la\too\o$ and for all $C^l$ norms
\begin{equation}\label{e:lim0}
\str\left(\Gs_{t,\la}\ii(\Rs_1 \Gs_{t,\la}\ii)^k -
\Hs_{t,\la}\ii(\Rs_2 \Hs_{t,\la}\ii)^k  \right)
\end{equation}
\begin{equation}\label{e:lim00}
= \str(\Rs^k_1 - \Rs^k_2)(-\la)^{-k-1} + O(|\la|^{-k-2}).t
\end{equation}
\begin{equation}\label{e:lim00a}
\ + \ (k+1) \str(\Rs(\Ds)^k(\Gs_{t,\la}\ii -
\Hs_{t,\la}\ii))(-\la)^{-k} + O(|\la|^{-k-1}).t \ .
\end{equation}

\noi Here $\Rs^{0}_i := \If$. In \eqref{e:lim00a}, the curvature
form $\Rs(\Ds)$ defined by the Calderon section $P(\Ds)\in
\Gr(\Ef^0)$, can be replaced by $ \Rs(\Ds) + \Sss$, where $\Sss$
is any vertical smoothing operator in
$\Aa(B,\Psi^{-\o}(\Ef,\Ef^{'}))$, and hence by any curvature form
$\Rs^{\pi_*(\Ef|\Pp)}$  with $\Pp\in \Gr(\Ef^0)$, without
affecting the constant order $t^0$ term.
\end{prop}
\begin{proof}
For $k=0$ both sides are equal to $\str(\Gs_{t,\la}\ii -
\Hs_{t,\la}\ii)$. To explain the proof for $k>0$ we prove first
the $k=1$ case. We have
\begin{equation*}
\str\left(\Gs_{t,\la}\ii \Rs_1 \Gs_{t,\la}\ii - \Hs_{t,\la}\ii
\Rs_2 \Hs_{t,\la}\ii \right)
\end{equation*}
\begin{equation*}
=\str\left(\Gs_{t,\la}\ii(\Rs_1 - \Rs_2)\Gs_{t,\la}\ii +
(\Gs_{t,\la}\ii - \Hs_{t,\la}\ii) \Rs_2 \Gs_{t,\la}\ii +
\Hs_{t,\la}\ii \Rs_2 (\Gs_{t,\la}\ii - \Hs_{t,\la}\ii)\right)
\end{equation*}
\begin{equation*}
=\str\left(\Gs_{t,\la}\ii(\Rs_1 - \Rs_2)\Gs_{t,\la}\ii\right) +
\str\left((\Gs_{t,\la}\ii - \Hs_{t,\la}\ii) \Rs_2
\Gs_{t,\la}\ii\right) + \str\left(\Hs_{t,\la}\ii \Rs_2
(\Gs_{t,\la}\ii - \Hs_{t,\la}\ii)\right)
\end{equation*}
\begin{equation}\label{e:k1a}
=\str\left(\Gs_{t,\la}\ii(\Rs_1 - \Rs_2)\Gs_{t,\la}\ii\right) +
\str\left(\Rs_2 \Gs_{t,\la}\ii(\Gs_{t,\la}\ii - \Hs_{t,\la}\ii)
\right) + \str\left(\Hs_{t,\la}\ii \Rs_2 (\Gs_{t,\la}\ii -
\Hs_{t,\la}\ii)\right) \ ,
\end{equation}
using the norm boundedness of $\Rs_2 \Gs_{t,\la}\ii$ and
\eqref{e:smooth2} to cycle the trace.

From the resolvent formulae on $\Aa(B,\End(\pi_*(\Ef)))$
\begin{equation}\label{e:k1b}
\Gs_{t,\la}\ii = (-\la)\ii - t\Ds_1^2 \Gs_{t,\la}\ii (-\la)\ii \ ,
\hskip 10mm  \Hs_{t,\la}\ii = (-\la)\ii - t\Ds_2^2 \Hs_{t,\la}\ii
(-\la)\ii \ ,
\end{equation}
where $\Ds_i := \Ds_{\Pp_i}$, we have
\begin{equation}\label{e:k1c}
\str\left(\Gs_{t,\la}\ii(\Rs_1 - \Rs_2)\Gs_{t,\la}\ii\right)
\end{equation}
\begin{equation*}
= \str\left(\Rs_1 - \Rs_2\right) \, (-\la)^2 -
 \str\left((\Rs_1 - \Rs_2)\Ds_1^2 \Gs_{t,\la}\ii\right) \, t(-\la)^2
\end{equation*}
\begin{equation*}
 - \ \str\left(\Ds_1^2 \Gs_{t,\la}\ii(\Rs_1 - \Rs_2)\right) \, t(-\la)^2
 +  \ \str\left(\Ds_1^2 \Gs_{t,\la}\ii(\Rs_1 - \Rs_2)
 \Ds_1^2 \Gs_{t,\la}\ii\right) \, t(-\la)^2 \ .
\end{equation*}
Since $(\Rs_1 - \Rs_2)\Ds_1^2 \in\Aa(B,\Psi^{-\o}(\Ef,\Ef^{'}))$
we have
\begin{eqnarray*}
|\str((\Rs_1 - \Rs_2)\Ds_1^2 \Gs_{t,\la}\ii)| & \leq &
\|(\Rs_1 - \Rs_2)\Ds_1^2 \Gs_{t,\la}\ii\|_1 \\
 & \leq & \|(\Rs_1 - \Rs_2)\Ds_1^2\|_1 \, \|\Gs_{t,\la}\ii\| \\
& = & O(t^0 |\la|\ii)  \hskip 5mm {\rm as } \ \ \la\too \o \ ,
\end{eqnarray*}
where $\| . \|_1$ is the trace norm, and hence
\begin{eqnarray*}
|\str((\Rs_1 - \Rs_2)\Ds_1^2 \Gs_{t,\la}\ii) \, t(-\la)^2 | = O(
|\la|^{-3}).t  \hskip 5mm {\rm as } \ \ \la\too \o \ .
\end{eqnarray*}
Similar estimates hold for the third and fourth terms on the
right-side of \eqref{e:k1c}, and hence
\begin{equation}\label{e:k1d}
|\str(\Gs_{t,\la}\ii(\Rs_1 - \Rs_2)\Gs_{t,\la}\ii) - \str(\Rs_1 -
\Rs_2) \, (-\la)^2 | = O(|\la|^{-3}).t \hskip 5mm {\rm as } \ \
\la\too \o \ .
\end{equation}
For the second term in \eqref{e:k1a} we have using \eqref{e:k1b}
\begin{equation*}
\str\left(\Rs_2 \Gs_{t,\la}\ii(\Gs_{t,\la}\ii - \Hs_{t,\la}\ii)
\right)
\end{equation*}
\begin{equation}\label{e:k1e}
 = \str(\Rs_2(\Gs_{t,\la}\ii - \Hs_{t,\la}\ii)) (-\la)\ii
 - \str(\Rs_2\Ds^2_1 \Gs_{t,\la}\ii(\Gs_{t,\la}\ii -
 \Hs_{t,\la}\ii)\, t(-\la)\ii \ .
\end{equation}
The second term of \eqref{e:k1e} is equal to
\begin{equation*}
 \str(\Rs_2\Ds^2_1 (\Gs_{t,\la}\ii - \Hs_{t,\la}\ii)) \,  t(-\la)^{-2}
-  \str(\Rs_2\Ds^4_1 \Gs_{t,\la}\ii(\Gs_{t,\la}\ii -
\Hs_{t,\la}\ii)) \, t^2 (-\la)^{-2}
\end{equation*}
which is $O(|\la|^{-3}).t$ as $\la\too \o$. Thus the second and,
by a similar analysis, the third terms of \eqref{e:k1a} are equal
to
\begin{equation}\label{e:k1f}
 \str(\Rs_2(\Gs_{t,\la}\ii - \Hs_{t,\la}\ii)) (-\la)\ii
+ O(|\la|^{-3}).t \hskip 5mm {\rm as } \ \ \la\too \o \ .
\end{equation}
Equations \eqref{e:k1d} and \eqref{e:k1f} prove \eqref{e:lim0} for
the sup-norm. Taking derivatives with respect to $B$ gives similar
estimates proving \eqref{e:k1} for each $C^l$ norm.

To explain in this case the statement on the curvature form in
\eqref{e:lim00a}, let $\Sss\in\Aa(B,\Psi^{-\o}(\Ef,\Ef^{'}))$.
Then
\begin{equation}\label{e:Sss}
 \str(\Sss (\Gs_{t,\la}\ii - \Hs_{t,\la}\ii)) =
 t(-\la)\ii \,(\Sss\Ds^2_2\Hs_{t,\la}\ii -
 \Sss\Ds^2_1\Gs_{t,\la}\ii)\ ,
\end{equation}
and hence
\begin{equation*}
 |\str(\Sss (\Gs_{t,\la}\ii - \Hs_{t,\la}\ii))|  \leq
 t(-\la)\ii\,(\|\Sss\Ds^2_2\|_1\,\|\Hs_{t,\la}\ii\| +
 \|\Sss\Ds^2_1\|_1\,\|\Gs_{t,\la}\ii\|)
  =  O( |\la|^{-2}).t \ .
\end{equation*}
This proves that
\begin{equation}\label{e:k1}
\str\left(\Gs_{t,\la}\ii \Rs_1 \Gs_{t,\la}\ii -  \Hs_{t,\la}\ii
\Rs_2 \Hs_{t,\la}\ii \right) = \str(\Rs_1 - \Rs_2)(-\la)^{-2} +
O(|\la|^{-3}).t
\end{equation}
$$ + 2\,\str(\Rs(\Ds)(\Gs_{t,\la}\ii - \Hs_{t,\la}\ii))(-\la)^{-2}
+ O(t|\la|^{-2}) \ .
$$

\vskip 3mm

To prove \eqref{e:lim0} for the general case, we have inductively
for $k\geq 1$, by \eqref{e:smooth2} and \eqref{e:smooth3},
\begin{equation}\label{e:lim0a}
\str\left(\Gs_{t,\la}\ii(\Rs_1 \Gs_{t,\la}\ii)^k -
\Hs_{t,\la}\ii(\Rs_2 \Hs_{t,\la}\ii)^k  \right) =
\end{equation}
\begin{equation*}
\sum_{i=0}^k \str\left((\Hs_{t,\la}\ii \Rs_2)^i (\Gs_{t,\la}\ii -
\Hs_{t,\la}\ii)(\Rs_1 \Gs_{t,\la}\ii)^{k-i} \right)
\end{equation*}
$$+$$
\begin{equation*}
 \sum_{j=0}^{k-1} \str\left((\Hs_{t,\la}\ii \Rs_2)^j
\Hs_{t,\la}\ii (\Rs_1 - \Rs_2)\Gs_{t,\la}\ii(\Rs_1
\Gs_{t,\la}\ii)^{k-j-1} \right) \ .
\end{equation*}
For the second sum in \eqref{e:lim0a}
\begin{equation*}
\str\left((\Hs_{t,\la}\ii \Rs_2)^j \Hs_{t,\la}\ii (\Rs_1 -
\Rs_2)\Gs_{t,\la}\ii(\Rs_1 \Gs_{t,\la}\ii)^{k-j-1} \right)
\end{equation*}
\begin{equation*}
= \str\left(((-\la)^{-j}(\Rs_2 - t\Ds^2_2 \Hs_{t,\la}\ii \Rs_2
)^j). \, (-\la)\ii(\If - t\Ds^2_1 \Hs_{t,\la}\ii).(\Rs_1 -
\Rs_2)\right.
\end{equation*}
\begin{equation*}
\left. \times \ (-\la)\ii (\If - t\Ds^2_1
\Gs_{t,\la}\ii).(-\la)^{-k+j+1}
 (\Rs_1 - t\Rs_1\Ds^2_1 \Gs_{t,\la})^{k-j-1} \right)
\end{equation*}
\begin{equation*}
= (-\la)^{-k-1}\str\left((\Rs_2 - t\Ds^2_2 \Hs_{t,\la}\ii \Rs_2
)^j(\If - t\Ds^2_1 \Hs_{t,\la}\ii)(\Rs_1 - \Rs_2)
  (\If - t\Ds^2_1 \Gs_{t,\la}\ii)
 (\Rs_1 - t\Rs_1\Ds^2_1 \Gs_{t,\la})^{k-j-1} \right)
\end{equation*}
\begin{equation*}
= (-\la)^{-k-1}\str\left((\Rs_2 - t\Ds^2_2 \Hs_{t,\la}\ii \Rs_2
)^j((\Rs_1 - \Rs_2)- t\Ds^2_2 \Hs_{t,\la}\ii(\Rs_1 - \Rs_2)\right.
 \end{equation*}
\begin{equation*}
\left.- t (\Rs_1 - \Rs_2)\Ds^2_1 \Gs_{t,\la}\ii
 + t^2 \Ds_2^2 \Hs_{t,\la}\ii(\Rs_1 - \Rs_2)\Ds^2_1
 \Gs_{t,\la}\ii
 (\Rs_1 - t\Rs_1\Ds^2_1 \Gs_{t,\la})^{k-j-1} \right) \ .
\end{equation*}
\begin{center}
----expanding the terms $(\Rs_2 - t\Ds^2_2
\Hs_{t,\la}\ii \Rs_2 )^j$ and $(\Rs_1 - t\Rs_1\Ds^2_1
\Gs_{t,\la})^{k-j-1}$----
\end{center}
\begin{equation}\label{e:lim0b}
= \str\left(\Rs_2^j \Rs_1^{k-j} - \Rs_2^{j+1}\Rs_1^{k-j-1}\right)
(-\la)^{-k-1} + \str({\rm
pol}(t,\Rs_1,\Rs_2,\Gs_{t,\la}\ii,\Hs_{t,\la}\ii)) \ ,
\end{equation}
where ${\rm pol}(t,\Rs_1,\Rs_2,\Gs_{t,\la}\ii,\Hs_{t,\la}\ii)$ is
a finite polynomial expression consisting of terms like
\begin{equation*}
{\rm const.} \, t^{j-p}\Rs_2^p (\Ds_2^2 \Hs_{t,\la}\ii
\Rs_2)^{j-p}. t(\Rs_1 - \Rs_2)\Rs_1 \Gs_{t,\la}\ii . \Rs_1^q
(\Rs_1 \Ds^2_1 \Gs_{t,\la}\ii)^{k-j-1-q} \ ,
\end{equation*}
where $0\leq p \leq j-1$, $0\leq q \leq k-j-2$. By an estimate
akin to the $k=1$ case, we obtain
\begin{equation}\label{e:lim0c}
|\str({\rm pol}(t,\Rs_1,\Rs_2,\Gs_{t,\la}\ii,\Hs_{t,\la}\ii))| =
O(|\la|^{-k-2}).t \hskip 5mm {\rm as } \ \ \la\too \o \ .
\end{equation}
From \eqref{e:lim0b}, \eqref{e:lim0c}, the second sum in
\eqref{e:lim0a} hence telescopes down to
\begin{equation}
 \sum_{j=0}^{k-1} \str\left(\Rs_2^j \Rs_1^{k-j} -
 \Rs_2^{j+1}\Rs_1^{k-j-1}\right)(-\la)^{-k-1}
  + O(|\la|^{-k-2}).t
\end{equation}
\begin{equation*}
 = \str\left(\Rs_1^k -
 \Rs_2^k\right)(-\la)^{-k-1}
  + O(|\la|^{-k-2}).t \ ,
\end{equation*}
which is \eqref{e:lim00}.

For the first sum in \eqref{e:lim0a}, since $(\Rs_1
\Gs_{t,\la}\ii)^{k-i}$ is a vertical family of bounded operators,
we have cycling the trace
\begin{equation*}
  \str\left((\Hs_{t,\la}\ii \Rs_2)^i (\Gs_{t,\la}\ii -
\Hs_{t,\la}\ii)(\Rs_1 \Gs_{t,\la}\ii)^{k-i} \right)
\end{equation*}
\begin{equation*}
  = \str\left((\Rs_1 \Gs_{t,\la}\ii)^{k-i}
  (\Hs_{t,\la}\ii \Rs_2)^i (\Gs_{t,\la}\ii -
\Hs_{t,\la}\ii)\right)
\end{equation*}
\begin{equation*}
 = (-\la)^{-k}\str\left((\Rs_1 - t\Rs_1\Ds^2_1 \Gs_{t,\la}\ii)^{k-i}
 (\Rs_2 - t\Ds^2_2 \Hs_{t,\la}\ii  \Rs_2)^{i}
 (\Gs_{t,\la}\ii - \Hs_{t,\la}\ii) \right)
\end{equation*}
\begin{equation*}
 = (-\la)^{-k}\str\left(\Rs_1^{k-i}\Rs_2^{i}
 (\Gs_{t,\la}\ii - \Hs_{t,\la}\ii)
 \right) + O(|\la|^{-k-1}).t
\end{equation*}
\begin{equation*}
 = (-\la)^{-k}\str\left(\Rs(\Ds)^k(\Gs_{t,\la}\ii -
 \Hs_{t,\la}\ii)
 \right) + O(|\la|^{-k-1}).t \ ,
\end{equation*}
for since $\Rs_1^{k-i}\Rs_2^{i} - \Rs(\Ds)^k
\in\Aa(B,\Psi^{-\o}(\Ef,\Ef^{'}))$ we can apply the argument of
\eqref{e:Sss} to make the replacement.

Since there are $k+1$ such terms in the first sum of
\eqref{e:lim0a}, we get \eqref{e:lim00a}.
\end{proof}

\begin{rem}
Essentially the point here is that if $Q_{\la}\ii = (\Delta
-\la)\ii$ is a resolvent for a differential operator $\Delta$, and
$\Sss$ a smoothing operator, there is an asymptotic expansion
\begin{equation}\label{e:asympsmooth}
\str(\Sss Q_{\la}\ii) \sim \str(\Sss)(-\la)\ii + \sum_{j\geq
1}\str(\Sss\Delta^j)(-\la)^{-j-1} \ ,
\end{equation}
as $\la\to\o$, which if $\Delta$ is a positive operator has the
heat trace realization as $t\to 0$
$$\str(\Sss e^{-t\Delta}) \sim \str(\Sss) + \sum_{j\geq
1}\frac{1}{j!}\str(\Sss\Delta^j) t^j \ ,$$ which follows from
\eqref{e:asympsmooth} via
$$\str(\Sss e^{-t\Delta}) = \frac{i}{2\pi}\int_{\Cc}e^{-\la}\str(\Sss
Q_{\la}\ii)\ d\la \ .$$
\end{rem}

We need the following result in order to deal with the term
\eqref{e:lim00a} in the Chern character. For convenience write
$$\Rs := \Rs(\Ds) \ , \hskip 10mm \Gs_{\la}\ii = \Gs_{1,\la}\ii\ ,
\hskip 10mm \Hs_{\la}\ii = \Hs_{1,\la}\ii \ . $$

\begin{prop}
As $\la\too \o$ there is an asymptotic expansion of the
differential form $\str(\Rs^k(\Gs_{\la}\ii - \Hs_{\la}\ii))$ on
$B$
\begin{equation}\label{e:expansioncanontrace}
\str(\Rs^k(\Gs_{\la}\ii - \Hs_{\la}\ii))(z) \sim \sum_{j\geq 0}
c_{j,k} (z)(-\la)^{-\frac{j}{2}-1} + \sum_{i\geq 1} d_{i,k}
(z)\log(-\la).(-\la)^{-\frac{i}{2}-1} \ .
\end{equation}
\end{prop}
\begin{proof}
For $k=0$ this is the asymptotic expansion as $\la\too\o$ of
\cite{Gr99,Gr99a}
\begin{equation}\label{e:k0}
\str((\Ds_{\Pp_1}^2-\la\If)\ii - (\Ds_{\Pp_2}^2-\la\If)) \sim
\sum_{j=0}^{\o} a_j (-\la)^{-1-j/2} \ ,
\end{equation}
with $a_j\in\G(B)$. The method of \cite{Gr99,Gr99a} applies also
to the case here. (The existence of a trace expansion for
$\str(\Rs^k\Gs_{t,\la}\ii)$ is a harder problem.) We refer to
\cite{Gr96} for a comprehensive account of the symbol spaces
below.

Write $\la = \mu^2, \mu\in\Cf\backslash \Rf$. From
\eqref{e:smooth6} or \cite{Gr99a} we obtain
\begin{equation}\label{e:smooth7}
\Rs^k(\Gs_{\mu^2}\ii -  \Hs_{\mu^2}\ii) = \mu\ii \Rs^k\Ks_{\mu}\g
\Ss_{\mu}\Ts_{\mu} \ ,
\end{equation}
where in the sgo calculus of \cite{Gr96,Gr99,Gr99a} $\Ks_{\mu}$ is
a vertical strongly polyhomogeneous Poisson operator of order $0$,
$\Ts_{\mu} $ is a vertical trace operator of order $-1$, and
$\Ss_{\mu}$ is a vertical weakly polyhomogeneous $\pdo$ acting on
$\Aa(B,\pi_*(\Ef^{'}))$ with symbol pointwise in $S^{-\o,0}$.
Since $F_{k,\mu} = \Rs^k(\Gs_{\la}\ii - \Hs_{\la}\ii)$   is the
combination of a vertical differential-operator of order $k$
valued form, vertical Poisson operators and singular Green's
operators, and since the relative inverse is a smooth family of
smoothing operators \eqref{e:smooth2}, then by the symbol calculus
and \lemref{lem:commutatortrace}[1]
$F_{k,\mu}\in\Aa(S,\Psi^{-\o}(\Ef,\Ef^{'}))$ with cycled trace
equal to the trace $\mu\ii \str(\Ss_{\mu}\Ts_{\mu}
\Rs^k\Ks_{\mu}\g)$ on $\Aa(B,\Psi^{-\o}(\Ef^{'}))$. Since
$\Ts_{\mu} \Rs^k\Ks_{\mu}$ is a strongly polyhomogeneous vertical
pdo  in $\Aa(B,\Psi^{k-1}(\Ef^{'}))$ of order $k-1$ and with
symbol in $S^{k-1,0} \cap S^{0,k-1}$, the composed operator is a
weakly polyhomogeneous vertical pdo with symbol pointwise in
$S^{-\o,k-1} \cap S^{-\o,k-2}$. Now by \cite{GrSe95} there is an
asymptotic expansion
\begin{equation*}
\str(F_{k,\mu})(z) \sim \sum_{-\o < j < 0} \alpha_{j,k}
(z)\mu^{k-2-j} + \sum_{j\geq 0} (\beta_{j,k} (z)\log \mu +
\beta^{'}_{j,k}  (z))\mu^{k-2-j} \ .
\end{equation*}
Next since $F_{k,\mu}\in\Aa(S,\Psi^{-\o}(\Ef,\Ef^{'}))$  we have
$|\str(F_{k,\mu})(z)| = O(|\mu^{-2}|)$ and hence all terms
$\a_{j,k} ,\b_{j,k} $ with $j >k$ and terms $\b^{'}_{j,k} $ with
$j> k+1$ vanish.  Finally, since the log coefficients come
entirely from the homogeneous symbol of $\Ss_{\mu}$ one has as in
\cite{GrSe95} Thm.2.1 that the terms with, respectively,  $j = k$
and $j = k+1$ vanish. This completes the proof.
\end{proof}

Let $\Delta_i = D_{\Pp_i}^2$. Then the operator
$\Rs^k(\Delta_1^{-z} - \Delta_2^{-z})$ is trace class for
$\re(z)>>0$. The relative superzeta function trace
$\str(\Rs^k(\Delta_1^{-z} - \Delta_2^{-z}))$ is holomorphic for
such $z$ and has a meromorphic continuation to all of $\Cf$ with
singularity structure determined by the asymptotic expansion
\eqref{e:expansioncanontrace}. Since there is no term
$(-\la)\ii\log(-\la)$ there is a well-defined pseudo-trace
functional $\tau_{\Delta_1,\Delta_2}(\Rs^k)$, as explained in
Section(1.2).

\begin{thm}\label{t:limitatzero} As forms in $\Aa^{2\bullet}(B)$
\begin{equation}\label{e:limitatzero}
\lim_{t\too 0} \left(\ch(\Af_{t,\Pp_1}) -
\ch(\Af_{t,\Pp_2})\right) = \eta^{[M]}(\Ps_1,\Ps_2) + \sum_{k\geq
0}\frac{k+1}{k!} \tau_{\Delta_1,\Delta_2}(\Rs^k) \ .
\end{equation}
Here $\Rs^0_i := \If$ and the degree zero part of the right-side
of \eqref{e:limitatzero} is the pointwise relative index
$\ind(D_{P_1}) - \ind(D_{P_2})$.
\end{thm}
\begin{proof}
From \eqref{e:relparametrix}, \eqref{e:simplerlimit},
\eqref{e:cherncharacter2} and \propref{p:lim0} we have as $t\to 0$
\begin{equation}\label{e:relch}
\ch(\Af_{t,\Pp_1}) - \ch(\Af_{t,\Pp_2})
\end{equation}
\begin{equation*}
= \sum_{k=0}^{\o} (-1)^k \frac{i}{2\pi}\int_{\Cc} e^{-\la}
\str\left(\Gs_{t,\la}\ii(\Rs_1 \Gs_{t,\la}\ii)^k -
\Hs_{t,\la}\ii(\Rs_1 \Hs_{t,\la}\ii)^k \right) \ d\la  + o(1)
\end{equation*}
\begin{equation}\label{e:sum1}
= \sum_{k=1}^{\o} (-1)^k \frac{i}{2\pi}\int_{\Cc} e^{-\la}
(-\la)^{-k}  \ d\la  \,. \, \str(\Rs^k_1 - \Rs^k_2)
\end{equation}
\begin{equation}\label{e:sum2}
+ \sum_{k=0}^{\o} (-1)^k (k+1)\frac{i}{2\pi}\int_{\Cc} e^{-\la}
(-\la)^{-k} \str(\Rs^k(\Gs_{t,\la}\ii - \Hs_{t,\la}\ii)) \ d\la
\,. \, + o(1) \ .
\end{equation}
(Specifically the $o(1)$ term is given by an asymptotic expansion
$\sum_{j\geq 1} a_k t^{k/2}$ as $t\too 0$.) The first sum
\eqref{e:sum1} determines the $\eta$-form part of
\eqref{e:limitatzero}.  The zeta trace terms come from
\eqref{e:sum2}, as follows. From \eqref{e:expansioncanontrace} and
\begin{equation}\label{e:t}
\str(\Rs^k(\Gs_{t,\la}\ii - \Hs_{t,\la}\ii)) = t\ii
\str(\Rs^k(\Gs_{t\ii\la}\ii - \Hs_{t\ii\la}\ii))
\end{equation}
we obtain an asymptotic expansion as $t\to 0$
\begin{equation*}
\str(\Rs^k((\Gs_{t,\la}\ii - \Hs_{t,\la}\ii)) \sim
\frac{1}{t}\sum_{j\geq 0} c_{j,k}
\left(-\frac{\la}{t}\right)^{-\frac{j}{2}-1} +
\frac{1}{t}\sum_{i\geq 1} d_{i,k}
\log\left(-\frac{\la}{t}\right).\left(-\frac{\la}{t}\right)^{-\frac{i}{2}-1}
\ .
\end{equation*}
Hence for any $\e>0$ the $k^{{\rm th}}$ term of \eqref{e:sum2} is
\begin{equation*}
(k+1)\frac{i}{2\pi}\int_{\Cc} e^{-\la} (-\la)^{-k}
\str(\Rs^k((\Gs_{t,\la}\ii - \Hs_{t,\la}\ii)) \ d\la
\end{equation*}
\begin{equation*}
= (k+1) \sum_{j=0}^{N-1} c_{j,k}  \left(\frac{i}{2\pi}\int_{\Cc}
e^{-\la} (-\la)^{-\frac{j}{2}-k-1} \ d\la \right) \,
t^{\frac{j}{2}} +  O(t^{N/2})
\end{equation*}
\begin{equation*}
+ (k+1) \sum_{j=1}^{N-1} d_{j,k}  \left(\frac{i}{2\pi}\int_{\Cc}
e^{-\la} (-\la)^{-\frac{j}{2}-k-1} \log(-\frac{\la}{t}) \ d\la
\right) \, t^{\frac{j}{2}}\log t + O(t^{N/2 - \e})
\end{equation*}
\begin{equation*}
= \frac{(k+1)}{\G(k+1)}c_{0,k}  + (k+1)\sum_{j=1}^{N-1}
\G(\frac{j}{2}+k+1)\ii c_{j,k} \, t^{\frac{j}{2}}
\end{equation*}
\begin{equation*}
+ (k+1)\sum_{j=1}^{N-1} \G^{'}(\frac{j}{2}+k+1)\ii d_{j,k} \,
t^{\frac{j}{2}} - (k+1)\sum_{j=1}^{N-1} \G(\frac{j}{2}+k+1)\ii
d_{j,k} \, t^{\frac{j}{2}}\log t  + o(1)
\end{equation*}
\begin{equation*}
= \frac{(k+1)}{\G(k+1)}c_{0,k} + o(1)
\end{equation*}
\begin{equation*}
= \frac{(k+1)}{\G(k+1)} \tau_{\Delta_1,\Delta_2}(\Rs^k)
 + o(1) \ .
\end{equation*}
For the degree $0$ pointwise index, observe that when $k=0$
\begin{eqnarray*}
c_{0,0} & =  & \lim_{t\too 0} \frac{i}{2\pi} \int_{\Cc} e^{-\la}
\str(\Gs_{t,\la}\ii - \Hs_{t,\la}\ii) \, d\la \\
& = & \lim_{t\too 0} \left(\str(e^{-t\Delta_1}) -
\str(e^{-t\Delta_2})\right) \\
& = & \ind(D_{P_1}) - \ind(D_{P_2}) \ .
\end{eqnarray*}
This completes the proof.
\end{proof}

\vskip 3mm


\noi The proof of Theorem (II) is completed by the following
identification.

\begin{prop}
The cohomology class of the Chern character form
$\ch(\Af_{t,\Pp})$ is equal to $\ch(\Ind(\Ds_{\Pp}))$.
\end{prop}
\begin{proof} The proof follows closely \cite{BeGeVe92} and
so we shall make only brief comments. We assume the kernels
$\Ker(D^z_{P_z})$ form a super vector bundle $\Ker(\Ds_{\Pp})$
over $B$, identified with a smooth family of finite-rank
projections $\Pi_0$ -- the general case follows by standard
arguments from this one. We then have $\ch(\na^0) =
\ch(\Ind(\Ds_{\Pp}))$, where $\na^0 =
\Pi_0\cdot\na^{\pi_*(\Ef)}\cdot\Pi_0$. By the Diagonalization
Lemma of \cite{BeGeVe92} there is an invertible
$g\in\Aa(B,\End(\pi_{*}(\Ef)))$ with $g- \If \in \sum_{r\geq
1}\Aa^r(B,\Psi^{-\o}(\Ef,\Ef^{'}))$ such that with respect to the
decomposition $\Ker(\Ds_{\Pp})\oplus\Ker(\Ds_{\Pp})\pp$
\begin{equation}\label{e:infty1}
(\Ff_{t,\Pp} - \la\If )\ii = \d_t(g)\ii\left(
\begin{bmatrix}
  t\d_t(\Uu) & 0 \\
  0 & t\d_t(\Vv) \\
\end{bmatrix} - \la\If\right)\ii\d_t(g) \ ,
\end{equation}
with
\begin{eqnarray}
\d_t(g)^{\pm 1} & = & \If + O(t^{-1/2}) \ , \label{e:infty2} \\
t\d_t(\Uu) & = & (\na^0)^2 + O(t^{-1/2}) \ , \label{e:infty3}
\\
t\d_t(\Vv) & = & t\Ds_{\Pp}^2 \ \ \mod(\sum_{r\geq
1}\Aa^r(B,\End(\pi_*(\Ef^{'})))) \ . \label{e:infty4}
\end{eqnarray}
It follows from \eqref{e:infty1},
\eqref{e:infty3},\eqref{e:infty4} that $e^{-\Ff_{t,\Pp}}$ is of
the form
$$\d_t(g)\ii(\frac{i}{2\pi}\int_{\Cc}e^{-\la}
\begin{bmatrix}
  \Pi_0((\na^0)^2 - \la\If)\ii \Pi_0 & 0 \\
  0 & \Pi_0\pp(t\Ds_{\Pp}^2 - \la\If)\ii \Pi_0\pp \\
\end{bmatrix}  \ d\la)
\begin{bmatrix}
  \If + O(t^{-1/2}) & 0 \\
  0 & O(t^{1/2}) \\
\end{bmatrix} \d_t(g)$$ and hence from \eqref{e:infty2} that
$$\str(e^{-\Ff_{t,\Pp}}) = \str(e^{-(\na^0)^2}) +
O(t^{-1/2}) \ .$$ This implies that $\lim_{t\too
\o}\ch(\Af_{t,\Pp}) = \ch(\na^0)$, which completes the proof.
\end{proof}

\vskip 1mm

For a general smooth family of elliptic operators there is an
alternative proof \cite{Sc02} that $\ch(\Af_t)$ is a
representative for the Chern character of the index bundle by
showing that the {\em zeta-form}
\begin{equation*}
\zeta(\Ff, s) := \str(\Ff\si)|^{{\rm mer}} \ ,
\end{equation*}
defined using the resolvent analysis approach, is exact at $s=0$.
The conclusion then follows from the families analogue of
\eqref{e:Deltapoles}, generalizing Seeley's zeta function formula
for the index of a single elliptic operator.

\end{document}